\documentclass[article,ij4uq]{ij4uq}              

\usepackage[hang]{footmisc}  
\fancypagestyle{plain}{%
  \fancyhf{}
  \fancyfoot[R]{\small\bf\thepage }
  \fancyfoot[L]{\fottitle}
  }

\usepackage{algorithmicx}

\def\J{\mathcal{J}}
\def\R{\mathbb{R}}

\def\E{\mathrm{E}}

\def\u{u}
\def\h{\mathsf{h}}

\renewcommand{\L}{\mathcal{L}}%
\renewcommand{\P}{\mathcal{P}}

\def\wbu{\widetilde{\textbf{u}}}

\def\bu{\textbf{u}}
\def\xik{\xi^{(k)}} 
\def\wk{w^{(k)}} 

\def\bQ{\mathbf{Q}}
\def\bA{\mathbf{A}}
\def\bi{\mathbf{i}}
\def\bff{\mathbf{f}}
\def\bu{\mathbf{u}}
\def\bS{\mathbf{S}}
\def\bq{\mathbf{q}}
\def\bc{\mathbf{c}}
\def\su{{u}}
\def\wJ{\widetilde{\mathcal{J}}}

\renewcommand{\myyear}{2017}
\renewcommand{\today}{}
\begin{document}

\volume{Volume x, Issue x, \myyear\today}
\title{An adaptive reduced basis collocation method based on PCM ANOVA decomposition for anisotropic stochastic PDEs}
\titlehead{Adaptive reduced basis collocation method for anisotropic stochastic PDEs}
\authorhead{H. Cho, \& H. Elman}
\author[1]{Heyrim Cho}
\author[2]{Howard Elman}
\corremail{hcho1237@math.umd.edu}
\address[1]{Department of Mathematics, University of Maryland, College Park, MD, 20742, USA}
\address[2]{Department of Computer Science, University of Maryland, College Park, MD, 20742, USA}


\dataO{mm/dd/yyyy}
\dataF{mm/dd/yyyy}

\abstract{

The combination of reduced basis and collocation methods 
enables efficient and accurate 
evaluation of the solutions to parameterized PDEs. 
In this paper, we study the stochastic collocation methods 
that can be 
combined with reduced basis methods to solve 
high-dimensional parameterized stochastic PDEs. 
We also propose an adaptive algorithm using a  
probabilistic collocation method (PCM) and ANOVA decomposition. 
This procedure involves two stages. 
First, the method employs an ANOVA decomposition to identify 
the effective dimensions, 
i.e., subspaces of the parameter space in which the contributions 
to the solution are larger, 
and sort the reduced basis solution 
in a descending order of error. 
Then, the adaptive search refines the parametric space 
by increasing the order of polynomials 
until the algorithm is terminated by a saturation constraint. 
We demonstrate the effectiveness of the proposed algorithm 
for solving 
a stationary stochastic convection-diffusion equation, 
a benchmark problem chosen because solutions contain 
steep boundary layers and anisotropic features. 
We show that
two stages of adaptivity are critical 
in a benchmark problem with anisotropic stochasticity. 

} 

\keywords{Reduced basis method, ANOVA decomposition, Probabilistic collocation method, Anisotropic stochasticity, High-dimensionality 
}

\maketitle

\section{Introduction} 

Numerical methods for high-dimensional anisotropic 
parameterized partial differential equations (PDEs) 
have received increasing interest in 
various areas of science and engineering. 
These problems are costly to solve because of requirements of 
high fidelity in space or high dimensionality in the parameterized 
random space. This study concerns some new approaches for reducing
costs, by combining reduced basis methods for the spatial 
component of the problem with specialized 
quadrature rules for the parameterized component. 
%
%

An approach to reduce spatial complexity is 
reduced basis methods, 
developed to cope with 
parameterized systems that require 
repeated use of a high-fidelity solver. 
Costly computations for high-resolution solutions 
for multiple parameter values 
are replaced by 
an approach with lower complexity, 
which computes parameterized solutions with 
a cost that is independent of 
the dimension of the original problem 
\cite{RBM_Ito,RBM_Peterson,RBM_IVP,RBM_Quatrni1,RBM_Rozza1,RBM_Patera1,RBM_JFE,RBM_YMaday1,RBM_YMaday2}. 
Instead of solving the full discrete PDE system, 
a cheaper system is obtained by a projection 
into a lower-dimensional parametric space determined by the 
full-system solution at a distinguished set of parameters. 
Convergence analysis is available 
for various types of PDEs 
\cite{RBM_Patera1,Quatrni_RBMbook,RBM_Rozza1} 
using so-called certified reduced basis methods. 
Extension to parameterized stochastic PDEs 
is straightforward \cite{PChen0,PChen1}, and 
applications in stochastic control \cite{PChen2} 
have been studied as well.

The two components of reduced basis methods are referred to 
as an offline step (construction of the reduced space 
by computation of required full-system solutions) 
and an online step (computations on the reduced space) 
designed to be cheap. 
In addition to requirements of spatial resolution, 
the other contribution to cost comes from 
the high dimensionality of the stochastic space. 
In order to parametrize noisy input data for a given PDE, 
random fields can often be expanded as an infinite combination 
of random variables by, for instance, 
Karhunen-Lo\`{e}ve \cite{GhanemSpanos,Papoulis} or 
polynomial chaos \cite{db1,db2} expansions. 
Such random fields can be well approximated 
using a finite number of variables \cite{Schwab12}, 
but the dimensionality may be high. 
For such high dimensions in the stochastic space, 
high-dimensional integration techniques 
\cite{Sparse_Novak,Gri05},  
based for example on 
Smolyak sparse-grid \cite{ZabarasSG,NobileTW_SINA}
or ANOVA decomposition \cite{Schwab09,CaoCG09,Foo1},  
can be employed 
to compute the statistical moments of the solution. 
Although these methods are often efficient 
compared to Monte-Carlo methods, 
the curse of dimensionality remains a challenge. 
Thus, new ideas are necessary 
to increase efficiency  
and facilitate rapid evaluation of quantities of interest.

The idea of reduced basis methods combined with 
high-dimensional collocation schemes \cite{QifangRBM} 
brings the aforementioned efforts together. 
Instead of searching arbitrary samples, 
the candidate reduced basis points are chosen to be 
the sparse-grid probabilistic collocation points. 
If the solution vector 
is not of full rank in the parametric space, 
the computational cost can be reduced 
with essentially no loss of accuracy.  
However, the online stage involves a cost that is proportional 
to the number of collocation points. Thus, the algorithm 
suffers from the curse of dimensionality with costs that 
are less than for using the full collocation scheme, 
but still depending on the stochastic dimensionality. 
This cost can be further reduced by using 
dimension adaptation approaches, 
for instance, anisotropic sparse-grid \cite{PChen1} and 
ANOVA decomposition \cite{JanH_RBMANOVA,QifangRBMA}.

In this paper, we study the reduced basis collocation method 
for stochastic PDEs, 
and propose an adaptive procedure both 
in the dimensions of the stochastic space and 
the degrees of freedom in each parametric direction. 
An outline of the paper is as follows. 
In section 2, 
we review the reduced basis collocation method based on 
ANOVA collocation methods 
and examine the quality of different interpolation points. 
Then, we introduce 
the adaptive ANOVA procedure employed in the reduced basis algorithm, 
emphasizing its advantage for reducing computational cost and 
automatic sorting. 
%
Then, in section 3, 
we propose an adaptive reduced basis algorithm 
based on a probabilistic collocation method (PCM). 
In addition to the adaptive dimension procedure described 
in section 2, 
resolution in the parametric space is managed by adapting 
the polynomial order $p$ discussed in section 3. 
We test these ideas on 
the stationary stochastic convection-diffusion equation chosen 
as a benchmark problem because of its interesting features 
such as solutions with steep boundary layers and discontinuities. 
The reduced basis ANOVA approach is suitable to study 
the interaction of these characteristics with the 
underlying noise that yields 
anisotropy in the stochastic space.  

\section{Reduced basis collocation methods based on ANOVA decomposition} 
\label{sec:RBM-SG} 
In this section, we review the reduced basis approximation for 
the solution computed using collocation in the stochastic space. 
Consider a stochastic problem consisting of 
a partial differential operator $\L$ and a boundary operator $b$, 
\begin{eqnarray} \label{eq:SPDEorg} 
	\L( x, \omega, u(x,\omega) ) = f(x), \quad ^\forall x \in D, \\ \nonumber 
	b( x, \omega, u(x,\omega) ) = g(x), \quad ^\forall x \in \partial D,	 
\end{eqnarray} 
on $D \times \Omega$, where $D \subset \R^d$ is the spacial domain, 
assumed to be a bounded and connected domain with polygonal boundary 
$\partial D$, 
and $\Omega$ is a complete probability space of 
$(\Omega, \Sigma, \mathcal{P})$ 
with $\sigma$-algebra $\Sigma$ and probability measure $\mathcal{P}$. 
Assume that the stochastic excitation can be represented 
by a finite number of random variables $\xi = [\xi_1,...,\xi_{M}]^T$ 
on $\Gamma = \prod_{m=1}^M \Gamma_m$, 
where $\Gamma_m = [a_m,\, b_m]$ is the image of $\xi_m(\omega)$. 
This could come from a variety of sources, for example, 
a truncated Karhunen-Lo\`eve (KL) expansion \cite{GhanemSpanos,Papoulis}, 
a partitioning of $D$ into subdomains, 
or uncertain boundary conditions. 

\begin{figure}
    \centerline{  \hspace{0.3cm} $1\times 16 $ \hspace{3.0cm} $4\times 4 $ \hspace{3.0cm} $16\times 1 $ } 
    \centerline{
    \includegraphics[width=4.0cm]{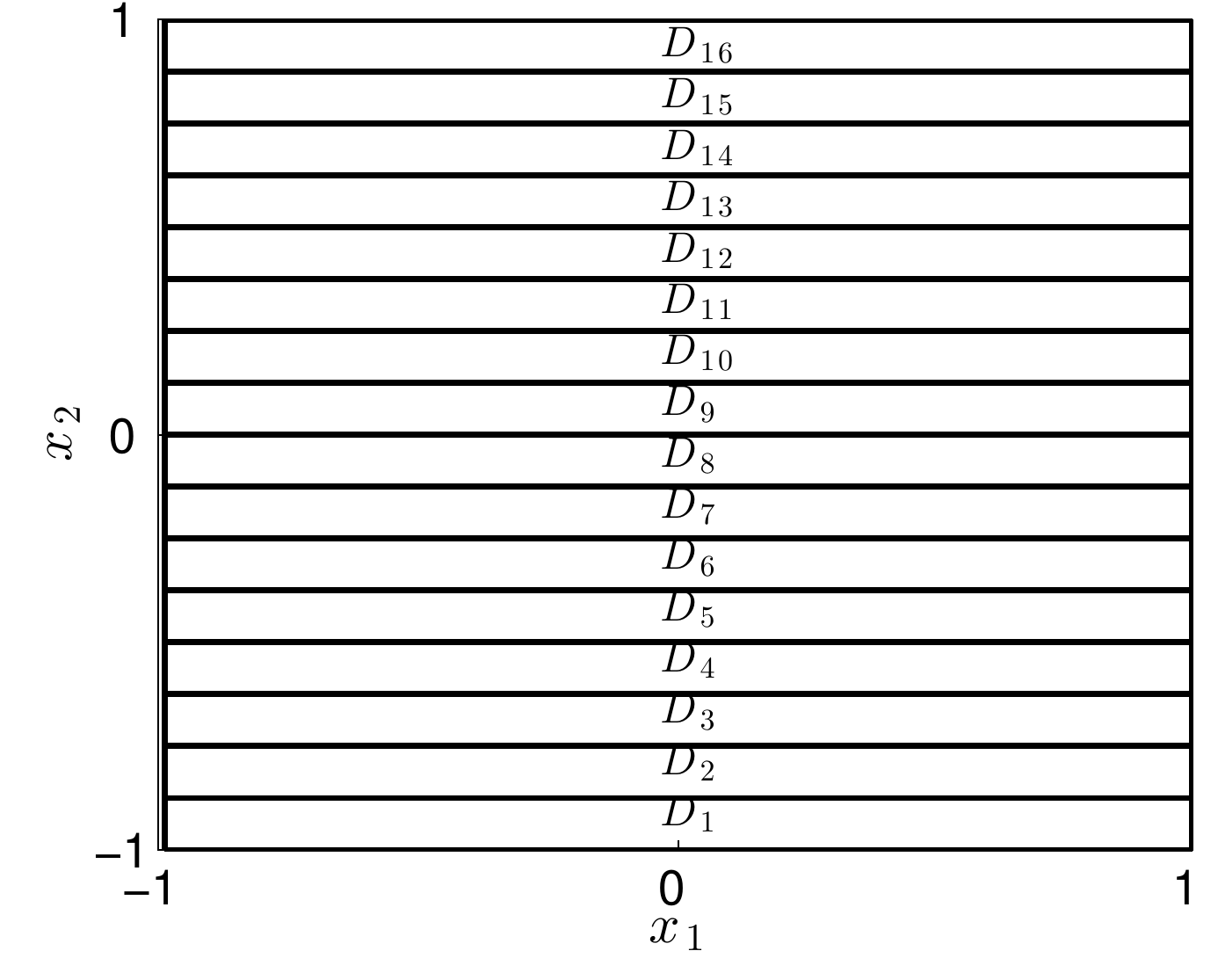}
    \includegraphics[width=4.0cm]{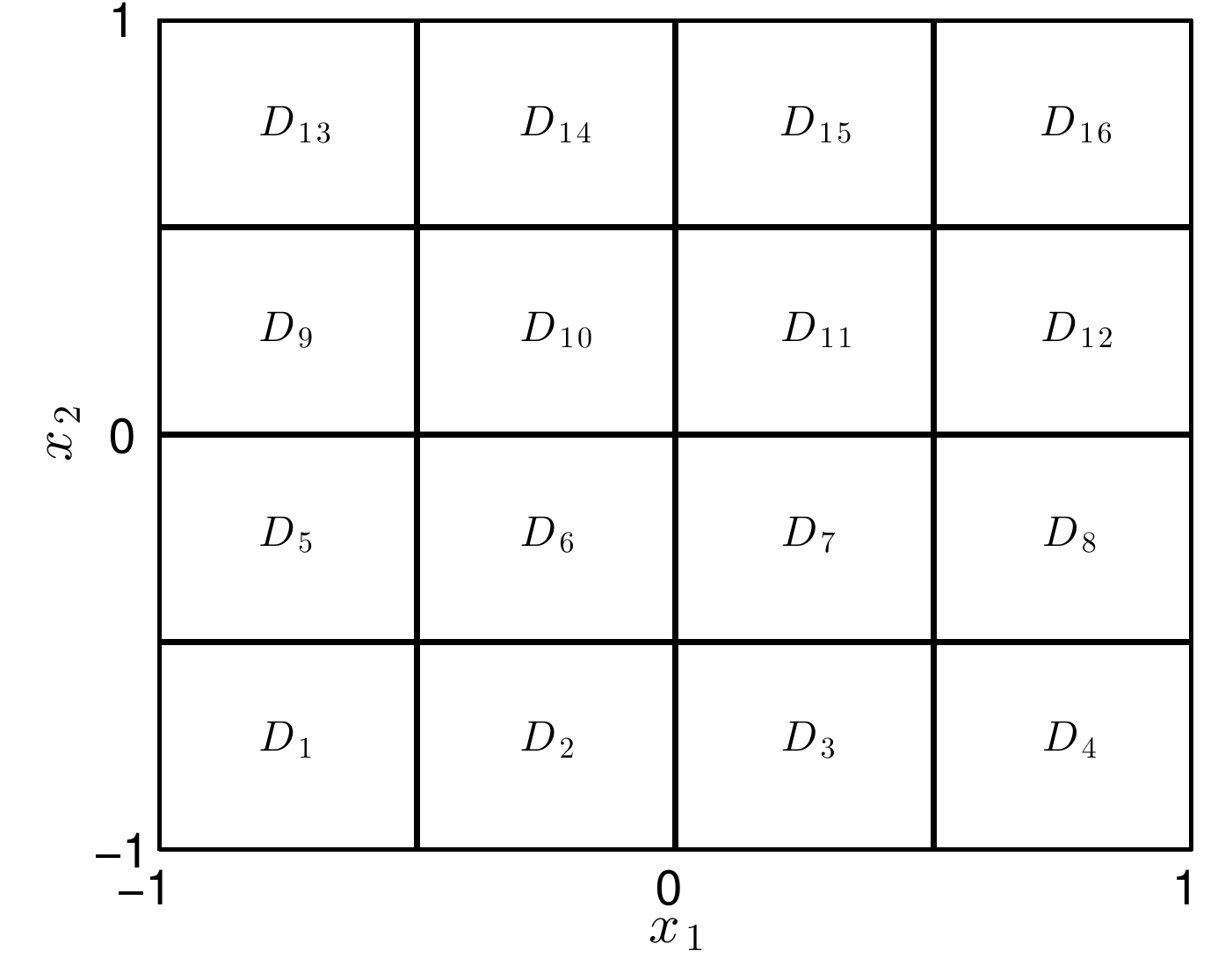}
    \includegraphics[width=4.0cm]{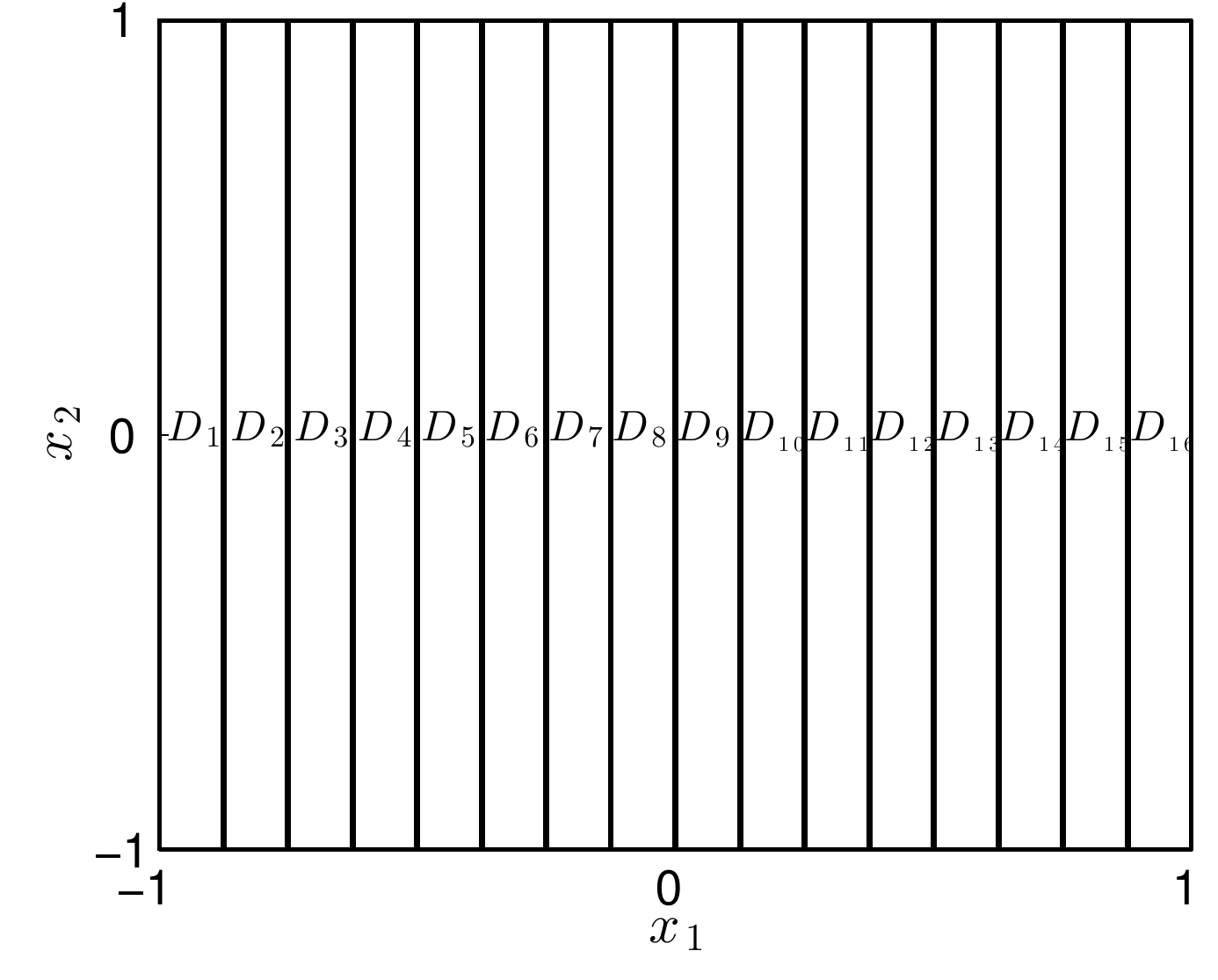}} 
   \caption{Examples of subdomains for test problem \eqref{eq:adiffcoeff} 
   of horizontal strips, squares, and vertical strips with $N_D = 16$.  } \label{fig:subD_rDim16}
\end{figure}

In this paper, 
we consider the steady-state convection-diffusion equation, 
\begin{eqnarray} \label{eq:convdiff}
	- \nabla \cdot a(x,\xi) \nabla u(x,\xi) + w \cdot \nabla u(x,\xi) &=& f(x) \quad\,  \textrm{ in } \,\,\quad D \times \Gamma \\ \nonumber
	u(x,\xi) &=& g_D(x) \quad\, \textrm{ on } \,\,\partial D \times \Gamma 
\end{eqnarray}
where $a(x,\xi)$ is a parametrized random diffusion coefficient, 
$w$ is a constant convective velocity, 
and $f(x)$ is the forcing term. 
The boundary conditions are imposed as 
Dirichlet conditions with $g_D(x)$. 
We will consider choices of 
$a(x,\xi)$ and $g_D(x)$ that involve discontinuity. 
%
The stochasticity in the diffusion coefficient $a(x,\xi)$ 
is taken as a piecewise constant function on partitions, i.e., 
$D$ is divided into $N_D$ equal-sized subdomains of horizontal 
strips, squares, or vertical strips. 
Figure \ref{fig:subD_rDim16} shows three examples 
for $N_D = 16$, namely, $1\times 16$, $4\times 4$, and $16\times 1$, 
on $D = [-1,\,1]\times [-1,\,1]$. 
On each subdomain, we consider 
\begin{equation} \label{eq:adiffcoeff} 
	a(\cdot;\xi)|_{D_m} = \nu \, \xi_m, \quad m = 1, ..., N_D, 
\end{equation}
with a constant $\nu$ and random variables   
$\{\xi_m\}_{m=1}^{N_D}$.  
The dimensionality of the random space is identical to the 
number of subdomains, $M = N_D$. 
The boundary conditions are imposed as Dirichlet conditions 
with 
\begin{equation}
	g_D(x_1,x_2) = \begin{cases} 	
 1, \quad x \in 	
	\{ x\in \partial D | x_1 = -1\} \cup \{ x\in \partial D | -1\leq x_1 \leq 0, x_2 = -1 \}.  \\ 
 0, \quad \textrm{ otherwise } 
\end{cases}. 
\label{eq:convdiffBC} 
\end{equation}
The conditions yield jump discontinuities at the points $(0,\, -1)$ 
and $(-1,\,1)$. 

\begin{figure}
    \centerline{ {\footnotesize \hspace{0.5cm} $2\times 2$, $\quad \nu=1/2$ \hspace{4.8cm} 
$2\times 2$, $\quad \nu=1/20$ } } 
    \centerline{  {\footnotesize  $\E[u(x)]$ \hspace{2.5cm} $\sigma[u(x)]$  \hspace{2.3cm} $\E[u(x)]$ \hspace{2.5cm}  $\sigma[u(x)]$  } } 
    \centerline{ 
    \includegraphics[width=3.5cm]{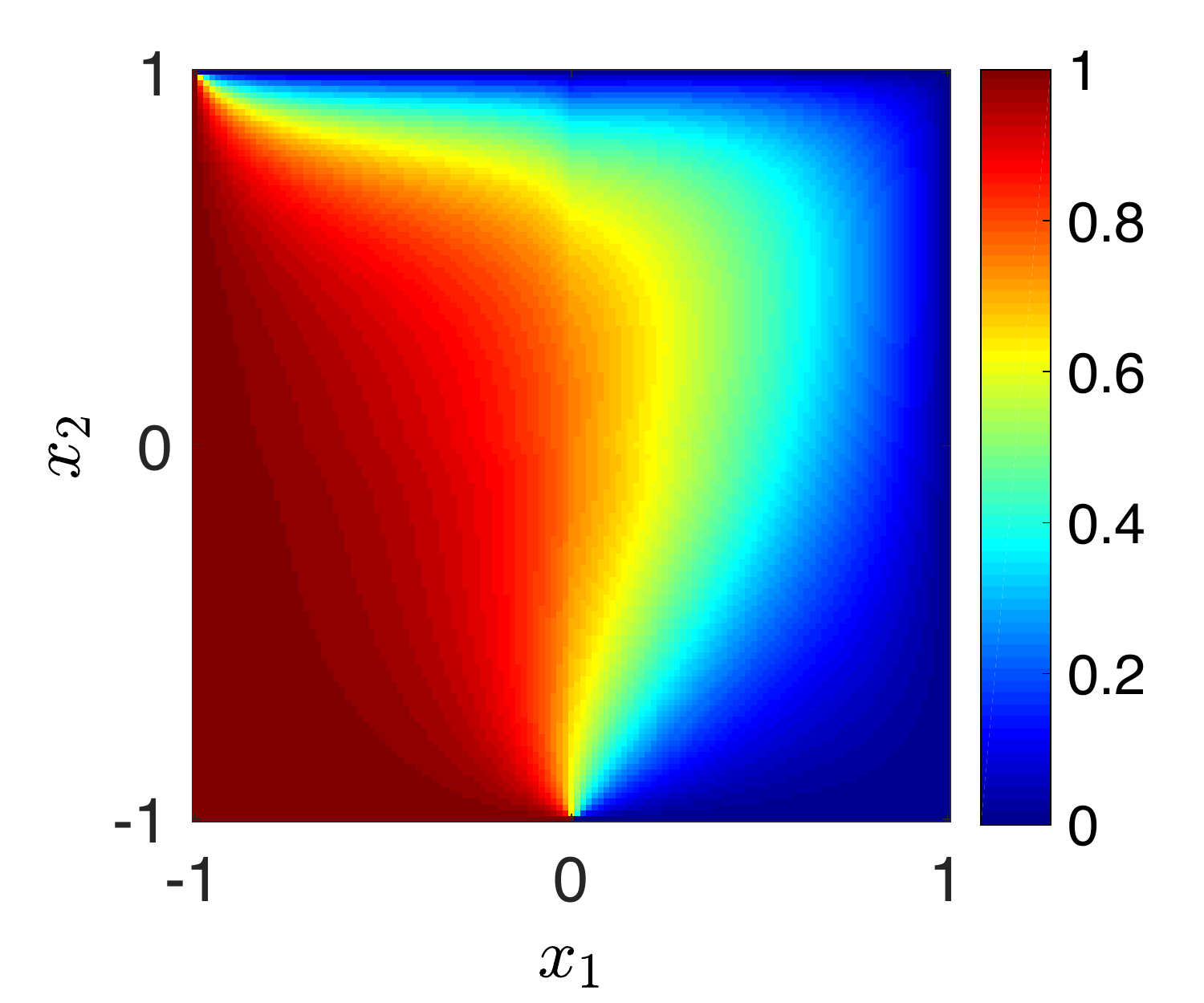}
    \includegraphics[width=3.5cm]{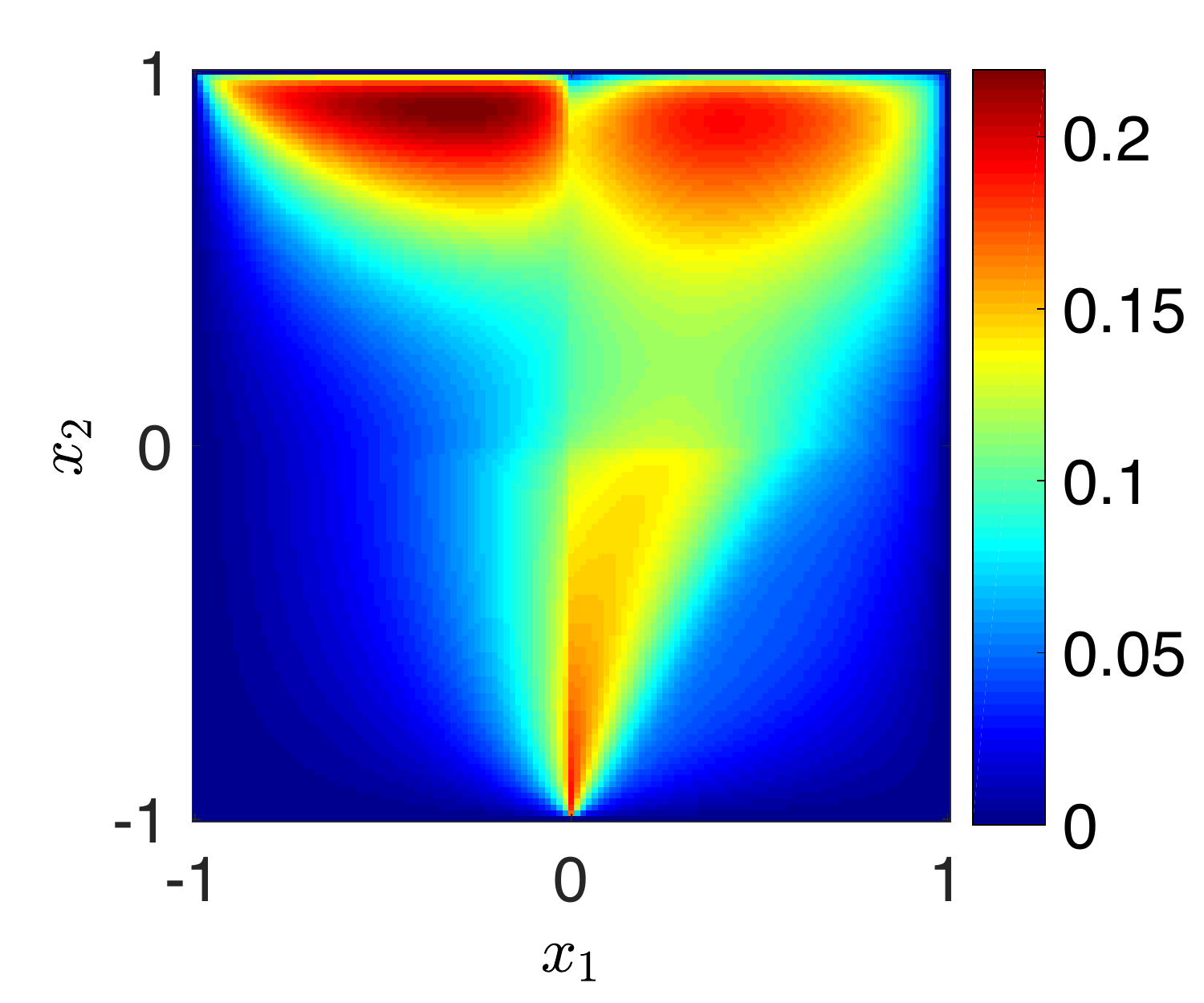}
    \includegraphics[width=3.5cm]{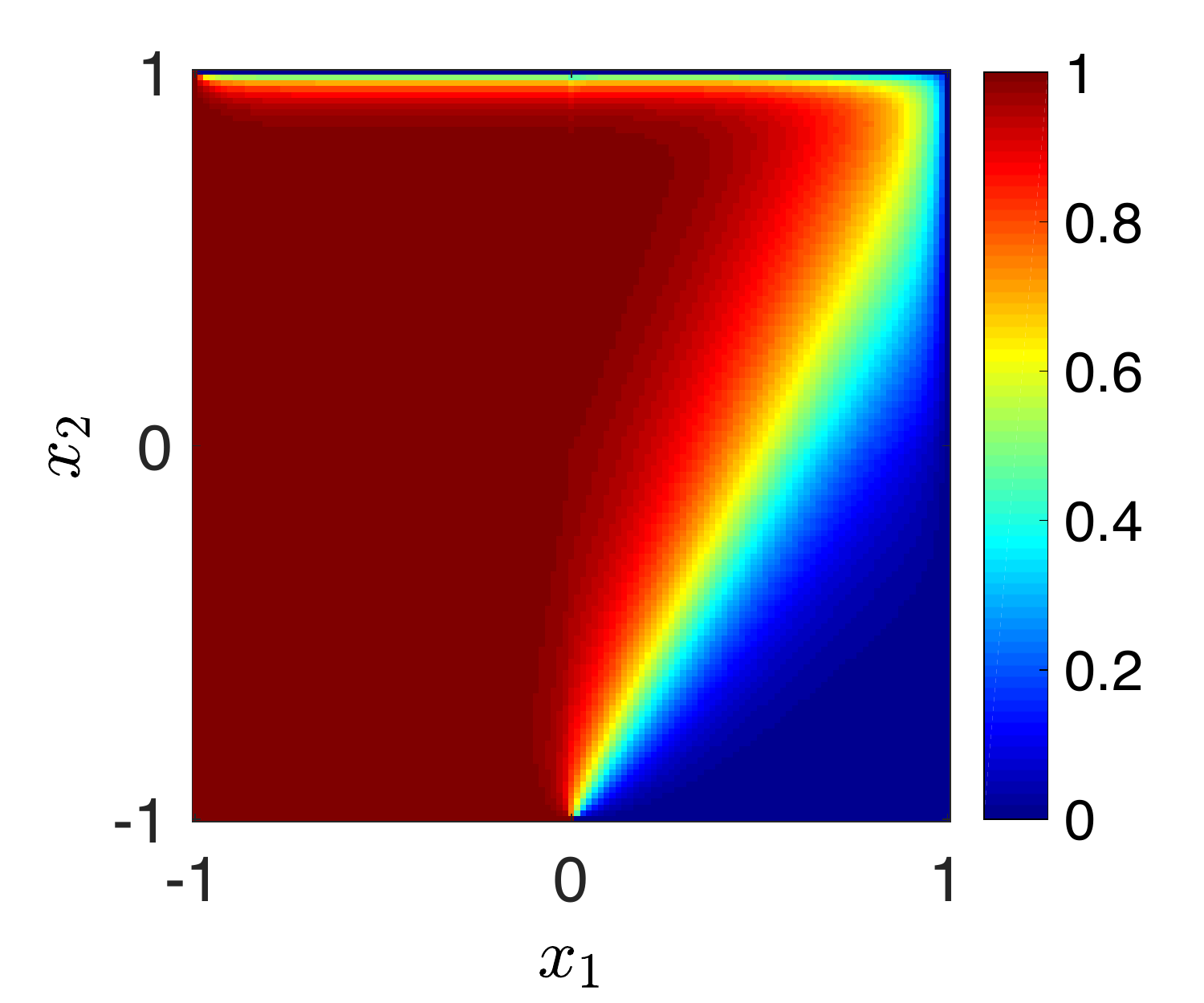}
    \includegraphics[width=3.5cm]{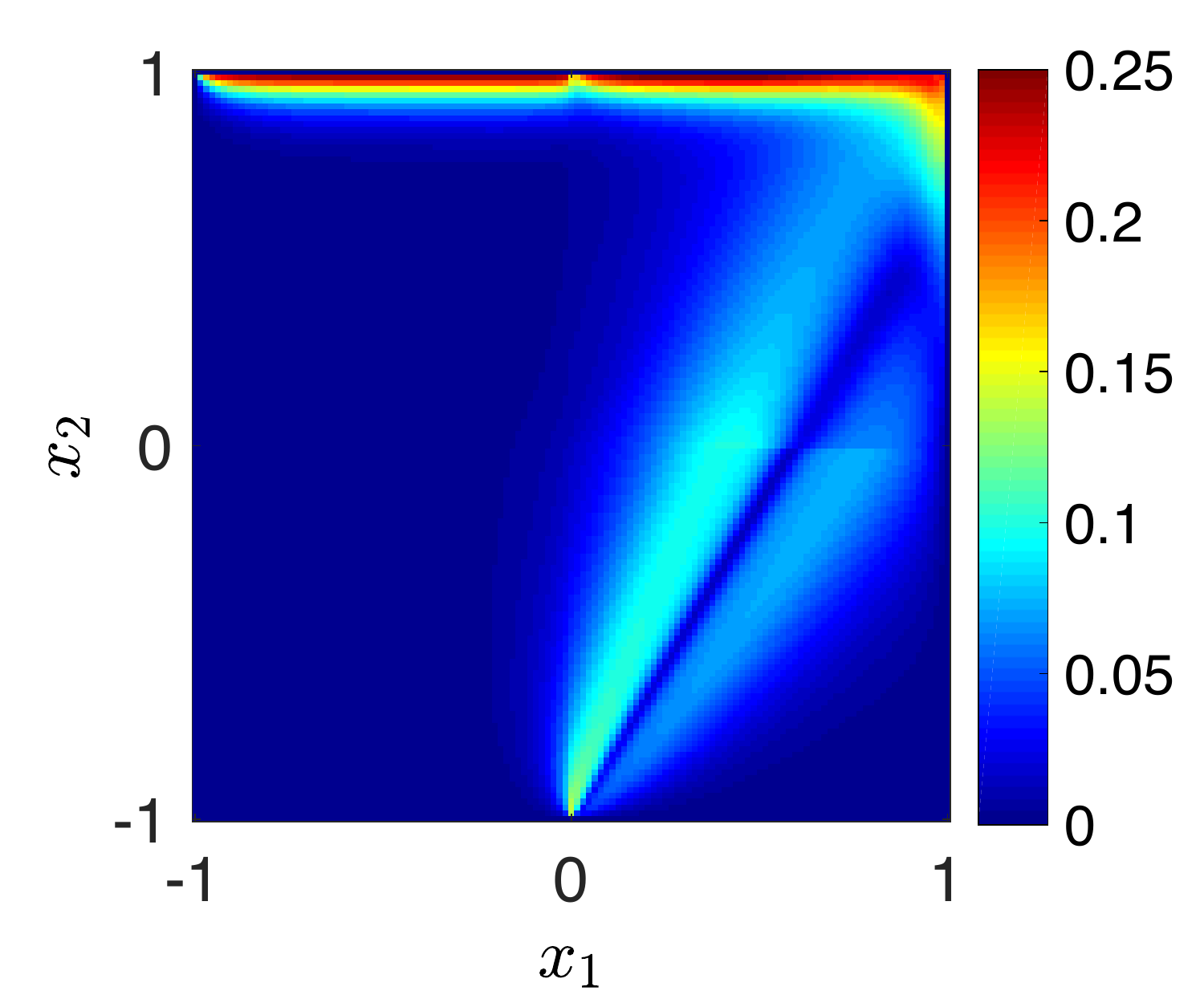}
    }    
   \caption{Examples of mean $\E[u(x)]$ and standard deviation $\sigma[u(x)]$ for diffusion parameter $\nu = 1/2$ and $\nu = 1/20$ on $2\times 2$ partition. } \label{fig:acd_mmt_refsol}
\end{figure}

Examples of solutions to this system with 
a constant convective velocity 
$w = \left( \sin(\pi/6),\, \cos(\pi/6) \right)$ and 
a constant forcing $f=1$ are 
shown in Figure \ref{fig:acd_mmt_refsol}. 
The mean and standard deviation are plotted 
for the case of $N_D = 4$ 
with a $2\times 2$ partition 
and uniform random variables $\{\xi_m\}_{m=1}^{N_D}$ 
on $\Gamma_m = [0.01,\, 1]$ 
with diffusion parameter $\nu = 1/2$ and $1/20$. 
There is an exponential boundary layer near 
the top boundary $x_2 = 1$. 
The discontinuity from the boundary condition 
is smeared by 
the presence of diffusion, producing an internal layer. 
In addition, the standard deviation 
is less regular on the edges of the subdomains 
associated with $a(\cdot,\cdot)$. 
The diffusion coefficient lies in the range $[0.01 \nu, \nu]$, 
and the images on the right show 
results for a smaller diffusion parameter $\nu = 1/20$. 
For this case, the overall variance on the domain tends to 
be smaller than for $\nu = 1/2$, but the magnitude is 
larger in a narrow region near the top boundary layer. 
We will further study the impact of 
randomness on each subdomain 
on the irregular features of the solution 
using the reduced basis method combined with 
ANOVA decomposition. 
This approach will reveal the effects of the anisotropic 
nature of the problem and show the effectiveness 
of the method 
for anisotropic systems while 
obtaining good spatial resolution  
with a reasonable amount of computational cost.

%

\subsection{Collocation methods} 

Collocation methods approximate the solution by interpolation 
at a set of Lagrange polynomials at collocation points, 
\begin{equation}\label{eq:CollocSol} 
	u\left( x,\, \xi \right) \approx 
	\sum_{\xik \in \Theta} 
	u_c\left( x, \xik \right) L_{\xik}(\xi) , 
\end{equation}
where $\Theta \subset \Gamma$ is a finite sample set 
(the collocation points), 
$\{ L_{\xik}(\xi) \}$ are the Lagrange 
interpolation polynomials on $\Gamma$, and 
the coefficient function $u_c\left( x, \xik \right)$ 
is the solution of a deterministic problem corresponding 
to a given realization $\xik$ of the random vector $\xi$. 
We denote the weak form of the system \eqref{eq:convdiff} 
at a given realization $\xi$ 
as $B_{\xi}( u(\cdot, \xi), v) = l(v)$, where
\begin{equation}
	B_{\xi}(u, v) \doteq \int_D \left( a \nabla u \cdot \nabla v  + (w \cdot \nabla u) v \right) dx, \quad l(v) \doteq  \int_D f v \, dx , 
\end{equation}
Let $X_{\h}$ be a finite-dimensional 
finite element space 
of dimension $N_{\h} = O({1}/{\h^d})$ 
defined on $D \subset \R^d$, 
where $\h$ is a grid parameter of the spatial domain. 
We consider a quadrilateral finite element space 
with the streamline diffusion method \cite{FEM_Elman} to 
enhance stability when boundary layers are present. 
Nonhomogeneous boundary conditions are interpolated using 
addition basis functions at the boundary. 
Then the finite element solution $u_{\h}(\cdot,\xi) \in X_{\h}$ 
can be computed as {
\begin{equation} \label{eq:fullwkprob}
	B_{\xi}^{(sd)}( u_{\h}(\cdot, \xi), v) = l(v), \quad ^\forall v\in X_{\h}. 
\end{equation}
where $B_{\xi}^{(sd)}$ incorporates streamline-diffusion terms 
(see \cite{FEM_Elman}, section 6.3 for details).}
For fixed $\xi$, we refer to this 
sample solution as a snapshot. 
Given a finite set of parameters $\Theta$, 
we can gather all the snapshots computed 
on $\Theta$ as 
\begin{equation}
	S_{\Theta} \doteq \{ u_{\h}(\cdot, \xi), \, \xi \in \Theta  \}. 
	\label{eq:CollocSolMat}
\end{equation}
This can be associated with a matrix 
$\bS_{\Theta} \in \R^{N_{\h} \times |\Theta|}$, 
where each column of $\bS_{\Theta}$ is 
a vector of coefficients of basis functions 
of the finite element solution. 
Here, we assume that the finite element approximation 
is sufficiently fine 
so that the spatial discretization error is acceptable. 

\begin{figure}
   \centerline{ {\footnotesize \hspace{0.6cm} $1 \times 4 $ \hspace{2.2cm} $4 \times 1 $ \hspace{2.2cm}  $1 \times 16 $ \hspace{2.2cm} $16 \times 1 $  } }
     \centerline{  
    	\rotatebox{90}{\hspace{0.8cm} {\footnotesize mean}}     	
    \includegraphics[width=3.0cm]{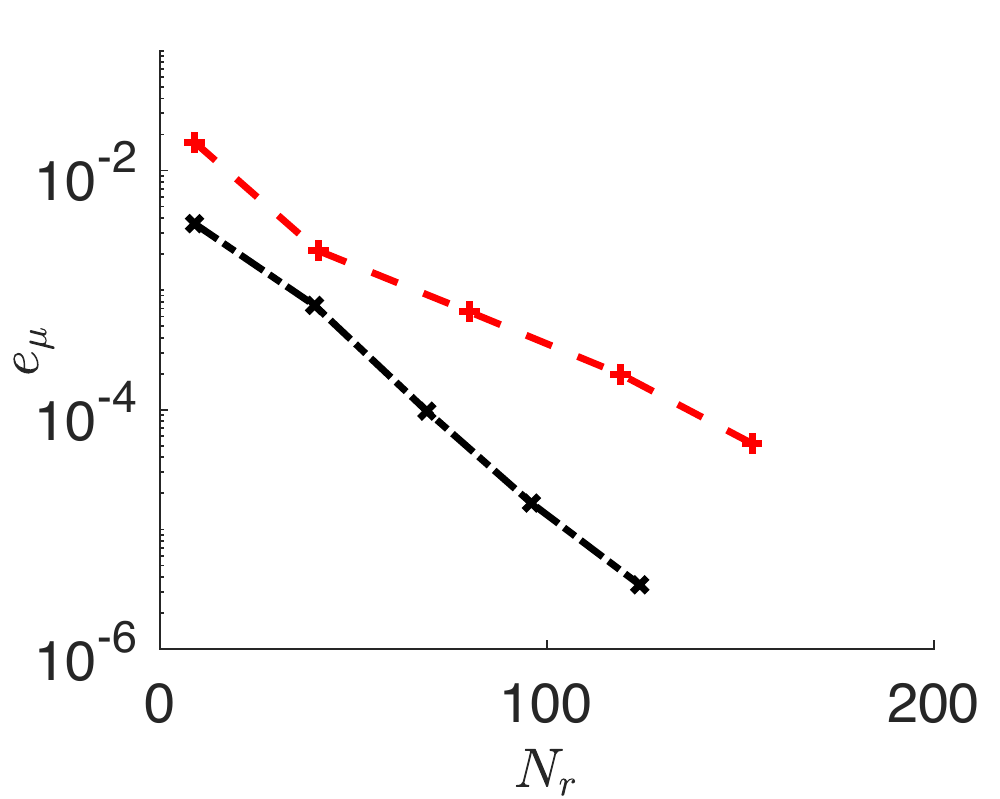} 
    \includegraphics[width=3.0cm]{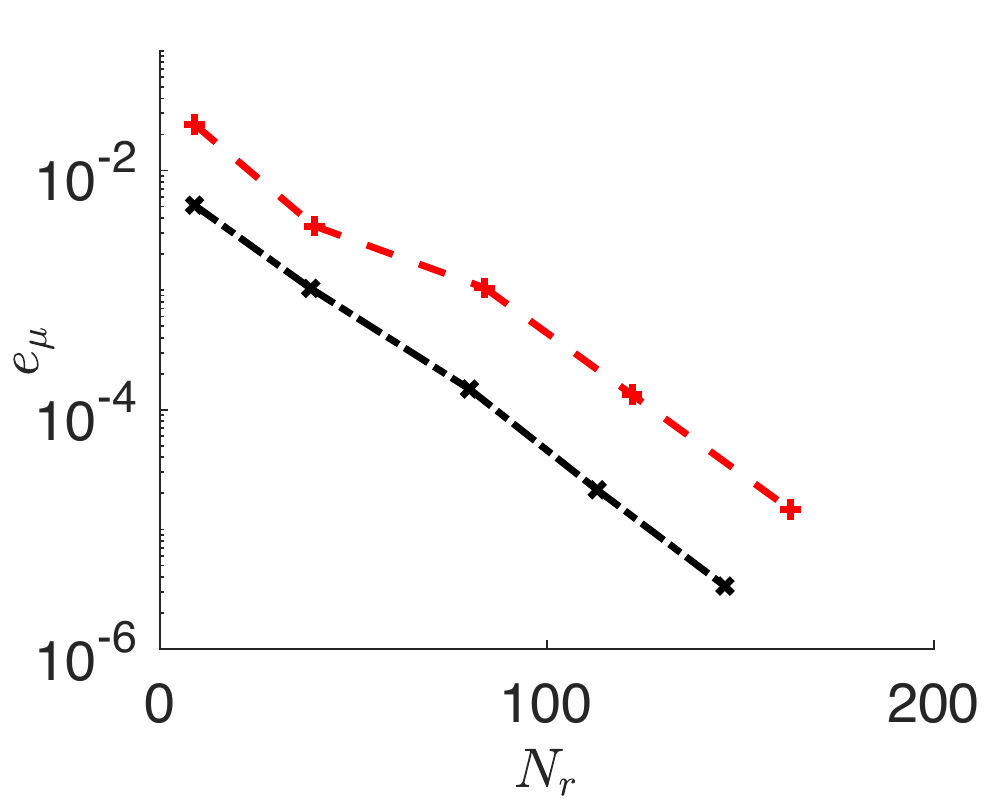} 
    \includegraphics[width=3.0cm]{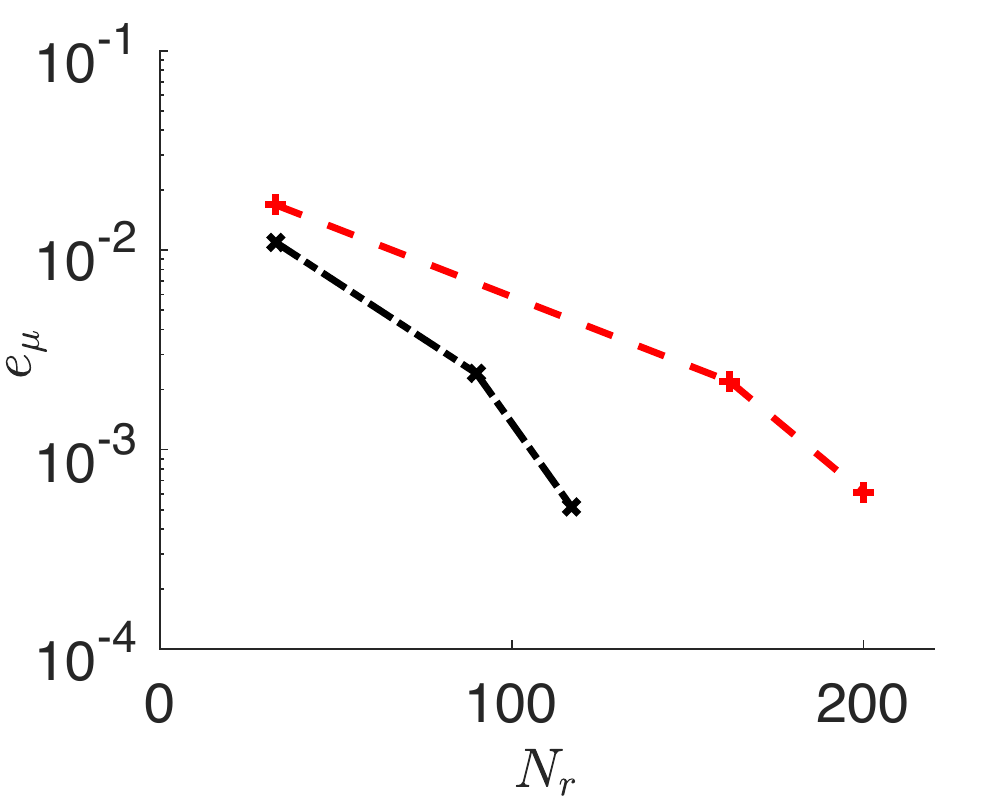} 
    \includegraphics[width=3.0cm]{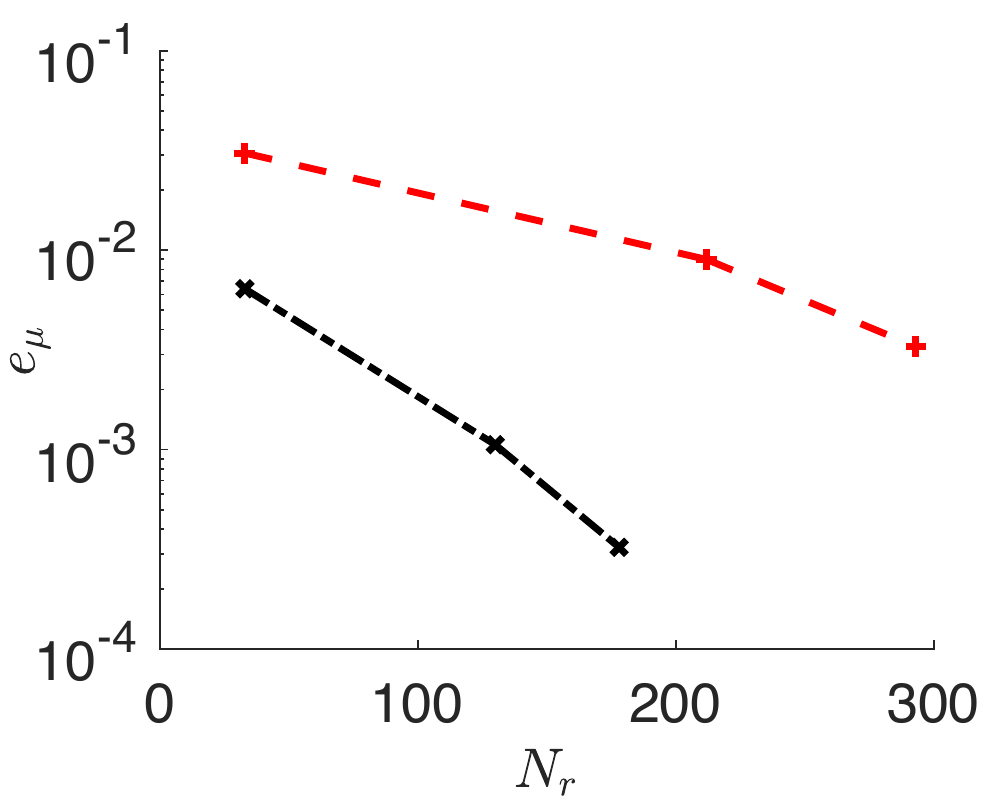} }
    \centerline{  
    	\rotatebox{90}{\hspace{1.1cm} {\footnotesize s.d. }} 
    \includegraphics[width=3.0cm]{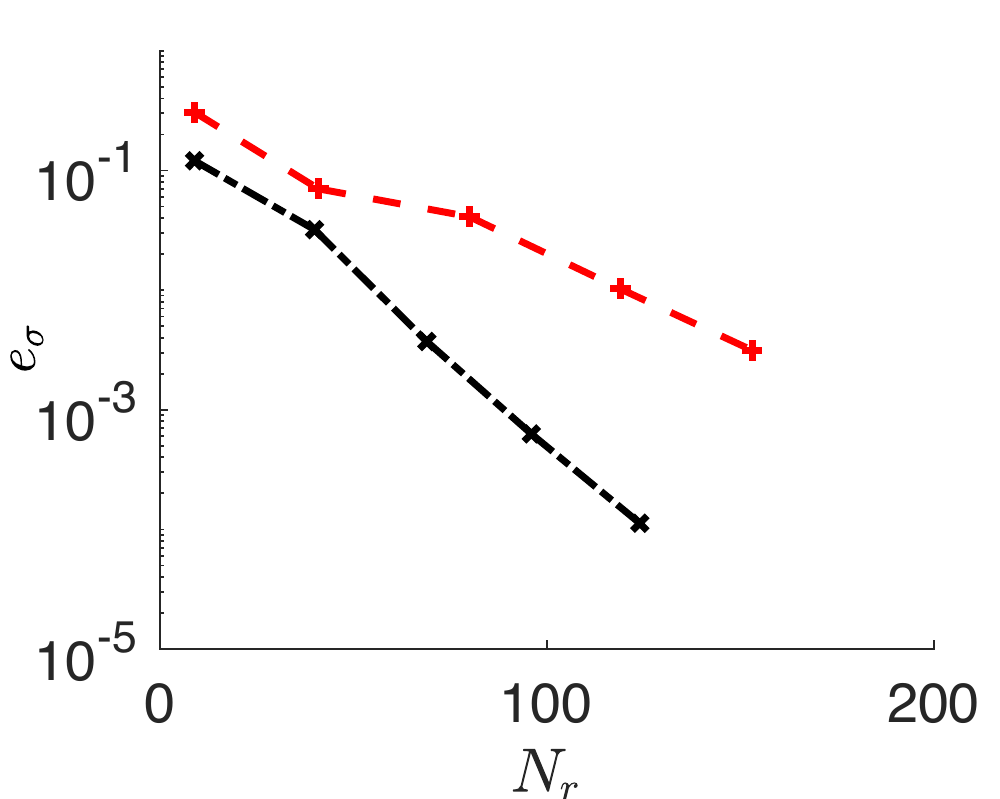} 
    \includegraphics[width=3.0cm]{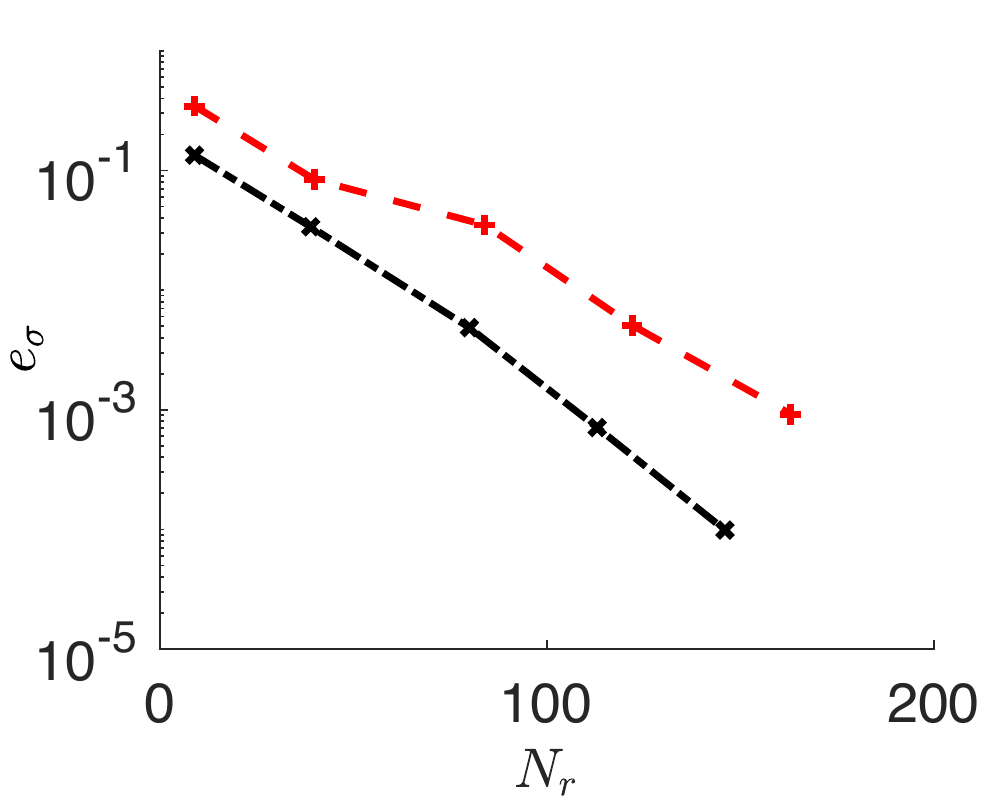} 
    \includegraphics[width=3.0cm]{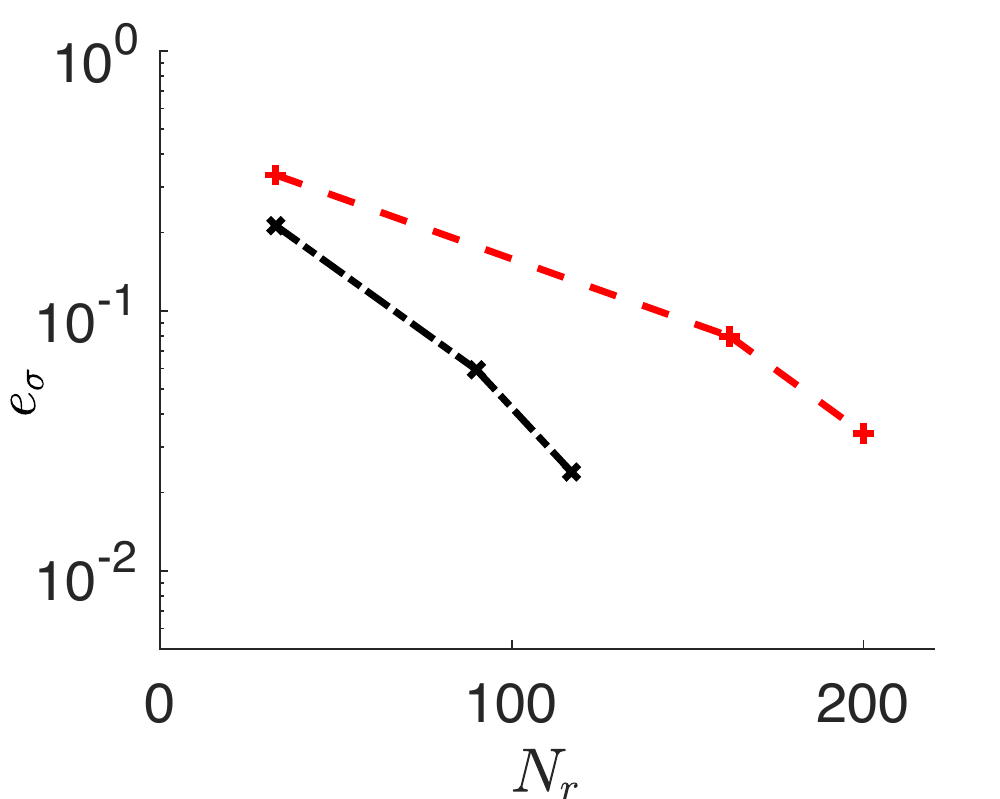}
    \includegraphics[width=3.0cm]{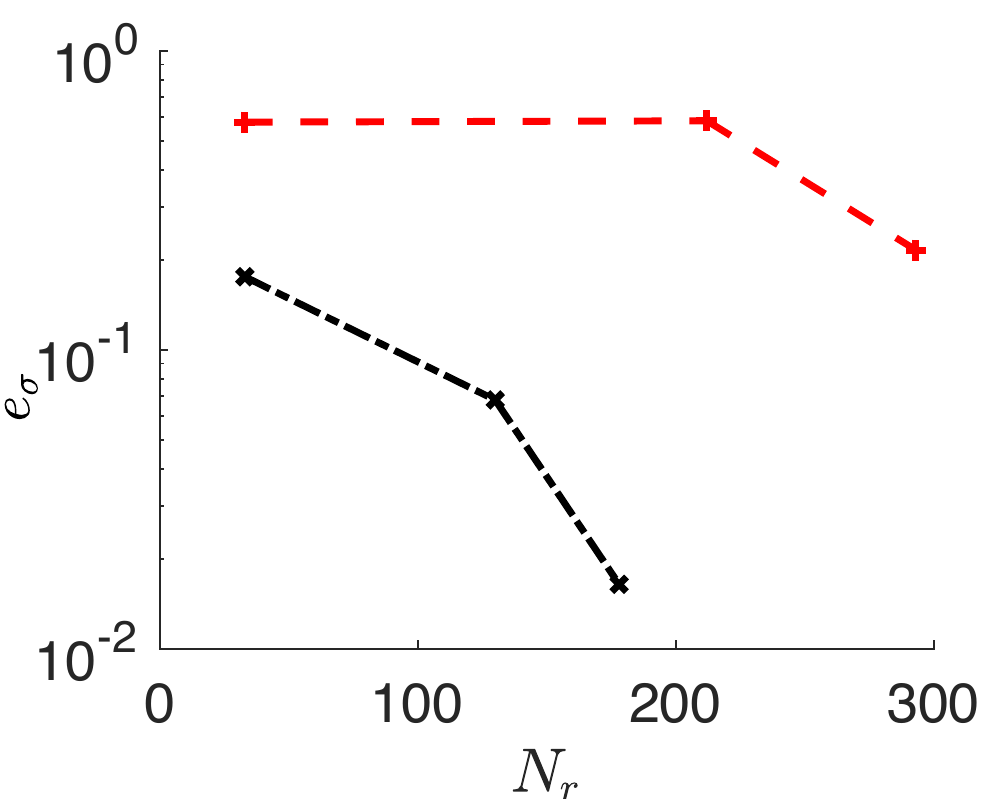}}
    \centerline{
    \includegraphics[width=7.2cm]{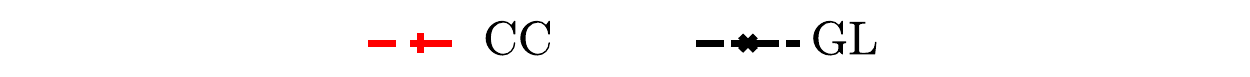}} 
   \caption{Comparison of Clenshaw-Curtis (CC) and Gauss-Legendre (GL) sparse-grid points in the reduced basis method. 
       The plots compare relative $L_2$-error in the mean and standard deviation versus 
       the size of reduced basis $N_r$ for the problem with $\nu = 1/2$. The results are computed with Algorithm \ref{Alg:RBM-} while increasing the sparse-grid level from $\ell=1$ to $5$ using tolerance $\epsilon_{RB} = 10^{-5}$ when $N_D=4$, and from $\ell=1$ to $3$ with $\epsilon_{RB} = 10^{-4}$ when $N_D=16$.}
       \label{fig:errallind_rDim16}
\end{figure}

The full collocation solution \eqref{eq:CollocSolMat} 
can be computed by applying a 
deterministic solver $|\Theta|$ times.
If $\Theta$ consists of a full tensor product of 
one-dimensional parameter values, i.e.,  
$\Theta = \Theta_1^{i_1} \otimes \cdots \otimes  \Theta_M^{i_M}$, 
where $\Theta_m^{i_m}$ is the set of collocation points in 
the $m$-th direction of order $i_m$, 
costs grow exponentially 
with respect to the dimension $M$. 
The {\em sparse-grid collocation} method \cite{Gri05,NobileTW_SINA,ZabarasSG,Sparse_Novak} 
reduces this cost using a subset of the full 
tensor grid collocation points 
$$
\Theta_{\ell} = \cup_{\ell+1 \leq |\bi| \leq \ell+M} \left( \Theta_1^{i_1} 
\otimes \cdots \otimes  \Theta_M^{i_M} \right), 
$$
where $|\bi| = i_1 + \cdots + i_M$ and the {\em level} $\ell$ 
is the parameter that 
truncates the level of interaction. 
%
The choice of sparse-grid collocation points are typically 
taken to be nodes used for quadrature in $M$-dimensional space.
The most common ones are Clenshaw-Curtis and Gauss abscissae. 
The convergence behavior of the reduced basis collocation method  
using Clenshaw-Curtis quadrature is studied in 
\cite{QifangRBM}. 
Here, we compare 
Clenshaw-Curtis and Gauss-Legendre points  
used in the  
reduced basis collocation method 
that will be described later.  
Figure \ref{fig:errallind_rDim16} shows 
the error in the moments with respect to the size of 
the reduced basis $N_r$. 
It is evident 
that the choice of the interpolation point is critical 
and Gauss-Legendre is the more effective choice 
for this example. 
In the sequel, we restrict our attention 
to Gauss-Legendre quadrature.

\subsection{ANOVA decomposition for collocation methods}

An alternative to sparse-grid collocation 
uses an ANOVA decomposition of 
an $M$-dimensional function. 
This decomposition has the form 
\begin{eqnarray}
\u (\,\cdot,\,\xi_1,...,\xi_M) &=& \su_0(\cdot) + \sum_{j_1=1}^{M} \su_{j_1}(\,\cdot,\,\xi_{j_1}) + 
\sum_{j_1<j_2}^{M} \su_{j_1j_2}(\,\cdot,\,\xi_{j_1}, \xi_{j_2}) + \cdots 
\label{eq:ANOVA} 
\end{eqnarray} 
This representation decomposes multivariate functions 
according to the degree of interaction among the variables. 
For instance, $\su_{j_1j_2}(\,\cdot,\,\xi_{j_1}, \xi_{j_2})$ 
describes the second-order or {\em two-body} interaction 
between the variables 
$\xi_{j_1}$ and $\xi_{j_2}$. 
Similar to sparse-grid collocation, 
by truncating the series at a certain level, 
an $M$-dimensional function can be approximated using 
a function on a lower-dimensional subspace of 
the $M$-dimensional parameter space. 
Such a decomposition 
can be computed by 
$\su_K(\,\cdot,\,\xi_K) \doteq \int_{\Gamma_{K'}} u(\,\cdot,\,\xi) d\mu ( \xi_{K'} )  
	- \sum_{S \subset K} \su_T(\,\cdot,\,\xi_S) $,  
where 
$K \subset \{1,\cdots,M\}$ and 
$K'$ is the complement of $K$ in $\{1,\cdots,M\}$. 

This construction can be put in a form 
similar to that of collocation 
by replacing the Lebesgue measure 
with a Dirac measure $d \mu = \delta(\xi - \bc)$. 
This is called the anchored-ANOVA method 
where $\bc$ is the anchor point 
\cite{Foo1,Schwab09}. 
The first few terms in the anchored-ANOVA decomposition become 
	$\su_0(\cdot) = u(\,\cdot,\,\bc) = u(\,\cdot,\,c_1, \cdots, c_M)$, $\su_j(\,\cdot,\,\xi_j) = u(\,\cdot,\,c_1, .., \xi_j, .., c_M) - \su_0(\cdot)$, and 
	$\su_{ij}(\,\cdot,\,\xi_i,\xi_j) = 	 u(\,\cdot,\,c_1, ..\xi_i, .. , \xi_j, .., c_M) - \su_i(\,\cdot,\,\xi_i) - \su_j(\,\cdot,\,\xi_j) - \su_0(\cdot)$. 
%
By considering the ANOVA approximation in the random space 
discretized by the probabilistic collocation method (PCM) \cite{BabuskaNT_SINM05,XiuH_SISC05} 
in each parametric space, 
the solution approximated up to level $\ell$ 
can be written as 
\begin{equation} 
	u(\cdot,\xi) = \sum_{ K\subset\{1,...,M\}, |K|\leq\ell} 
	\varkappa_{M,|K|,\ell} \left[ \sum_{\xik \in \Xi_{K}^p} 
	u_K\left( \cdot, \xik \right) L_{\xik}(\xi) \right], 
\label{eq:sol_ANOVA}
\end{equation} 
where $\varkappa_{M,|K|,\ell}$ 
is a signed number 
counting the frequency of the sample solutions at $\xik \in \Xi_{K}^p$ appearing in the ANOVA  
expansion up to level $\ell$, that is, 
$\varkappa_{M,j,\ell} = \sum_{r=j}^{\ell} \left( -1 \right)^{r-j} \frac{(M-j)!}{(r-j)!(M-r)!}$ 	
and  $\Xi_{K}^p$ is the set of PCM-ANOVA collocation points 
\cite{Foo1,Handy_SISC_ANOVA_2012} in  
direction $K$ up to polynomial order $p$ for 
all the parameters. 
For instance, 
if $\{\xi_{m}^{\kappa_m}\}_{\kappa_m=1}^{p_m}$ are the PCM points 
in the $m$-th parameter of polynomial order $p_m$, then 
\begin{eqnarray}
\label{eq:ANOVApt}
	\Xi_{1}^{p_1} &=& \{ (\xi_{1}^{\kappa_1}, c_2,...,c_M)\in \Gamma \,|\,  1\leq \kappa_1 \leq {p_1}\},  \\ 
	\nonumber 
	\Xi_{(1,2)}^{(p_1,p_2)} &=& \{ (\xi_{1}^{\kappa_1}, \xi_{2}^{\kappa_2}, c_3,...,c_M)\in \Gamma \,|\, 1\leq \kappa_1 \leq {p_1},\, 1\leq \kappa_2 \leq {p_2}\}. 
\end{eqnarray} 
%
%
We denote the set of 
collocation points up to level $\ell$ and its size as 
\begin{equation}
\Xi_{\ell} \doteq \cup_{|K| \leq \ell} \Xi_K^p, \qquad 
C( M, \ell) \doteq |\Xi_{\ell}| = \sum_{l = 0}^{\ell}  
\left( \begin{matrix} M \\ l  \end{matrix}  \right) p^l, 
	\label{eq:ANOVAcost} 
\end{equation}
respectively. 
%
For stochastic problems, 
it was shown in \cite{ZabarasANOVA} that an ANOVA 
decomposition of level two has good accuracy 
if the anchor point is chosen as the mean of the probability 
density function considered, and we will use this technique. 

An error indicator for the ANOVA truncation can be 
constructed using the relative size of 
either the mean or variance in each term \cite{xiuchoigk,ZabarasSG}. 
In our work, we take the ANOVA indicator for the $K$-th 
term using the mean as 
\begin{equation} 
\label{eq:anovaIndc} 
\gamma_K \doteq \frac{\| \E[ \su_K ] \|}{ \sum_{|S| < |K|} \| \E[ \su_S ] \| }.
\end{equation} 

\subsection{Reduced basis collocation method} 
\label{sec:RBMANOVA}
Reduced basis algorithms aim to reduce the computational 
cost by approximating the solution 
using full discrete solutions computed only 
on a small subset of parameters $\Theta_r \subset \Xi_{\ell}$. 
Let $N_r \doteq |\Theta_r|$. 
The goal is to find a set of finite element solutions 
$Q = \{q_1, \cdots, q_{N_r} \} $ computed 
on $\Theta_r$ 
that give an accurate representation of 
the column space of $\bS_{\Xi_{\ell}}$. 
Such an approximate solution 
$u_{rb}(\cdot,\xi) \in span(Q)$ satisfies 
\begin{equation} \label{eq:rbwkprob} 
	B_{\xi}\left( u_{rb}(\cdot, \xi), v \right) 
	= l(v),\quad ^\forall v \in span (Q). 
\end{equation} 
The size $N_r$ of the reduced problem \eqref{eq:rbwkprob} 
is typically much smaller 
than $N_{\h}$ from \eqref{eq:fullwkprob}, 
and the utility 
of reduced basis methods 
derives from this fact 
\cite{RBM_Patera1,Quatrni_RBMbook}. 
Here we summarize the procedure.

Assume that the operator $\L$ is linear and is affinely dependent 
on the vector of parameters $\xi$, 
that is, $\L = \sum_{i=1}^{\hat{n}_a} 
\phi_i(\xi) \L_i$. 
Then, the weak form \eqref{eq:fullwkprob} 
requires solution of a linear system of order $N_{\h}$, 
\begin{equation} \label{eq:fullSys} 
	 \bA_{\xi} \bu_{\xi} = \bff, 
\end{equation} 
where $\bA_{\xi} = \sum_{i=1}^{\hat{n}_a} \phi_i(\xi) \bA_i $, 
and $\bA_i$ are parameter-independent matrices. 
For each $q_i \in Q$, let $\bq_i$ denote the vector of 
coefficients associated with $q_i$. 
Then $\bQ = [\bq_1,\cdots,\bq_{N_r}] 
\in \R^{N_h \times N_r}$ is the matrix representation of $Q$, 
and the linear system 
obtained from a Galerkin condition 
corresponding to the reduced problem 
\eqref{eq:rbwkprob} can be written as 
\begin{equation}
	\bQ^T \bA_{\xi} \bQ \wbu_{\xi} = \sum_{i=1}^{\hat{n}_a} \phi_i(\xi) \left( \bQ^T \bA_i \bQ \right) \wbu_{\xi} = \bQ^T \bff, 
\end{equation}
where $\bQ \wbu_{\xi}$ is the reduced basis solution 
approximating $\bu_{\xi}$. 
If the parameter-independent quantities $\bQ^T \bA_i \bQ$ and 
$\bQ^T \bff$ are precomputed, 
the cost to assemble the reduced system for each $\xi$ 
becomes $O(N_r^2)$ times $\hat{n}_a$. 
The reduced basis algorithm computed 
with collocation points $\Theta$ 
is summarized in Algorithm \ref{Alg:RBM-}. 
Either sparse-grid or ANOVA collocation 
can be employed by updating the reduced basis matrix 
while increasing the level from $0$ to $\ell$, 
that is, 
$\bQ$ = RBM($\bQ$, $\Xi_l$, $\epsilon$), for $0\leq l \leq \ell$ 
(see \cite{QifangRBM,QifangRBMA} for details). 
Here, $\epsilon_{RB}$ is an acceptance 
tolerance {for the reduced basis 
solution based on an error indicator.} 
%
We use a residual based error indicator, 
\begin{equation} \label{eq:rbmindc} 
	\eta_{Q,\xik} \doteq \frac{ \|\bA_{\xik} \bQ \wbu_{\xik} - \bff\|_2 }{\| \bff \|_2 } \wk, 
\end{equation} 
where $\wk = \mathbb{P}(\xi = \xik )^\alpha$, 
$0 \leq \alpha \leq 1$ \cite{XWGK_JCP05,PChen2}. 
An alternative error indicator is 
a dual-based indicator developed in 
\cite{RBM_Rozza1,RBM_Patera1}, 
which depends on having an 
a posteriori error estimate for the associated problem.

\begin{algorithm}
\caption{ $(\bQ,\, \bS_{\Theta_r})$ = RBM($\bQ$, $\Theta$, $\epsilon_{RB}$) : Reduced basis algorithm to update the reduced basis matrix $\bQ$ considering the collocation points $\Theta$ and tolerance $\epsilon_{RB}$. 
} 
$\,\,\,$1: \textbf{for} $\,$ $1 \leq k \leq |\Theta|$  \textbf{do}  \\ 
$\,\,\,$2: $\quad$ Compute $\wbu_{\xik}$ by solving \eqref{eq:rbwkprob} and the error indicator $\eta_{\bQ,\xik}$ of \eqref{eq:rbmindc}.  \\ 
$\,\,\,$3: $\quad$ \textbf{if} $\,$ $\eta_{\bQ,\xik} < \epsilon_{RB} $ \textbf{then} \\ 
$\,\,\,$4: $\quad\quad$ Use the reduced solution $u_R(\cdot, \xik) = \bQ \wbu_{\xik}$ to serve as $u_c(\cdot,\xik)$ in \eqref{eq:CollocSol}. \\ 
$\,\,\,$5: $\quad$ \textbf{else} \\  
$\,\,\,$6: $\quad\quad$ Compute the full solution vector $\bu_{\xik}$ by solving \eqref{eq:fullSys}.  \\ 
$\,\,\,$7: $\quad\quad$ Use $u_h(\cdot,\xik)$ derived from $\bu_{\xik}$ to serve as $u_c(\cdot,\xik)$ in \eqref{eq:CollocSol}. \\ 
$\,\,\,$8: $\quad\quad$ Update $\bS_{\Theta_r}$ and $\bQ$ by augmenting $\bu_{\xik}$ and the orthogonal complement. \\ 
$\,\,\,$9: $\quad\quad$ Reconstruct the offline reduced matrices and vectors, 
$\bQ^T \bA_i \bQ$, $\bQ^T \bff$, $\bQ^T \bA_i \bA_j \bQ$, and $\bQ^T \bA_i \bff$. \\ 
10: $\quad$ \textbf{end} \\
11: \textbf{end} 
\vspace{.1cm}
\label{Alg:RBM-} 
\end{algorithm} 

The total cost of a reduced basis algorithm 
can be divided into two parts as 
\begin{eqnarray} \label{eq:CostGreedy}
 C(M, \ell) \cdot 
 (\textit{cost of reduced system solves and residual calculation}) + \hspace{0.4cm} \\ \nonumber
 N_r \cdot  
(\textit{cost of full system solves and construction of offline quantities}),  
\end{eqnarray}
where 
the former and latter terms are 
on the order of $N_r$ and $N_h$, respectively. 
Here,  
$C(M, \ell)$ is the number of 
sparse-grid nodes or ANOVA collocation points. 
The efficiency of the reduced basis algorithm comes from 
avoiding the full system solve 
on some collocation points.  
%
However, the cost 
becomes large as $M$ and $\ell$ increase.  
Dimension-adaptive approaches 
can further improve the efficiency for anisotropic problems. 
The ANOVA indicator $\gamma_K$ \eqref{eq:anovaIndc} can be used 
to select only the important directions. 
After each level $l$,  
we choose the ANOVA terms 
$u_K$ \eqref{eq:sol_ANOVA} 
such that $\gamma_K$ is greater than 
a tolerance $\epsilon_A$, 
and refer to the component of these terms as
the {\em effective dimensions} $\J_l$, i.e., 
\begin{equation}	
\J_{l} = \{ |K| = l \,|\, \gamma_K > \epsilon_{A} \}. 
\label{eq:EffecitveDim} 
\end{equation}
We use these indices to build 
the next level ANOVA approximation 
(see lines 19--20 of Algorithm \ref{Alg:RBM-ANOVA_adaptp} 
for details). 
The number of search points is 
reduced to $C(\widetilde{M}, \ell)$, 
where $\widetilde{M}$ is the number of the effective dimensions. 
This is a standard procedure in adaptive ANOVA schemes 
\cite{xiuchoigk,ZabarasSG} 
and has been studied in the context of reduced basis methods 
in \cite{QifangRBMA,JanH_RBMANOVA}.

Here, we emphasize that we also employ 
the ANOVA indicator $\gamma_K$ \eqref{eq:anovaIndc} 
to sort the search direction $\J_{\ell}$ and 
the newly added reduced basis solution $\bS_{N_r}$. 
Both are sorted in each level 
according to the descending order of $\gamma_K$. 
Since the ANOVA indicator implies  
the priority of the reduced basis, 
this procedure enhances the quality of the reduced basis and 
saves the additional cost that occurs from sorting 
the reduced basis afterwards. 
\begin{algorithm}
\caption{ $(\bS,\J)$ = RB-Sort($\{ \gamma_K \}$, $\{ \bS_{K} \}$, $\J_{\ell}$) : Sort the reduced basis and ANOVA index in the descending order of ANOVA indicator $\gamma_K$. 
} 
$\,\,\,$1: Initialize $\bS =[\,]$ and $\J = [\,]$. \\ 
$\,\,\,$2: \textbf{while} $\,$  $\J_\ell \neq \emptyset $  \textbf{do} \\ 
$\,\,\,$3: $\quad$ Find the index with the maximum ANOVA indicator $L  = \arg\max \{ \gamma_K \, | \, K \in \J_\ell \}$.  \\ 
$\,\,\,$4: $\quad$ Save the sorted reduced basis $\bS = [\bS, \bS_L ]$ and ANOVA index $\J = [\J, L ]$.  \\ 
$\,\,\,$5: $\quad$ $\J_\ell = \J_\ell \backslash L$.  \\  
$\,\,\,$6: \textbf{end} 
\vspace{.1cm}
\label{Alg:RBM-sort} 
\end{algorithm}

\subsection{Numerical results} \label{sec:ResultANOVA} 

In this section, we use the reduced basis algorithm based on ANOVA 
to examine the features of the benchmark problem 
\eqref{eq:convdiff} 
and the accuracy of the method. 
We consider the spatial domain $D = [-1,1]\times[-1,1]$, 
with the system parameters 
described in the beginning of section \ref{sec:RBM-SG}. 
We take the random variables $\xi_m$ in 
the piecewise constant diffusion coefficient 
$a(x,\xi)$ of \eqref{eq:adiffcoeff} 
to be independent and uniformly distributed 
on $\Gamma_m = [0.01,\,1]$. 
The magnitude of the diffusion parameter  
$\nu$ is either $1/2$ or $1/20$ 
for examples of moderate and small diffusion. 
We choose the convective velocity $w$ to correspond to a 
30-degree angle right of vertical, 
and the Dirichlet boundary condition $g_D(x_1,x_2)$ as in 
\eqref{eq:convdiffBC}. 
For each sample point, 
we find a weak discrete solution $u_h(x,\xi)\in X_{\h}$. 
For the finite element discretization of $D$, we use a 
bilinear ($Q_1$) finite element approximation \cite{FEM_Braess}. 
The deterministic solver is based on the IFISS software package 
\cite{IFISS,IFISS_sw} 
and the spatial discretization is done 
on a uniform $128\times 128$ grid, 
so that $N_{\h} = 128^2$. 

\begin{figure}  
    \centerline{ 
    \includegraphics[width=4.0cm]{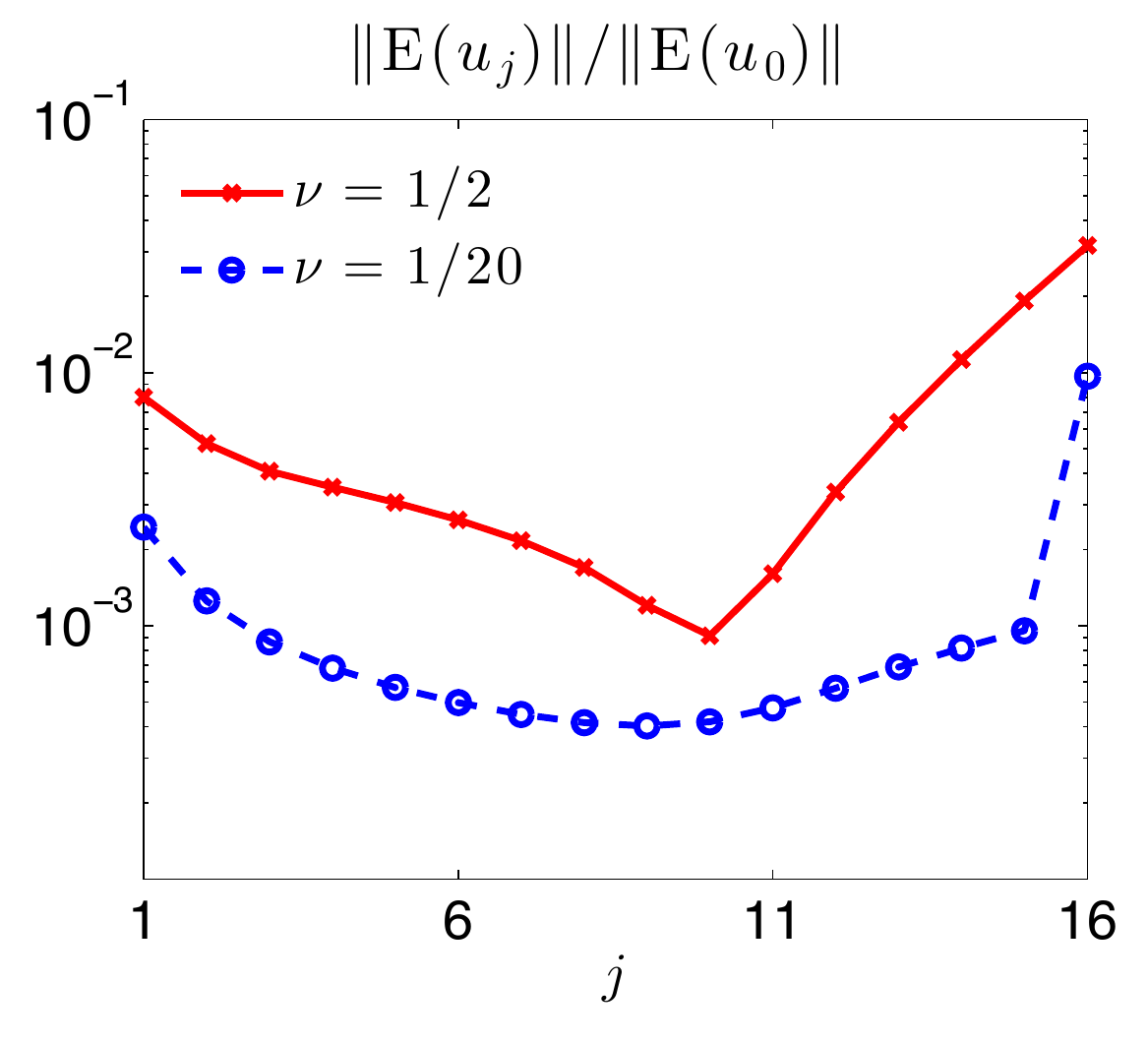} \hspace{0.2cm}
    \includegraphics[width=4.5cm]{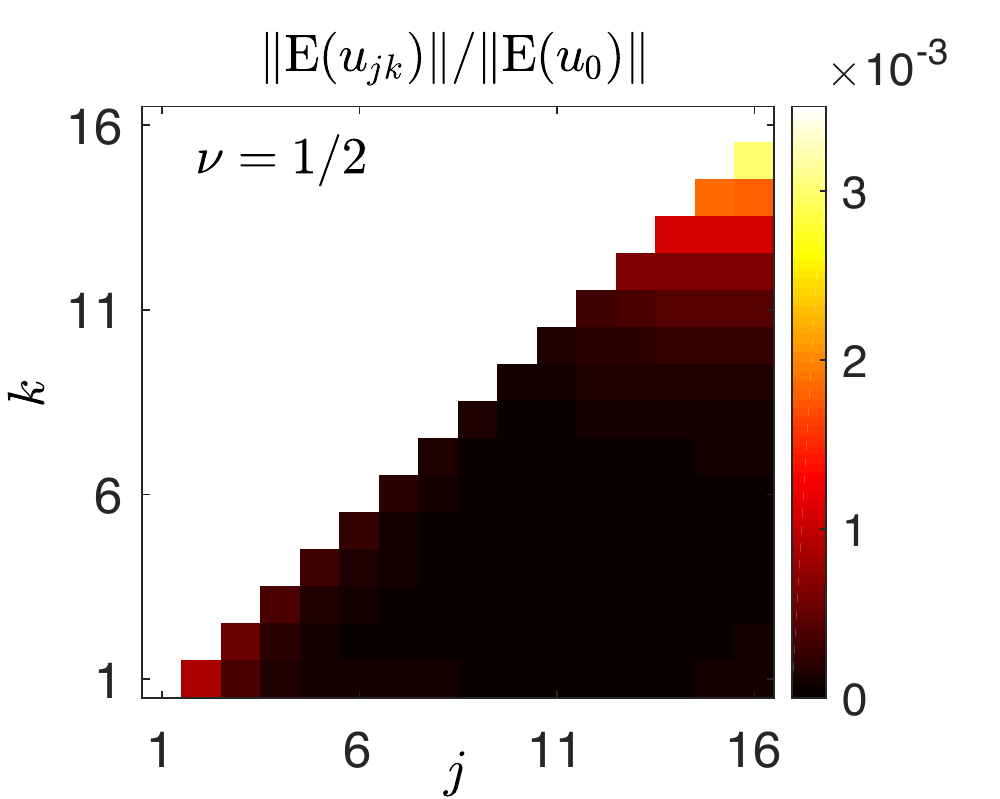}
    \includegraphics[width=4.5cm]{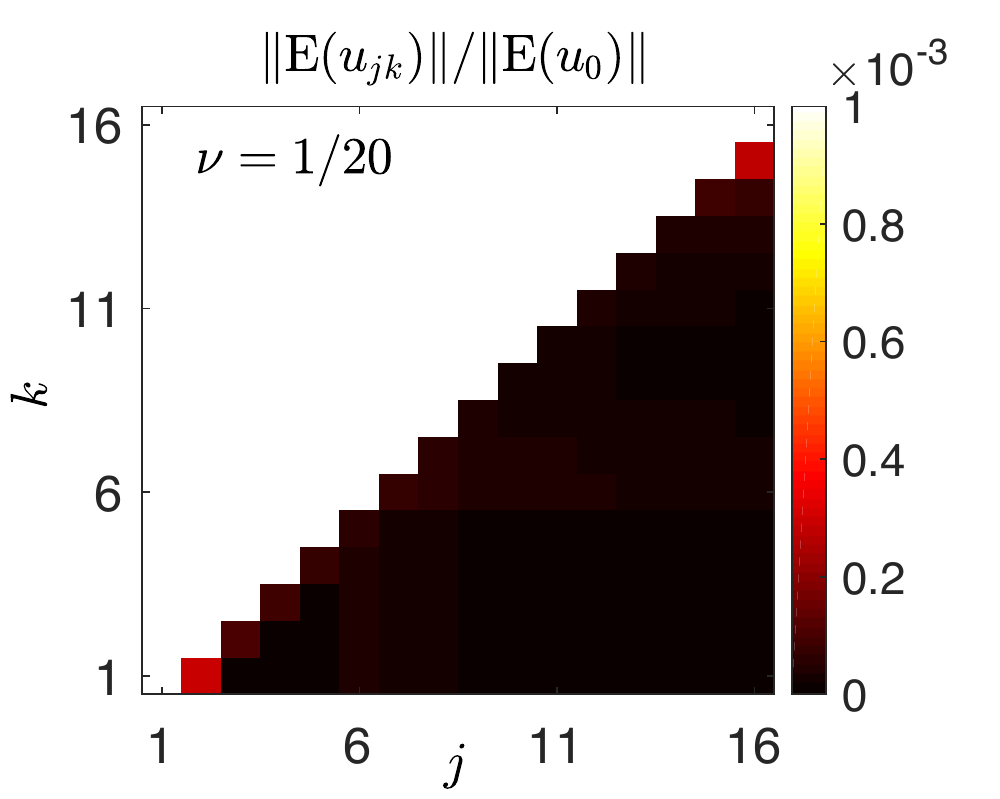}} 
   \caption{ Relative magnitudes of means of the first and second order ANOVA expansion terms $u_j(x,\xi_j)$ and $u_{jk}(x,\xi_j,\xi_k)$ in the case of $1\times 16$ subdomain, comparing $\nu=1/2$ and $\nu=1/20$. The solution for the example with smaller diffusion parameter $\nu=1/20$ is more anisotropic but has smaller variance. 
   }
   \label{fig:cANOVA_dim1x16_all}
\end{figure}

To explore accuracy 
with respect to parameterization, 
we compute a reference solution 
using a quasi Monte-Carlo (MC) method with 
$10^5$ samples of a Halton sequence, and denote  
{the computed mean and standard deviation as 
$\E[u^{mc}]$ and $\sigma[u^{mc}]$, respectively.} 
We then estimate the relative $L_2$-error 
of reduced basis solutions $u_{rb}$ 
in the mean and standard deviation as 
\begin{equation}
	e_{\mu} = \frac{\|\E[u^{mc}] - \E[u_{rb}] \|_2}{\| \E[u^{mc}] \|_2}, \qquad 
	e_{\sigma} = \frac{\|\sigma[u^{mc}] - \sigma[u_{rb}] \|_2}{\| \sigma[u^{mc}] \|_2}. 
\end{equation} 
%

%
\begin{figure}
    \centerline{ 
    \includegraphics[width=14cm]{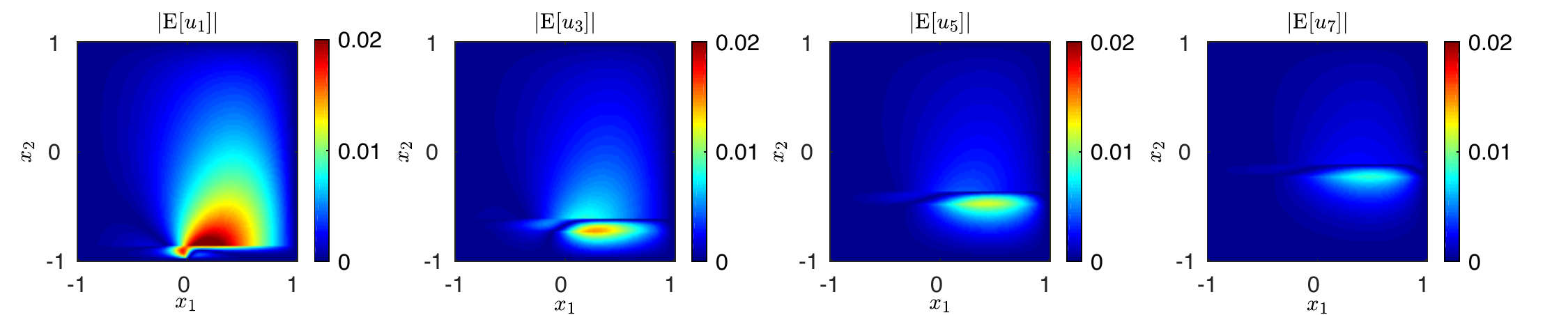}}  
    \centerline{ 
    \includegraphics[width=14cm]{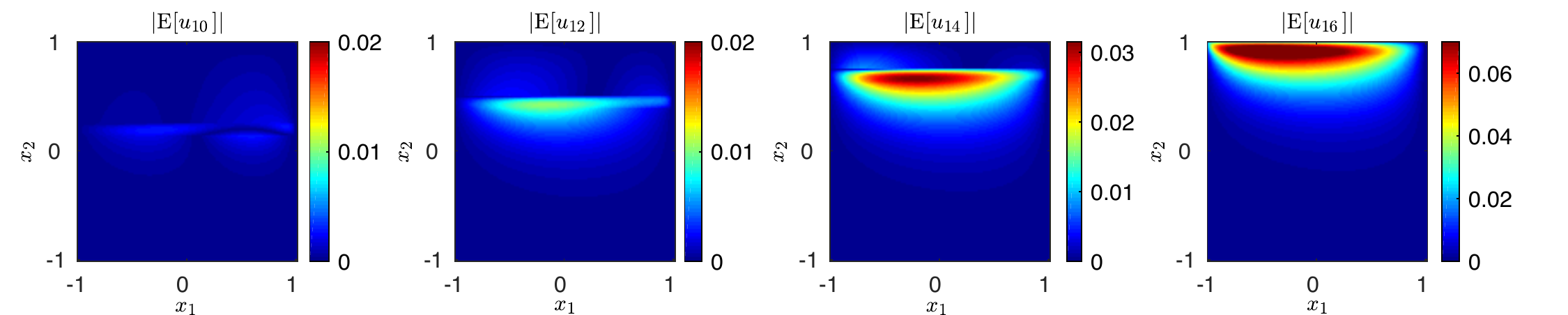}}    
    \caption{ Absolute values of the means of the first order ANOVA expansion terms $|\E[u_j(x,\xi_j)]|$ 
    for $\nu=1/2$ on $1 \times 16$ partition. The ANOVA terms with $\xi_j$ attached close to the top and bottom discontinuity have large values. } \label{fig:cANOVA_acd_dim1x16}
\end{figure}
\begin{figure}
    \centerline{ 
    \includegraphics[width=14cm]{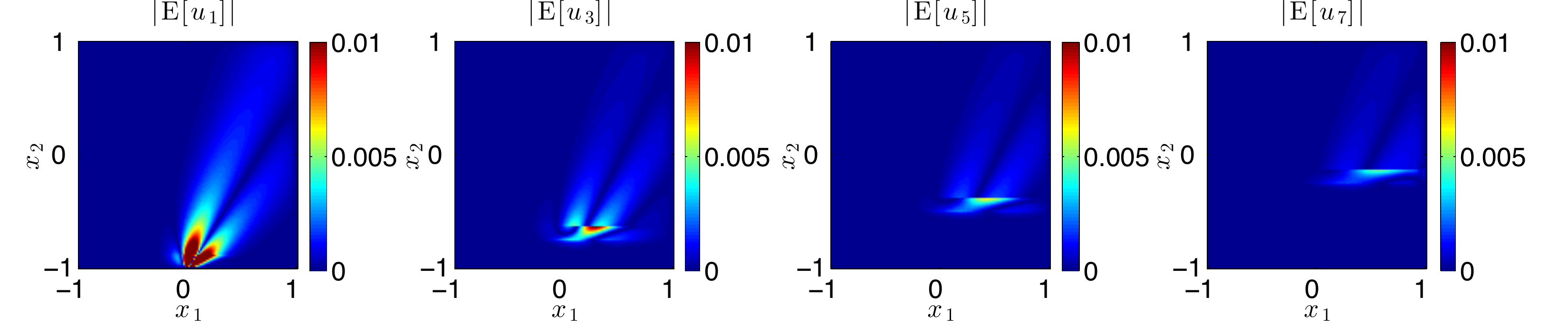}}  
    \centerline{ 
    \includegraphics[width=14cm]{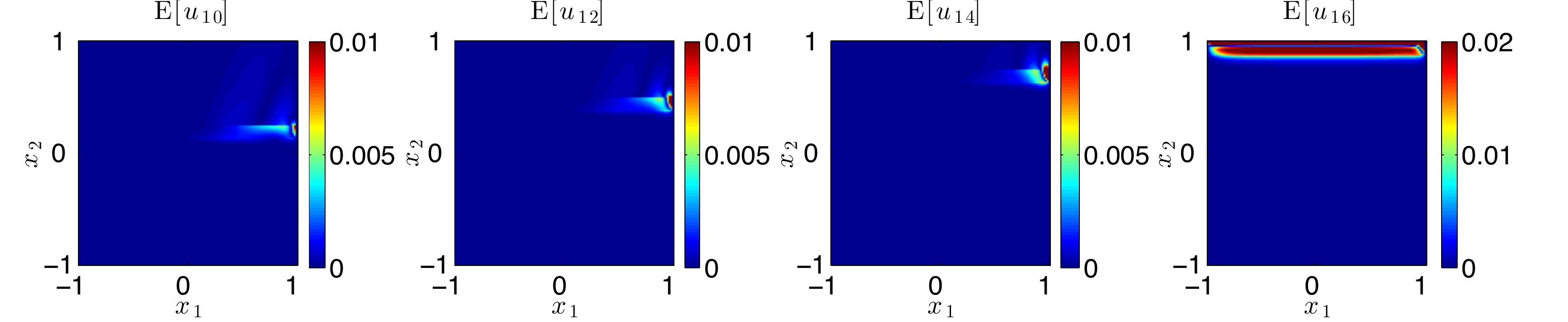}}  
    \caption{ Absolute values of the means of the first order ANOVA expansion terms $|\E[u_j(x,\xi_j)]|$ 
    for $\nu=1/20$ on $1 \times 16$ partition. The effect of the noise from each subdomain is more localized than the case of $\nu=1/2$. 
    } \label{fig:cANOVA_acd_dim1x16_}
\end{figure} 

\begin{figure}
    \centerline{ 
    \includegraphics[width=14cm]{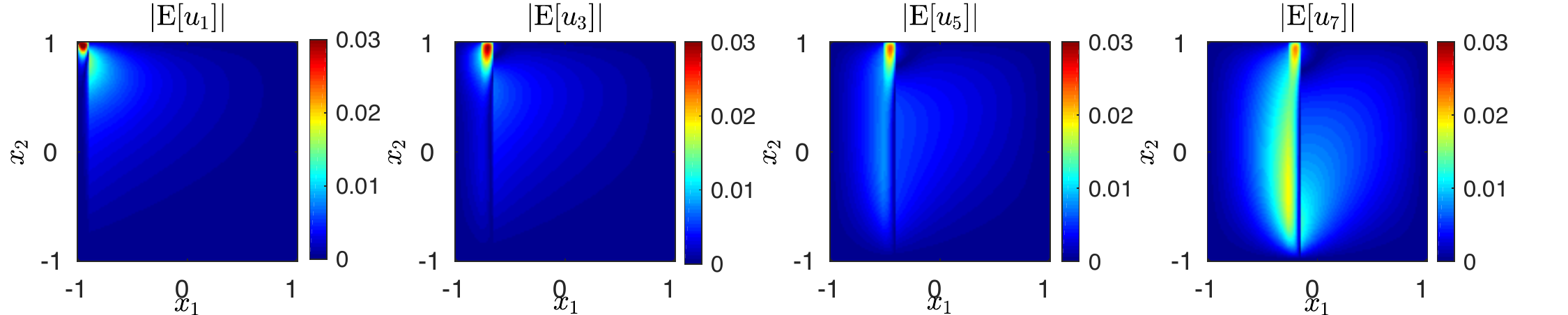}}  
    \centerline{ 
    \includegraphics[width=14cm]{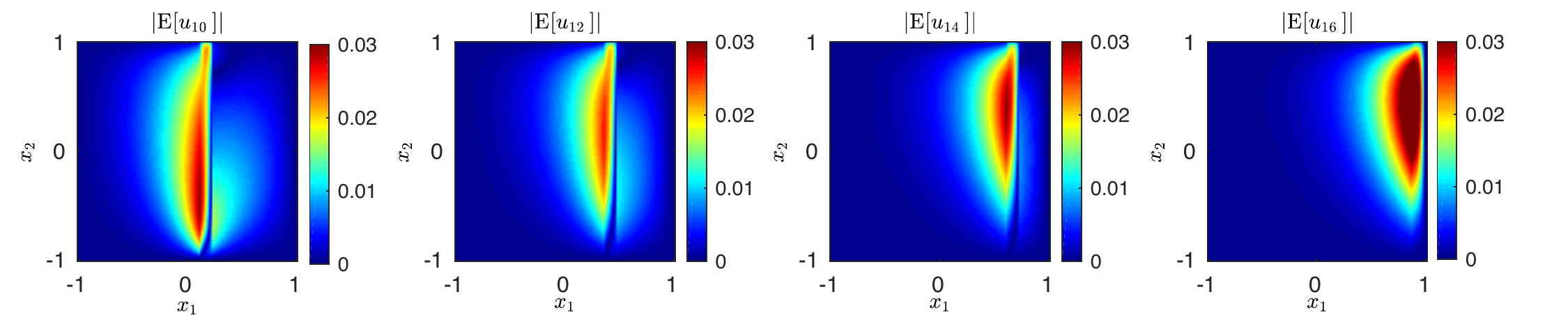}}    
    \caption{ {Absolute values of the means of the first order ANOVA expansion terms $|\E[u_j(x,\xi_j)]|$ 
    for $\nu=1/2$ on $16 \times 1$ partition. The magnitudes of the ANOVA terms are less anisotropic than the case of $1\times 16$ since all of the directions are attached to the top discontinuity.} } \label{fig:cANOVA_acd_dim16x1}
\end{figure}
\begin{figure}
    \centerline{ 
    \includegraphics[width=14cm]{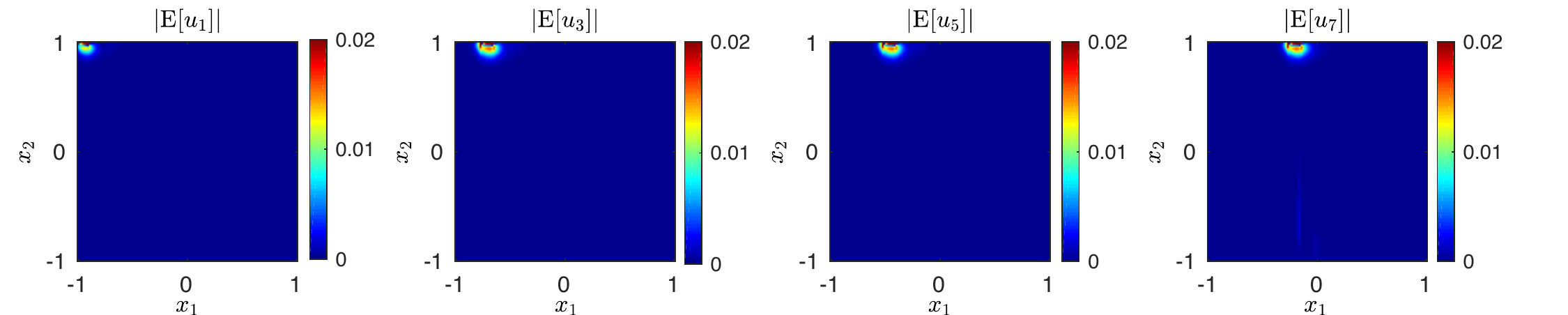}}  
    \centerline{ 
    \includegraphics[width=14cm]{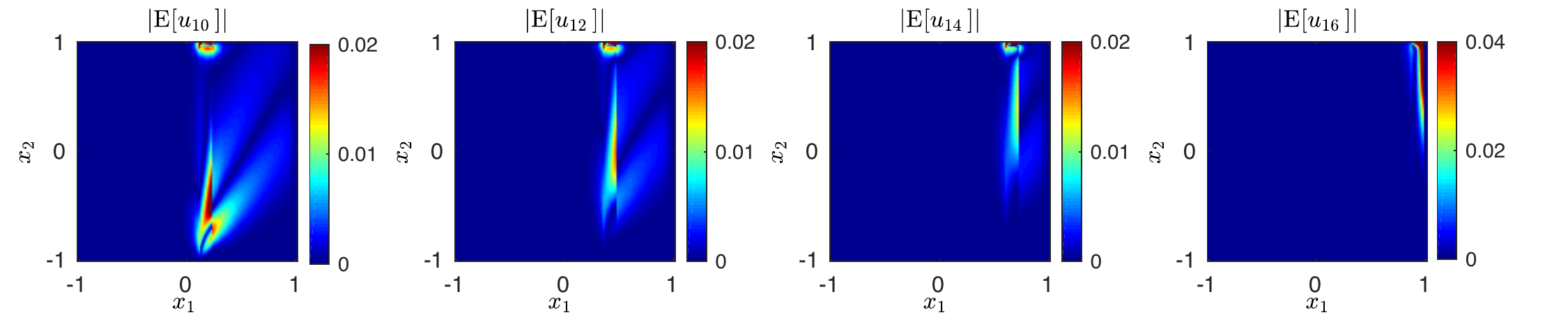}}  
    \caption{ {Absolute values of the means of the first order ANOVA expansion terms $|\E[u_j(x,\xi_j)]|$ for $\nu=1/20$ on $16 \times 1$ partition. The effect of the noise from each subdomain is more localized close to the top and bottom discontinuity than the case of $\nu=1/2$.} } 
    \label{fig:cANOVA_acd_dim16x1_}
\end{figure} 

We first present the results where the ANOVA decomposition 
helps us 
to understand the impact of 
anisotropy on features of the solution. 
We use the reduced basis tolerance $\epsilon_{RB} = 10^{-4}$ 
(see Algorithm \ref{Alg:RBM-}) 
and the full ANOVA decomposition of the solution  
without truncation by setting $\epsilon_A = 10^{-15}$ 
( see \eqref{eq:EffecitveDim}). 
Figure \ref{fig:cANOVA_dim1x16_all} plots 
magnitudes of the means of the first and second-order 
ANOVA terms, $\|\E[u_j(x,\xi_j)]\|$ and 
$\|\E[u_{jk}(x,\xi_j, \xi_k)]\|$, 
for a $1\times 16$ partition 
(see Figure \ref{fig:subD_rDim16}). 
%
%
The image on the left 
shows that the ANOVA terms for components of $\xi_j$ 
near the top boundary (large $j$), 
where there is a boundary layer, 
and near the bottom (small $j$), 
where there is a discontinuity 
along the inflow boundary, are large. 
Among the second order terms, the pairs involving 
the aforementioned directions 
and the terms with physically adjacent subdomains, 
e.g, $u_{j,j+1}(x,\xi_j,\xi_{j+1})$, 
have large magnitudes as well. 
These points are further substantiated 
in Figure \ref{fig:cANOVA_acd_dim1x16}, 
where the absolute values of the means 
$|\E[u_j(x,\xi_j)]|$ are plotted for several $j$. 
In particular, 
there is significant activity near discontinuities, 
$u_{16}$, $u_{14}$, and $u_{1}$. 
This is also true in the case of smaller diffusion parameter $\nu=1/20$ 
plotted in Figure \ref{fig:cANOVA_acd_dim1x16_}. 
In addition, the system is more anisotropic compared to $\nu=1/2$ 
and the contribution is more localized 
in each subdomain. 
We emphasize that this feature will make 
the reduced basis algorithm based on ANOVA more efficient, 
since the existence of 
relatively important parameters indicates 
a possible truncation in the ANOVA decomposition. 
In constrast, Figures 
\ref{fig:cANOVA_acd_dim16x1} and \ref{fig:cANOVA_acd_dim16x1_} 
show results 
for a $16\times 1$ partition. 
The maximum values of the ANOVA terms are all of the same order, 
in other words, the system is less anisotropic 
with respect to the parameters 
than the case of $1\times 16$, 
since all subdomains 
are attached to the top boundary layer. 
This is more apparent when $\nu=1/20$ 
in Figure \ref{fig:cANOVA_acd_dim16x1_}, 
where the effect  of the noise from each 
subdomain is more localized. 
%
%

\begin{figure}
    \centerline{ {\footnotesize \hspace{0.5cm}$ 1\times 16$ \hspace{3.2cm} $ 16\times 1$ } }
    \centerline{ \rotatebox{90}{\hspace{1.2cm} mean }
    \includegraphics[width=4.0cm]{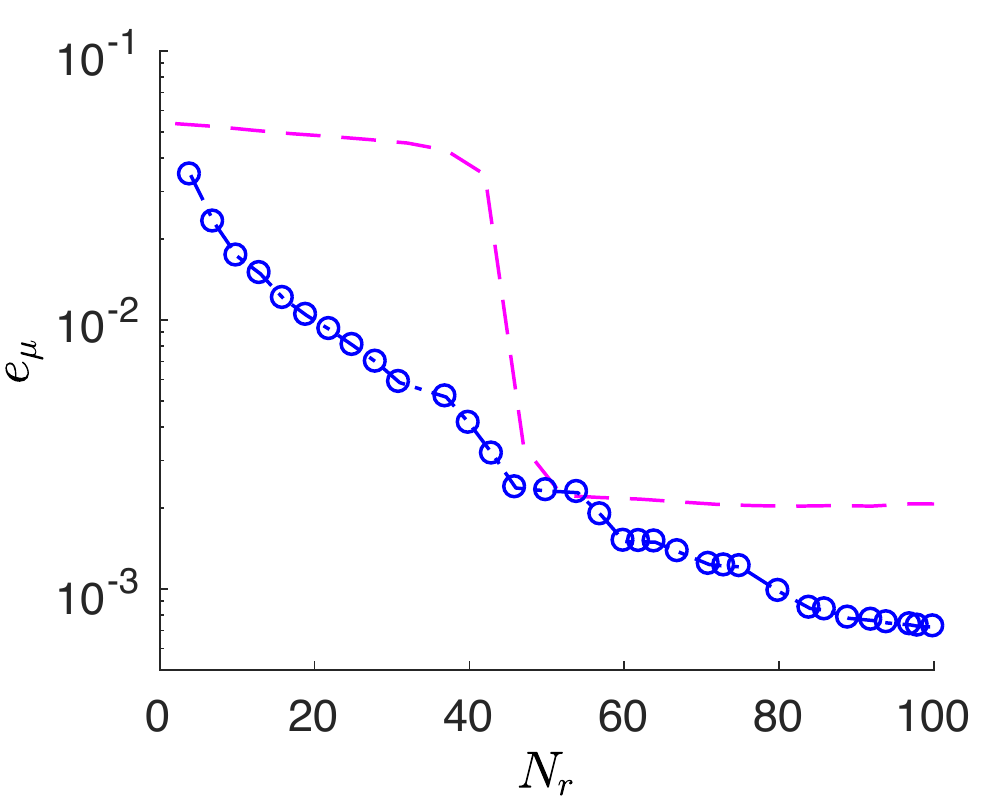}
    \includegraphics[width=4.0cm]{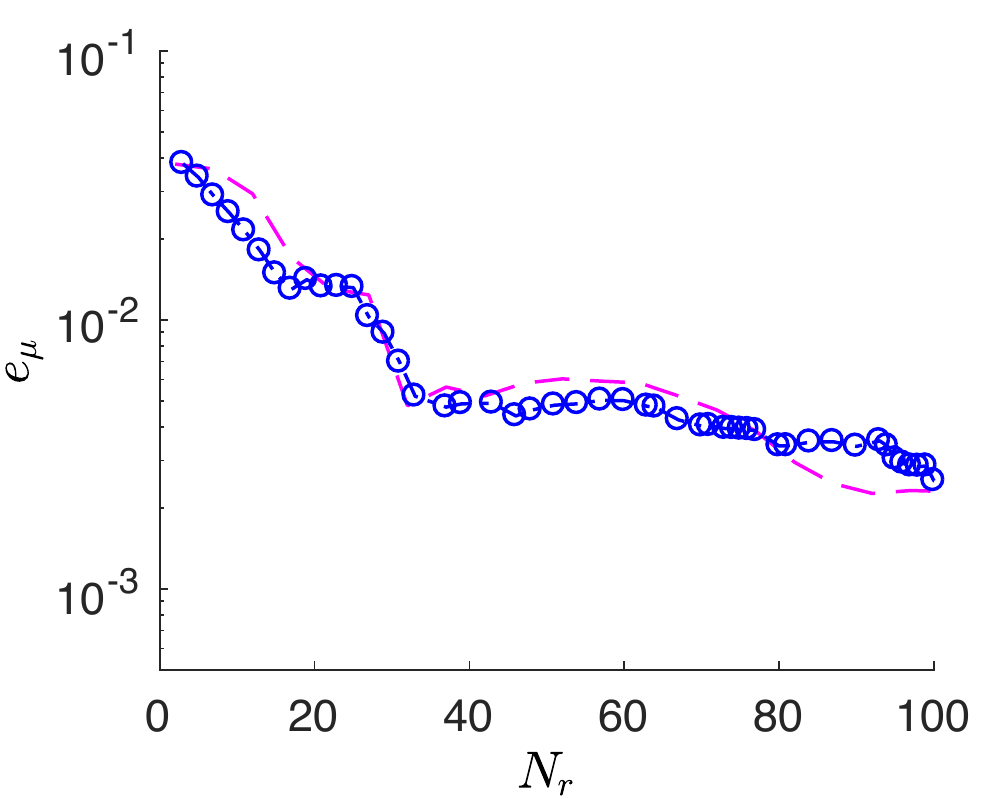}}
    \centerline{ \rotatebox{90}{\hspace{1.5cm} s.d.}
    \includegraphics[width=4.0cm]{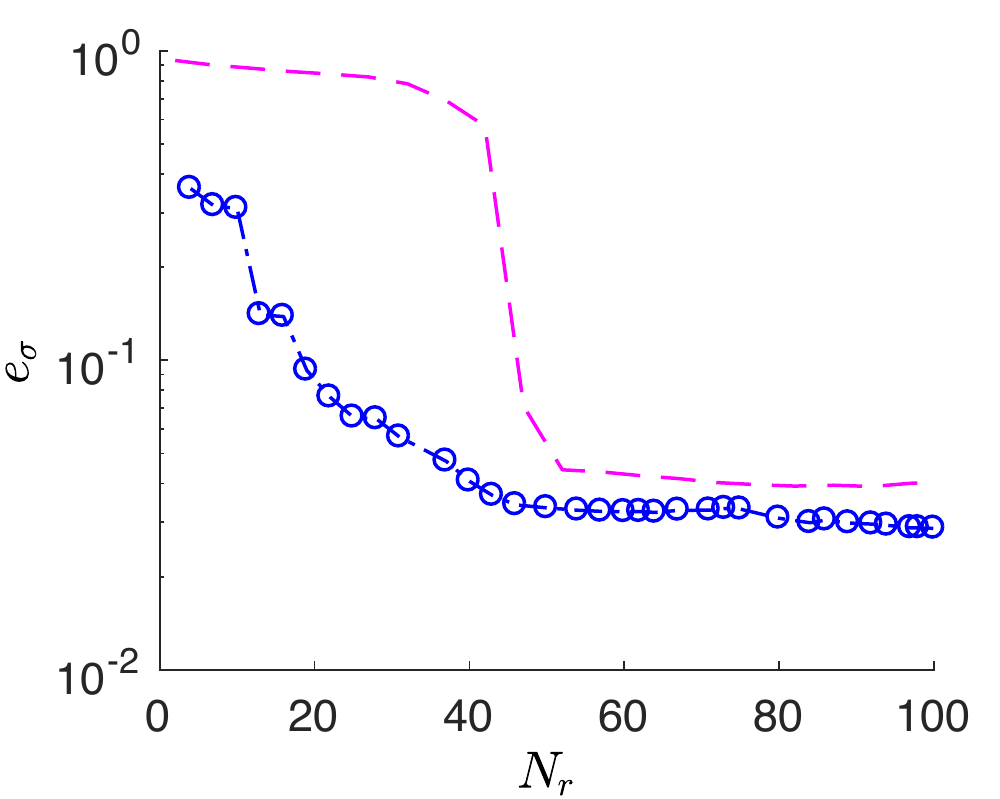}
    \includegraphics[width=4.0cm]{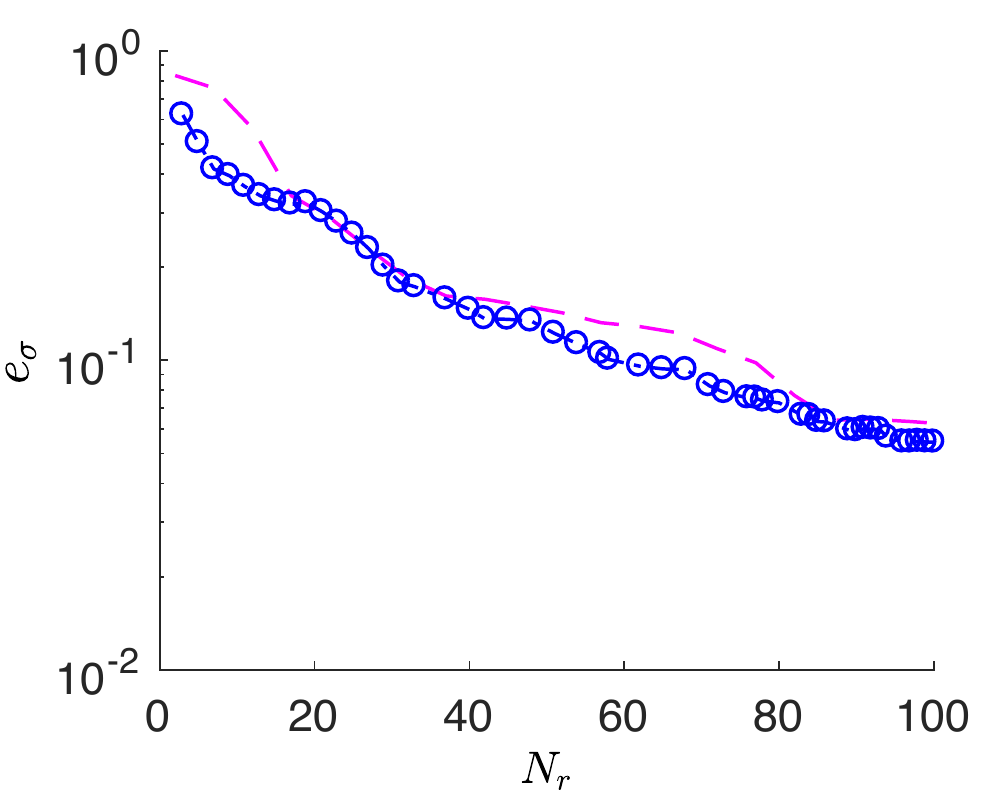}}
    \centerline{ 
    \includegraphics[width=7.5cm]{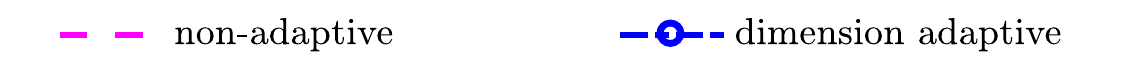}} 
   \caption{Error in the mean and standard deviation of the solution computed using reduced basis of size $N_r$. Shown results compare the reduced basis algorithm with and without the dimension-adaptive procedure for the problems on $1\times 16$ and $16\times 1$ partitions.  
   The dimension-adaptive procedure using ANOVA decomposition has an advantage for anisotropic problems, e.g., $1\times 16$. } 
	\label{fig:errtbt_rDim16}
\end{figure}

Next, we examine the quality of the reduced basis sorted 
using the ANOVA indicator $\gamma_K$ \eqref{eq:anovaIndc} 
as in Algorithm \ref{Alg:RBM-sort}.  
Figure \ref{fig:errtbt_rDim16} 
plots the relative $L_2$-error in the moments 
with respect to the size of the reduced basis $N_r$ 
for the problems on 
$1\times 16$ and $16\times 1$  partitions 
with $\nu = 1/2$. 
The errors are computed 
every two indices up to 100. 
The pink dashed lines 
show the error using the sparse-grid points of level two 
without the dimension-adaptive procedure, 
typically, generated in increasing order of index, 
e.g., $\xi_1$, $\xi_2$, $\cdots$, 
and $\xi_{12}$, $\xi_{13}$, $\cdots$.  
The initial saturation of the error 
in the $1\times 16$ 
case is due to this ordering, 
and the error drops suddenly 
after obtaining the reduced basis 
that has 
a large contribution to the solution, 
for instance, when the points using 
index $\xi_{16}$ are considered. 
%
However, 
our algorithm 
sorts the solution 
by the magnitude of $\gamma_K$ 
and the error effectively decays 
while avoiding an additional computation  
induced either by adopting the greedy procedure 
or sorting the reduced basis afterwards.  
This computation becomes particularly 
expensive for large $M$ and $N_r$, 
which makes our algorithm more preferable. 
This advantage is not apparent in the $16\times 1$ case 
shown in the second column, 
due to the fact that all directions 
are influenced in essentially the same way 
by the top boundary layer. 
However, our approach will show its effectiveness  
once the random variables $\{\xi_j\}$ have 
distinctive distributions.

\begin{table} 
\centering
\caption{ Number of full solves $N_r$ in the reduced basis collocation method using PCM-ANOVA points for $1\times 4$ partition and $\nu = 1/20$. } 
\label{Tbl:nRB_ANOVA}
\begin{tabular}{c|ccccc|cccc}
\hline 
$\ell$ & \multicolumn{5}{c|}{ $\ell = 3$ }  & 1 & 2 & 3 & 4  \\ \hline 
$p$ & 3 & 5 & 7 & 9 & 11 & \multicolumn{4}{c}{ $p = 9$ }  \\ \hline 
$\epsilon_{RB} = 10^{-3}$ &   3 & 4 & 4 & 4 & 4 & 4 & 4 & 4 & 4 \\
$\epsilon_{RB} = 10^{-4}$ &  16 & 24 & 31 & 35 & 37  & 15 & 25 & 35 & 38 \\
$\epsilon_{RB} = 10^{-5}$ &  33 & 63 & 80 & 90 & 97  &  25 & 74 & 90 & 103 \\
\hline 
\end{tabular}
\end{table}  
\begin{figure}
    \centerline{ {\footnotesize \hspace{.0cm} $\ell = 3$  \hspace{3.7cm} $p=9$ \hspace{1.7cm} } }
    \centerline{ \rotatebox{90}{\hspace{1.5cm} {\footnotesize  mean} }
    \includegraphics[width=4.5cm]{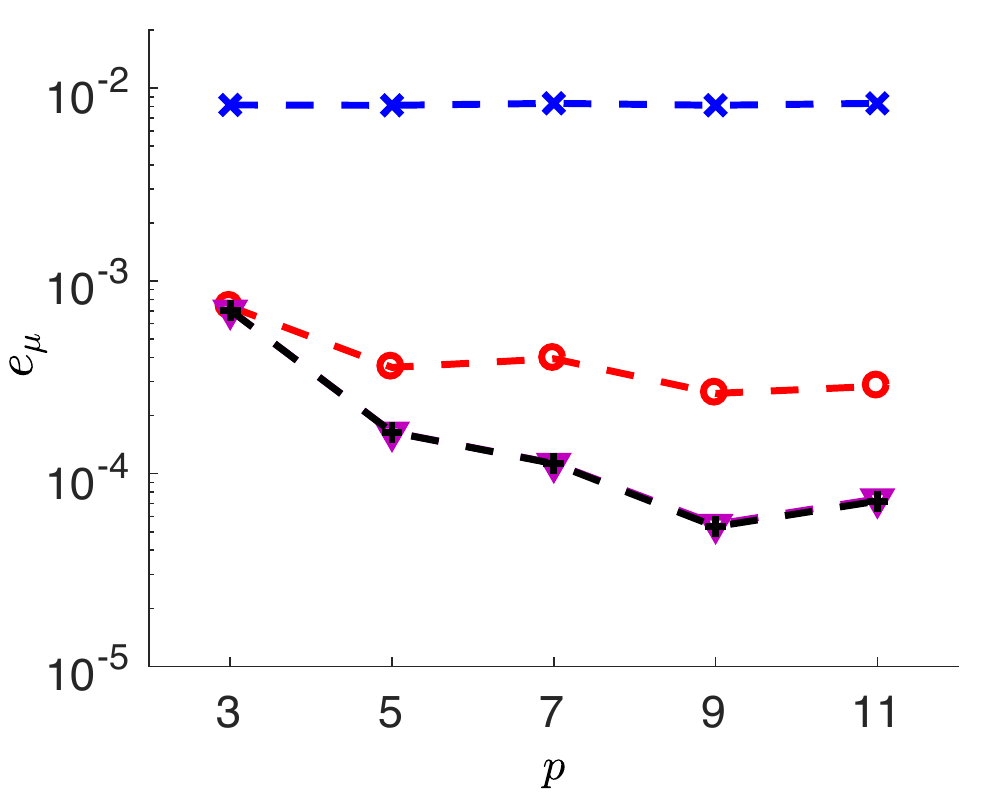} 
    \includegraphics[width=4.5cm]{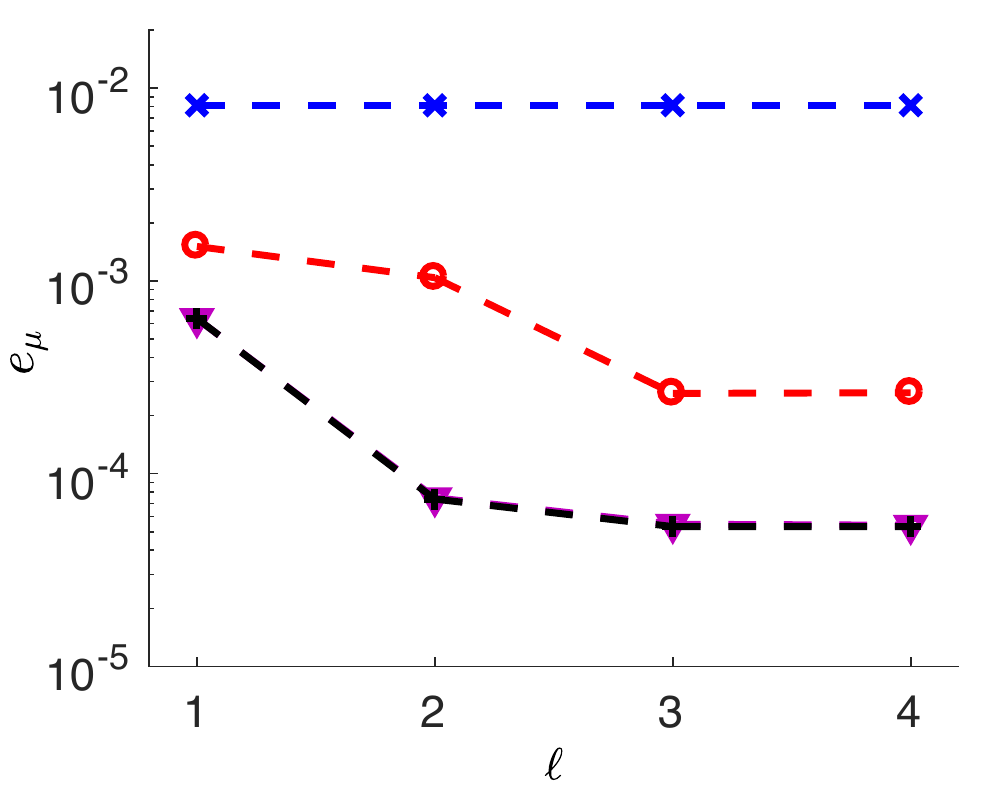} 
    \includegraphics[width=2.2cm]{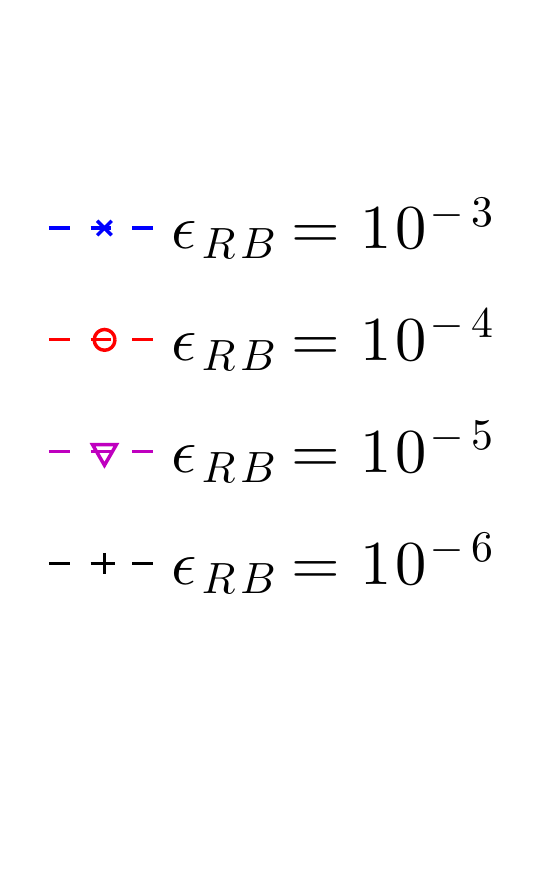} } 
   \caption{ Relative $L_2$-error of the mean 
   using the PCM-ANOVA collocation points with level $\ell$ and $p$ and reduced basis tolerance $\epsilon_{RB}$, for $1\times 4$ partition and $\nu = 1/20$. } \label{fig:RBM_ANOVA_nHnlev}
\end{figure}

We close this section by presenting 
a convergence study of the PCM-ANOVA points 
in Figure \ref{fig:RBM_ANOVA_nHnlev} and Table \ref{Tbl:nRB_ANOVA}, 
which motivates the necessity of 
adapting the resolution in each direction.  
Unlike the sparse-grid points, the ANOVA points 
have another independent parameter, 
the polynomial order $p$, in addition to the level $\ell$. 
Figure \ref{fig:RBM_ANOVA_nHnlev} shows the relative $L_2$-error 
of the mean with respect to these two parameters 
for a $1\times 4$  partition and $\nu = 1/20$. 
The results are computed by either fixing the level $\ell=3$ or 
the polynomial order $p=9$. Although in both cases, 
the errors decrease with respect to $p$ and $\ell$, 
they tend to stagnate after certain maximal 
values are reached. 
Thus, the polynomial order and the level should be increased 
together with a carefully selected tolerance 
regarding the system \cite{Foo1}.

\section{Adaptive algorithm for RBM-ANOVA} 
\label{sec:RBM-ANOVAp}

In this section, we introduce an adaptive reduced basis algorithm 
by increasing the number of collocation points used to compute each ANOVA term, 
in other words, 
the polynomial order $p$ in each effective dimension 
$\J_{\ell}$ \eqref{eq:EffecitveDim}. 
As shown in Figure \ref{fig:RBM_ANOVA_nHnlev} and also later 
in Figure 
\ref{fig:errRBmmt_dim16_A}, fixed $p$ may limit accuracy. 
However, it is not practical to choose large $p$ from the beginning,  
since that leads to a large number of collocation points. 
Thus, we propose an adaptive PCM-ANOVA reduced basis algorithm 
that increases $p$ only among the active dimensions $\J_{\ell}$ 
and with an appropriate criterion that limits the increases 
in $p$.

The idea is shown in Algorithm \ref{Alg:RBM-ANOVA_adaptp}. 
The method initializes the collocation points with polynomial 
order $p_0$, typically $p_0 = 3$, and level $\ell_0$. 
This choice will make the cost in the initial step 
to be smaller than 
that using a sparse-grid with the same level $\ell_0$. 
%
Then, among the effective dimensions $\J$ computed by 
the adaptive ANOVA criterion $\gamma_K$ \eqref{eq:anovaIndc}, 
we increase $p$ in specific directions until a 
termination criterion holds. %
For this, we use a saturation condition that measures 
the relative difference of ANOVA terms in direction $K$ 
after increasing the polynomial order, that is, 
\begin{equation} 
	\rho_K^{p_K} \doteq \frac{ \|\E[ \su_K^{\bar{p}_{K}} ] - \E[ \su_K^{p_K} ] \| }{ \| \sum_{|S|\leq |K|} \E[ \su_S ] \| } < \epsilon_p, 
	\label{eq:AdaptpCondi}
\end{equation} 
where $\su_K^{p_K}$ is the ANOVA term in direction $K$ 
computed with the increased polynomial order $p_K$, 
and $\su_K^{\bar{p}_K}$ is from the previous iteration, with 
$\bar{p}_{K}<{p}_{K}$. 
If this quantity is smaller than a prescribed tolerance 
$\epsilon_p = O(\epsilon_{RB})$, 
the reduced basis solution 
added in the current iteration 
is excluded from $\bS_{N_r}$ 
and the direction $K$ is excluded 
from the set of effective dimensions $\J$. 
This condition is based on 
the error estimate of the stochastic spectral solution 
between large polynomial orders \cite{XWGK_JCP05}. 
In addition, if all the 
error indicators \eqref{eq:rbmindc} 
computed from 
the collocation points in the direction $K$ 
are less than the 
reduced basis tolerance $\epsilon_{RB}$, 
we exclude $K$ from $\J$. 
%
A stronger terminating criterion can be imposed by 
restricting the maximum order of polynomial 
at level $K$, $|K|=l$, from the previous level 
as $p_{K} \leq \max_{|S|=l-1} p_{S}$. 
 
\begin{algorithm} 
\caption{Adaptive-RBM PCM-ANOVA( $\{\Xi_{l}^{p}\}_{l=0}^{\ell}$, $\epsilon_{RB}$, $\epsilon_{A}$, $\epsilon_{p}$ ) : Reduced basis method based on dimension and $p$-adaptive PCM-ANOVA collocation points up to level ${\ell}$. The parameters are typically chosen as $\ell_0 = 1$ or $2$ and $p_0=3$. } 

$\,\,\,$1: Compute the solution $\bu_{\xi^{(0)}}$ at the anchor point $\xi^0 \in {\Xi_0}$. \\ 
$\,\,\,$2: Initialize $\bS_{\Theta_r} = \bu_{\xi^{(0)}}$ and $\bQ = \bu_{\xi^{(0)}}/ \|\bu_{\xi^{(0)}}\|^2 $ and 
construct the offline reduced matrices and vectors. \\
$\,\,\,$3: Initialize active dimensions $\J = \cup_{l=1}^{\ell_0} \J_l$, 
 $\J_l = \{ K \subset \{ 1,..., M \} \, | \, |K| = l \}$, 
and $p_K = p_0$ for all $K \in \J$. \\ 
$\,\,\,$4:  \textbf{for} $\,$ $ \ell_0 \leq l \leq \ell$  \textbf{do}  \\ 
$\,\,\,$5: $\quad$  \textbf{while} $\,$ $\J \neq \emptyset$   \textbf{do}  \\ 
$\,\,\,$6: $\qquad$  Store $\bar{\J} = \J$.  \\		
$\,\,\,$7: $\qquad$ \textbf{for} $K \in \J $  \textbf{do}  \\
$\,\,\,$8: $\qquad\quad$ Store $\bar{\su} = \su_K$ and  $\bar{\bQ} = \bQ$.  \\
$\,\,\,$9:  $\qquad\quad$  Update $\bQ$ in direction $K$, $(\bQ, \bS_{K})$ = RBM($\bQ$, $\Xi_K^{p_K}$, $\epsilon_{RB}$) in Algorithm \ref{Alg:RBM-}.  \\
10:  $\qquad\quad$  Compute the ANOVA term $\su_K$ and indicator $\gamma_K$ in \eqref{eq:anovaIndc}.  \\
11:  $\qquad\quad$   \textbf{if} $\bQ = \bar{\bQ}$ or $\gamma_K < \epsilon_A$ or  $\rho_K^{p_K}  < \epsilon_p$ in \eqref{eq:AdaptpCondi} \textbf{then}  \\ 
12: $\qquad\qquad$  Exclude term $K$, i.e., $\J = \J \setminus K$, and reset $\su_K = \bar{\su}$, $\bQ=\bar{\bQ}$. \\ 
13: $\qquad\quad$  \textbf{else} \\  
14: $\qquad\qquad$ Increase polynomial order $p_K$. \\ 
15: $\qquad\quad$  \textbf{end} \\ 	
16: $\quad\quad$  \textbf{end} \\ 
17: $\quad\quad$   Update $\bS_{\Theta_r} = [\bS_{\Theta_r}, \bS]$, where 	$(\bS,\J_l)$ = RB-Sort($\{ \gamma_K \}_{K\in\bar{\J}}$, $\{ \bS_{K} \}_{K\in\bar{\J}}$, $\bar{\J}$) in Algorithm \ref{Alg:RBM-sort}.  \\ 
18: $\quad$ \textbf{end} \\ 
19: $\quad$ 
Select active dimensions $\wJ_{l} = \{ K \in \J_l \,|\, |K| = l,\, \gamma_K > \epsilon_{A} \}$.  \\
20: $\quad$ Set next level indices $\J_{l+1} = \{ |K| = l+1 \,| \,  K = \cup_{T \in \wJ_l} T \} $ and $\J = \J_{l+1}$. \\
21: \textbf{end} \\ 
\vspace{0.1cm}
\label{Alg:RBM-ANOVA_adaptp} 
\end{algorithm}

\subsection{Numerical results}

\begin{figure}
    \centerline{ {\footnotesize $1\times 16$\hspace{2.8cm}$\qquad 4\times 4$\hspace{2.8cm}$\qquad 16\times 1$} }
    \centerline{ \rotatebox{90}{\hspace{1.4cm} mean}
    \includegraphics[width=4.0cm]{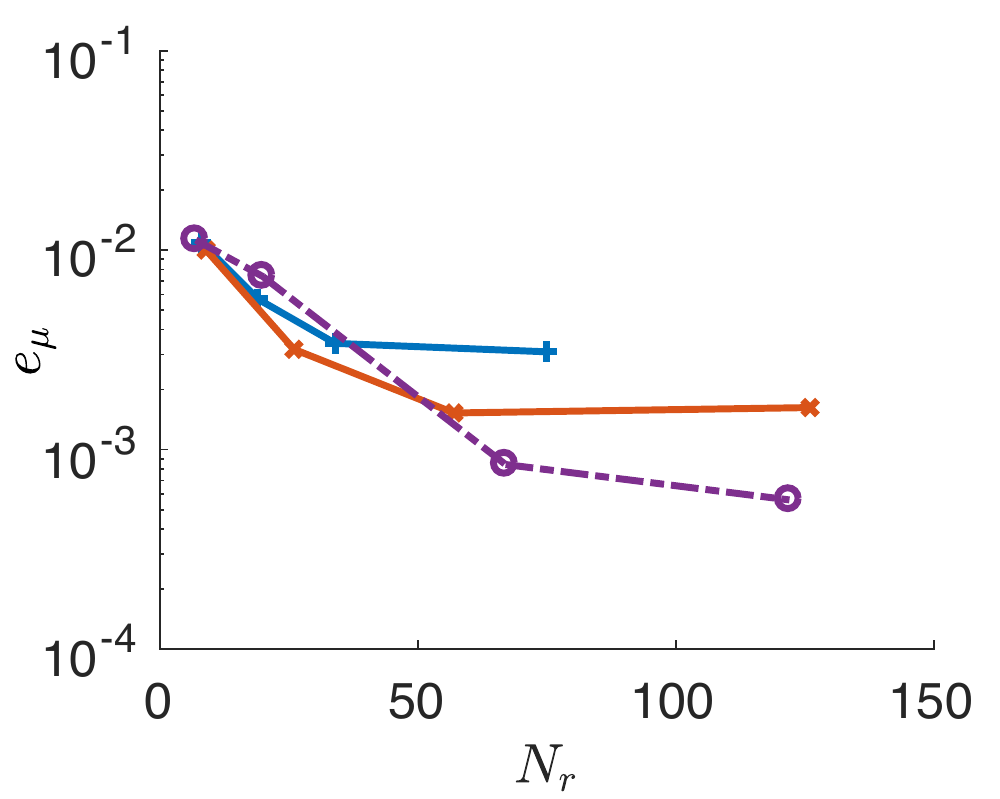}
    \includegraphics[width=4.0cm]{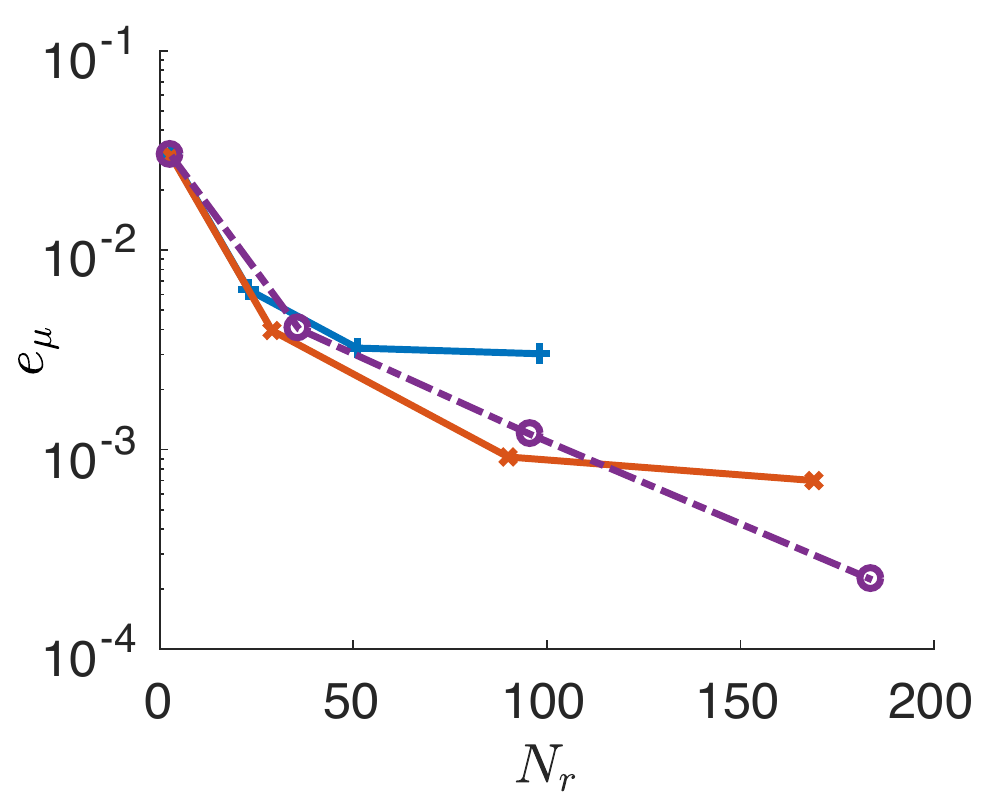}
    \includegraphics[width=4.0cm]{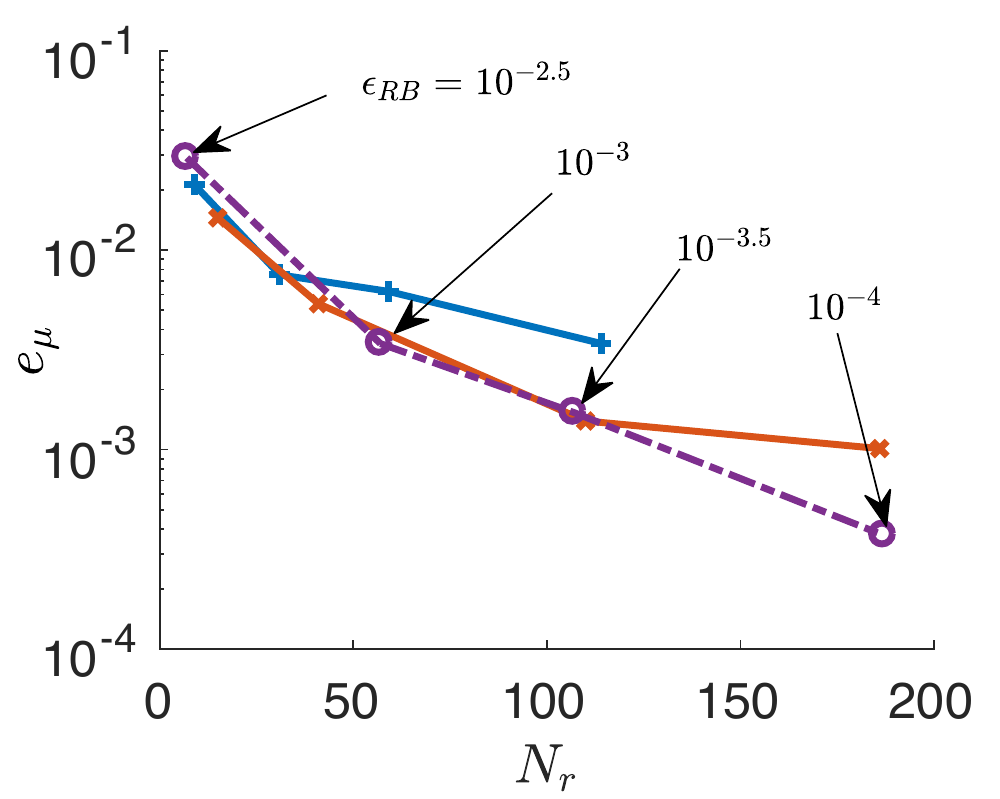}} 
    \centerline{ \rotatebox{90}{\hspace{1.4cm} s.d.}
    \includegraphics[width=4.0cm]{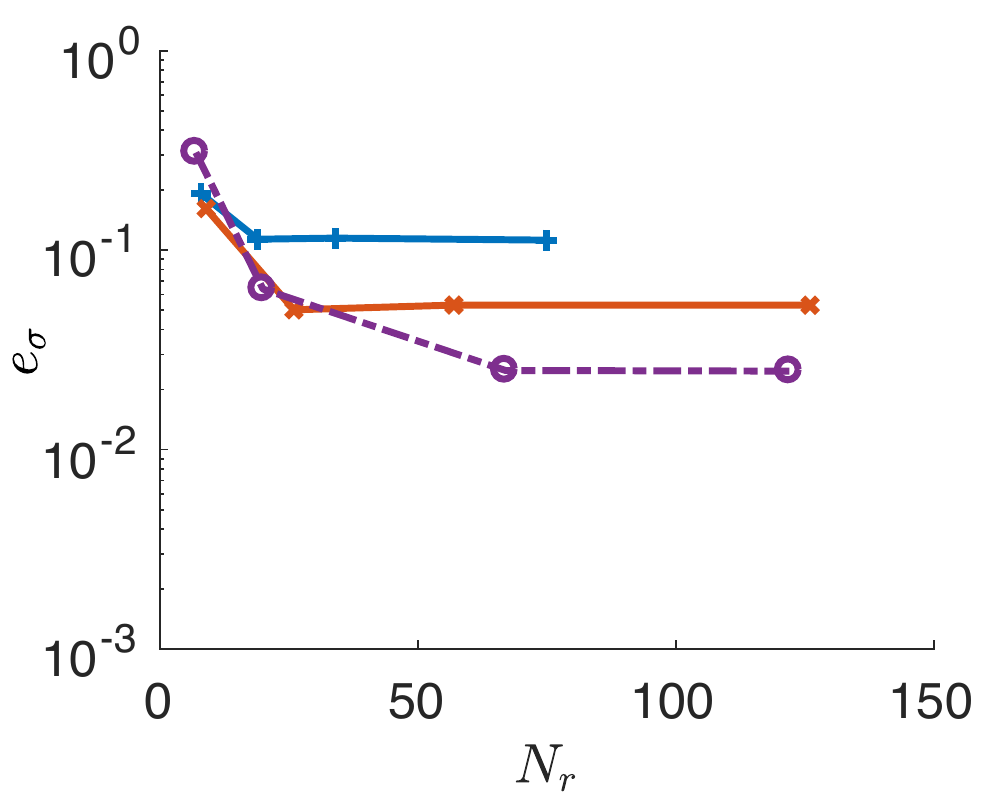}
    \includegraphics[width=4.0cm]{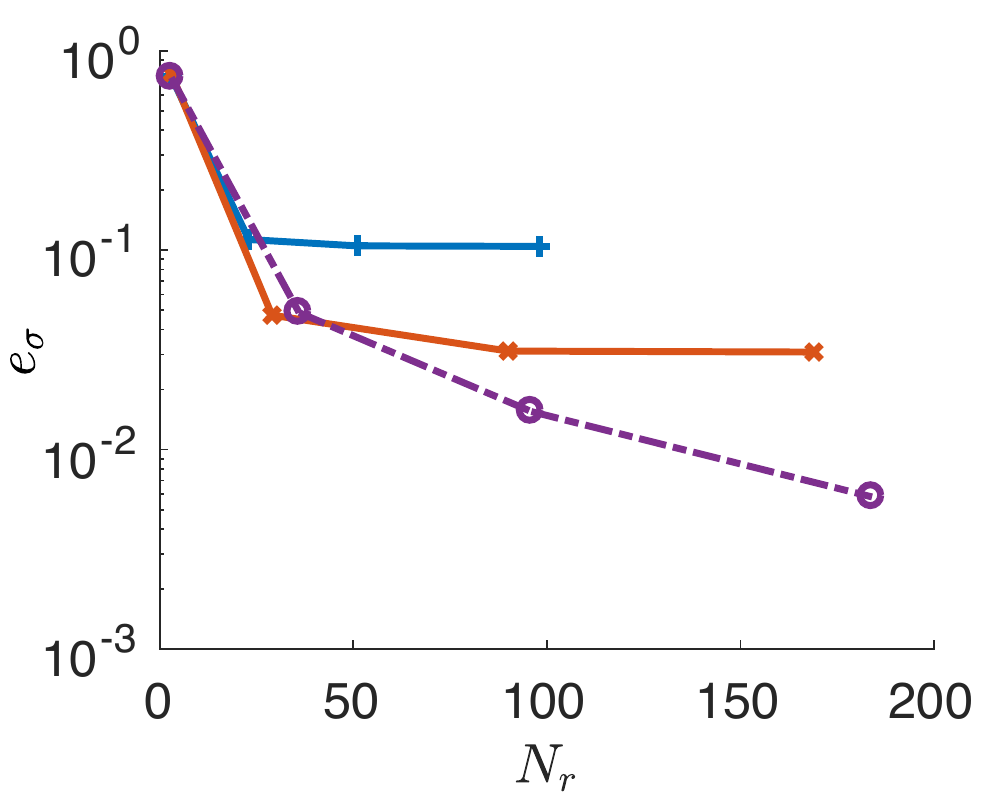}
    \includegraphics[width=4.0cm]{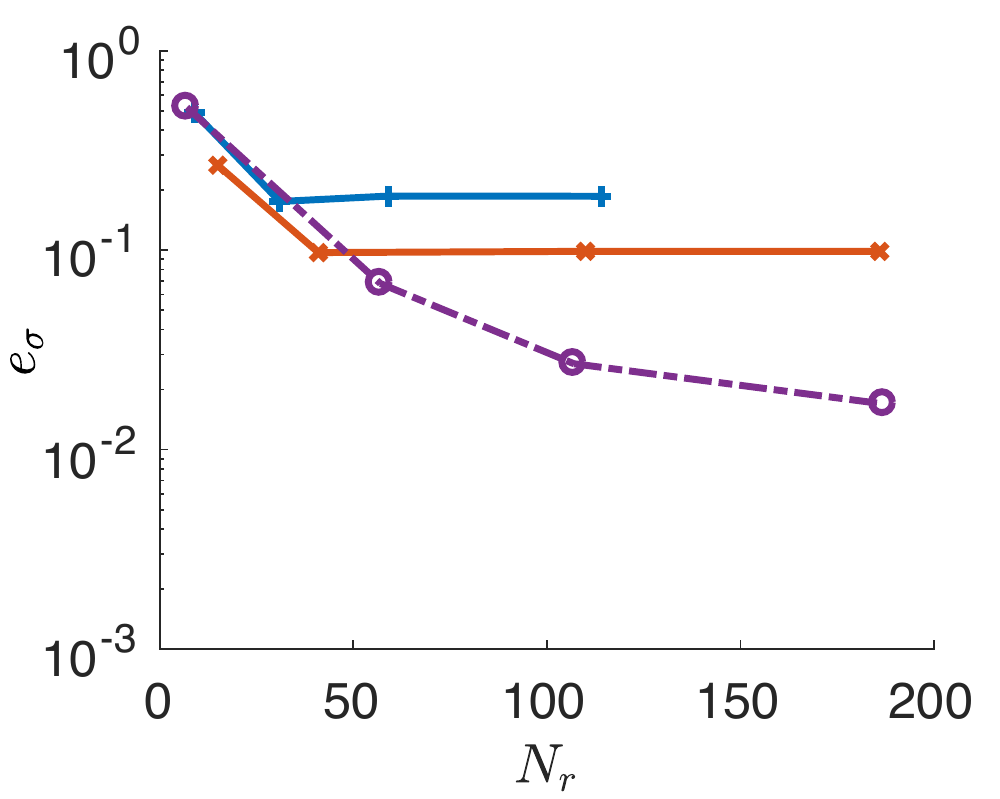}} 
    \centerline{ 
    \includegraphics[width=11.5cm]{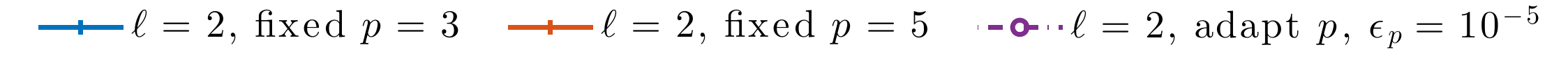}} 
   \caption{Error in the mean and standard deviation versus the size of reduced basis $N_r$ for the case of diffusion parameter $\nu = 1/2$ and $N_D=16$. The solution is computed by RBM using PCM-ANOVA with fixed $p$ or an adaptive $p$ with $\epsilon_p=10^{-5}$. The bullets correspond to varying $\epsilon_{RB}$ from $10^{-2.5}$ to $10^{-4}$. }
	\label{fig:errRBmmt_dim16_A}
\end{figure}

\begin{figure}
    \centerline{ {\footnotesize $1\times 16$\hspace{2.8cm}$\qquad 4\times 4$\hspace{2.8cm}$\qquad 16\times 1$} } 
    \centerline{ \rotatebox{90}{\hspace{1.9cm} mean}
    \includegraphics[width=4.0cm]{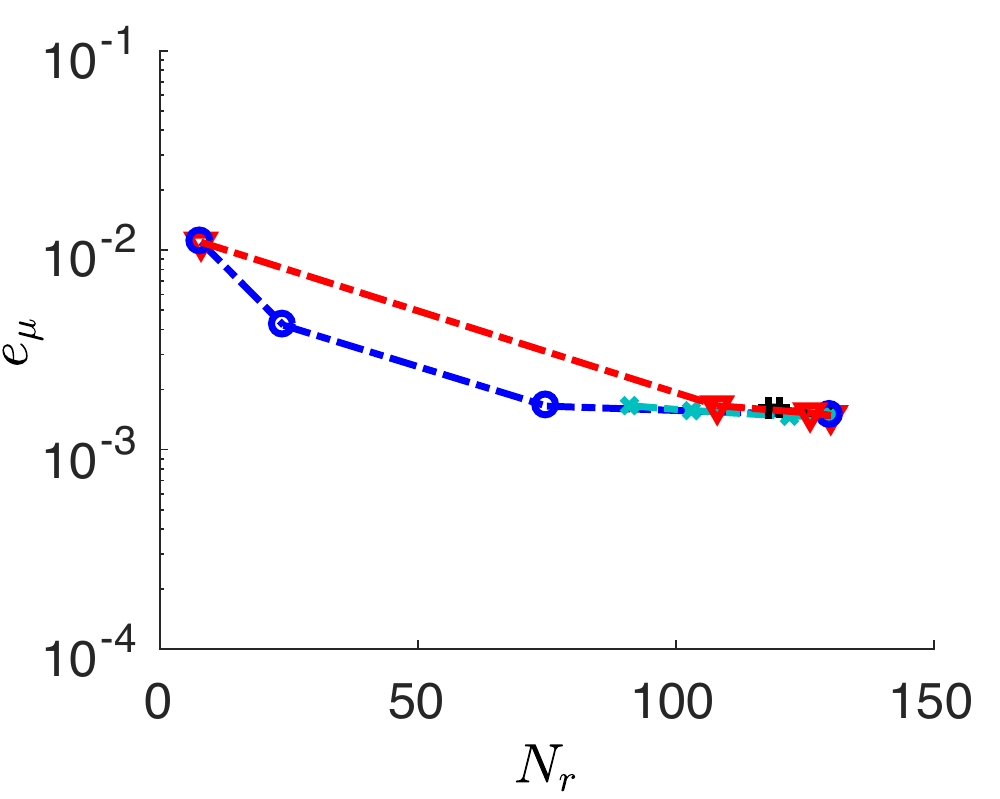}
    \includegraphics[width=4.0cm]{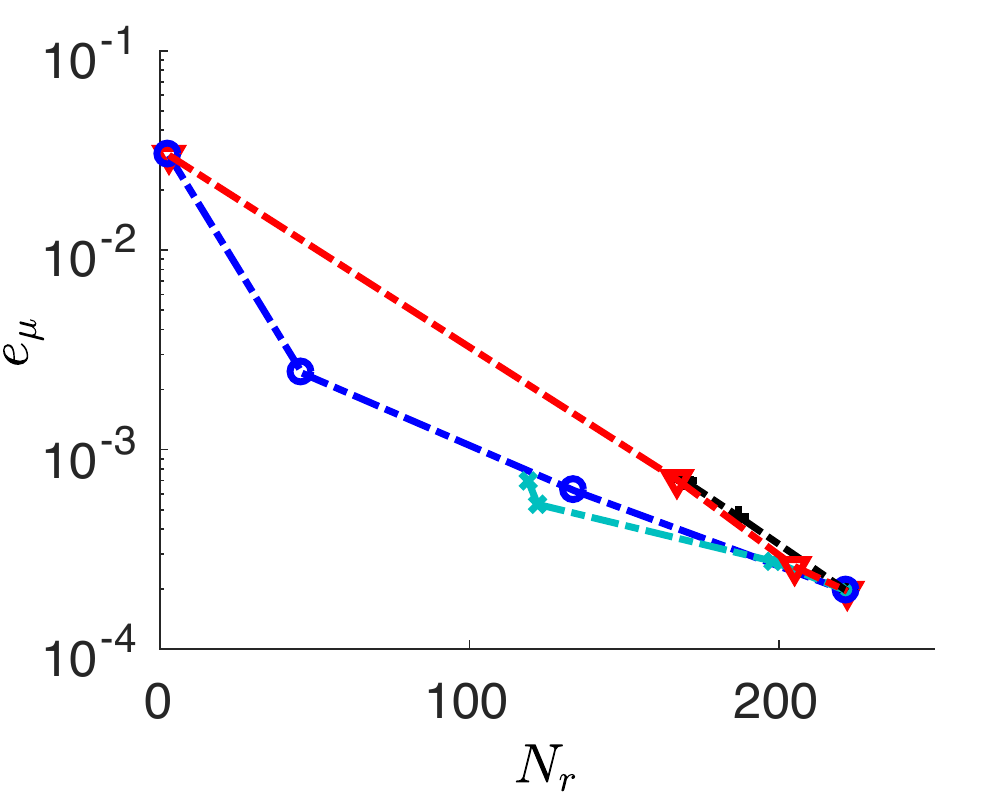}
    \includegraphics[width=4.0cm]{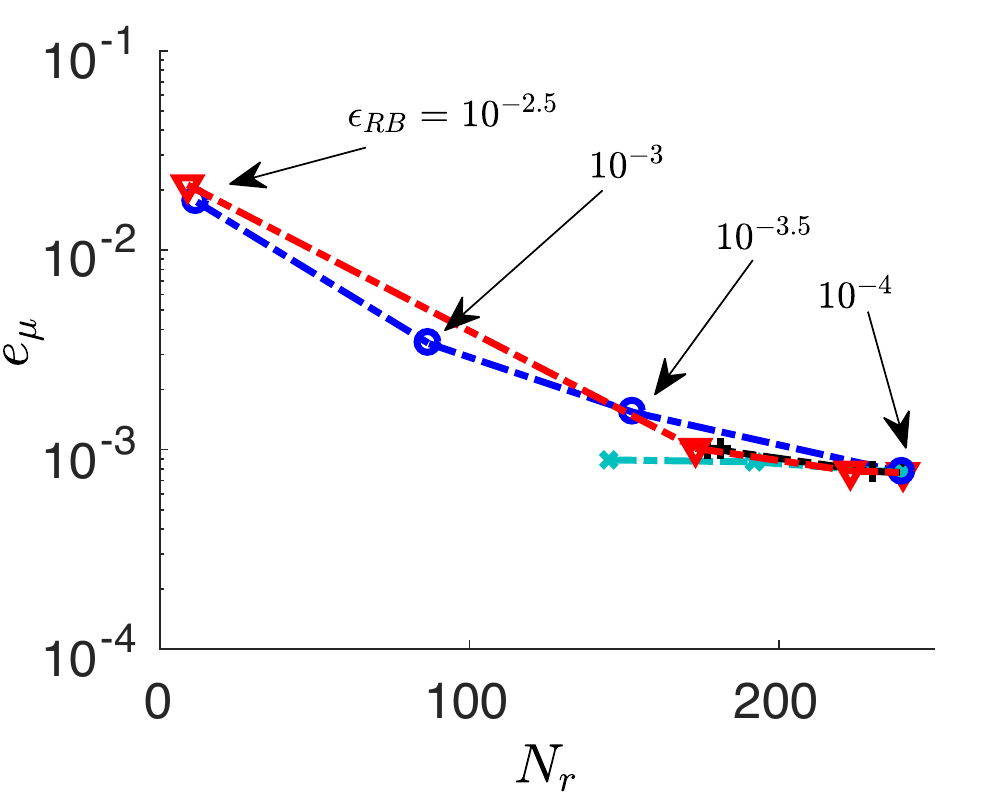}} 
    \centerline{ \rotatebox{90}{\hspace{1.9cm} s.d. }
    \includegraphics[width=4.0cm]{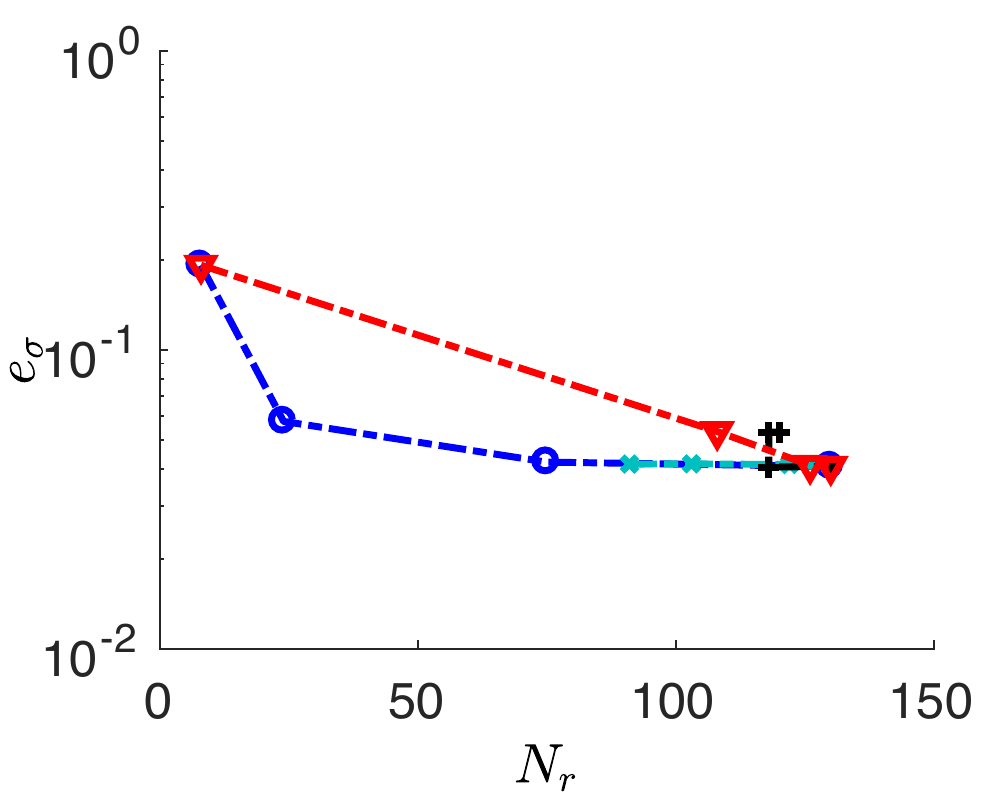}
    \includegraphics[width=4.0cm]{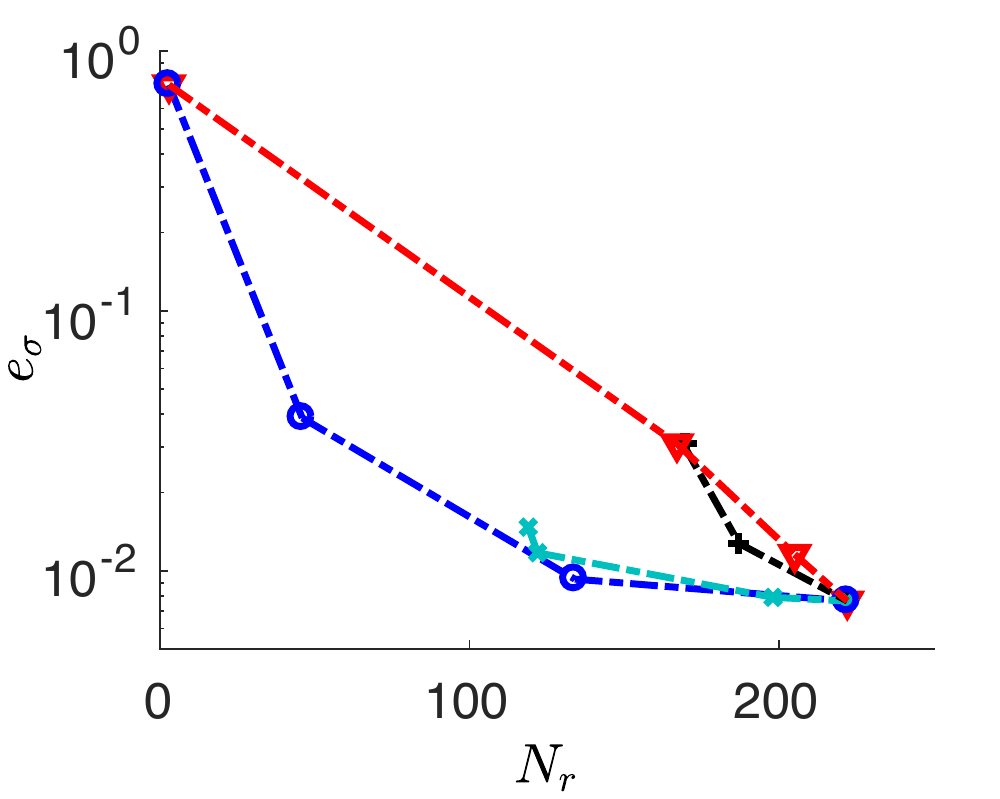}
    \includegraphics[width=4.0cm]{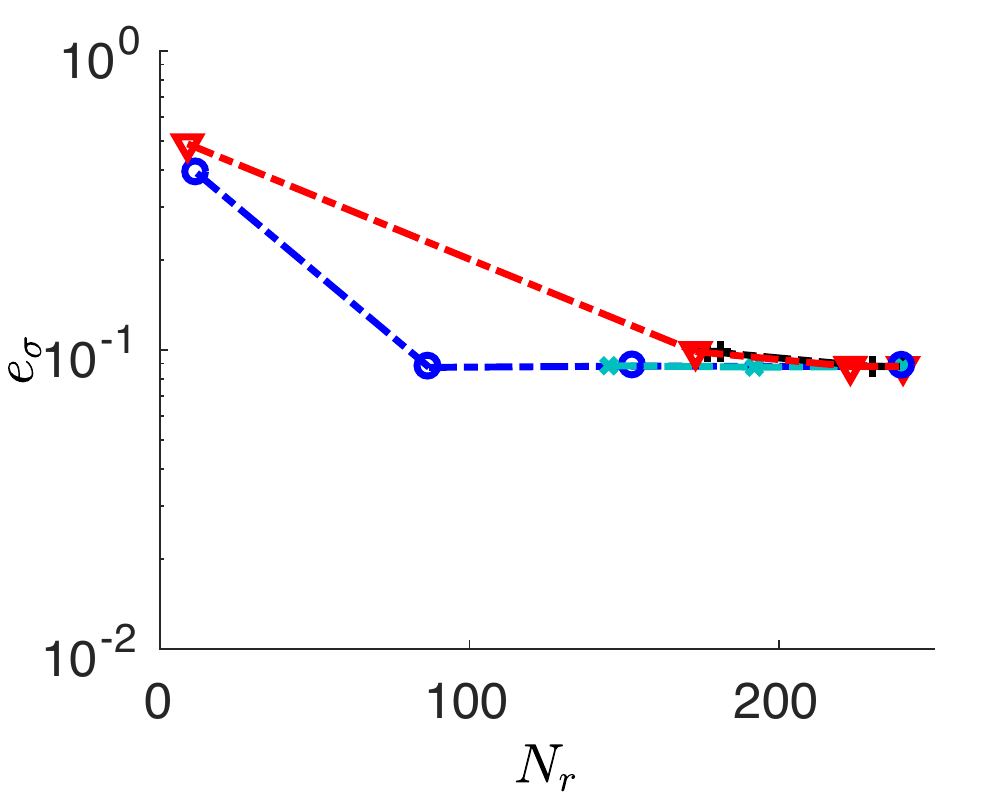}}
    \centerline{
    \includegraphics[width=10cm]{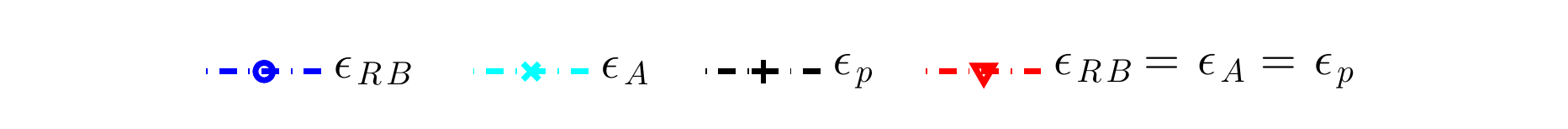}} 
   \caption{Error in the mean and standard deviation versus the size of reduced basis $N_r$ for the case of diffusion parameter $\nu = 1/2$  and $N_D=16$. 
    The bullets correspond to solutions computed by varying each tolerance $\epsilon_{RB}$, $\epsilon_{A}$, and $\epsilon_{p}$ from $10^{-2.5}$ to $10^{-4}$, while others are fixed as $10^{-4}$. The results for $\epsilon_{RB} = \epsilon_{A} = \epsilon_{p}$ are shown as well. The reduced basis solution is more accurate when $\epsilon_{A},\epsilon_{p} < \epsilon_{RB}$. } \label{fig:errRBmmt_dim16_alltol}
\end{figure}


We first study 
the adaptive strategies in Algorithm \ref{Alg:RBM-ANOVA_adaptp} 
for benchmark problems 
on $N_D = 16$ subdomains 
with $1\times 16$, $4\times 4$, and $16\times 1$ partitions. 
The initial parameters are 
selected as $p_0=3$ and $\ell_0=2$, 
and we successively increase $p_K$ by 2.  
We have three different tolerances to investigate, 
the reduced basis acceptance $\epsilon_{RB}$, 
the truncation of the ANOVA decomposition $\epsilon_{A}$, 
and the polynomial order $\epsilon_{p}$.

We first point out in Figure \ref{fig:errRBmmt_dim16_A} 
the advantage of adapting the polynomial order $p$. 
The tolerance $\epsilon_{RB}$ varies from $10^{-2.5}$ to $10^{-4}$, 
while $\epsilon_{p}=10^{-5}$, and 
$\epsilon_{A}=10^{-15}$ to prevent ANOVA truncation. 
The error of the reduced basis solution computed 
with fixed polynomial orders $p=3$ or $5$ 
does not decay after reaching a certain level. 
The accuracy is not improved 
while the size of the reduced basis increases. 
The adaptive procedure for $p$ overcomes this problem.

Next, we study the sensitivity of the method 
with respect to the tolerance. 
One of the tolerances among 
$\epsilon_{RB}$, $\epsilon_{A}$, 
and $\epsilon_{p}$ is 
varied from $10^{-2.5}$ to $10^{-4}$, 
while others are fixed as $10^{-4}$. 
We also test the case where all the tolerances are 
reduced simultaneously as 
$\epsilon_{RB} = \epsilon_{A} = \epsilon_{p}$. 
The results are plotted in Figure \ref{fig:errRBmmt_dim16_alltol} 
which shows that all the errors 
eventually converge to the level attained by 
$\epsilon_{RB} = \epsilon_{A} = \epsilon_{p} = 10^{-4}$. 
%
Moreover, the selected reduced basis solution 
is more accurate 
when $\epsilon_{p}$ is less than $\epsilon_{RB}$.  
In further simulations,  
we take 
$\epsilon_{p} = \epsilon_{RB}/C$, 
for some constant $C$, e.g., $C=2$. 
On the other hand, results for 
$\epsilon_{A} < 10^{-3.5}$ show no difference 
since the effective dimensions $\J_\ell$ include all directions 
(see Figure \ref{fig:cANOVA_dim1x16_all} for the magnitude of 
the ANOVA terms in the $1\times 16$ case). 
Thus, $\epsilon_{A}$ should be chosen carefully 
by taking into consideration the trend of $\gamma_K$ 
to balance accuracy and computational complexity.


\begin{figure}
    \centerline{   $\| \E[u_j] \|$ \hspace{3.7cm}  $p_j$ \hspace{3.7cm}  $N_{r_j}$ \hspace{0.4cm} } 
    \centerline{ 
    \includegraphics[width=4.5cm]{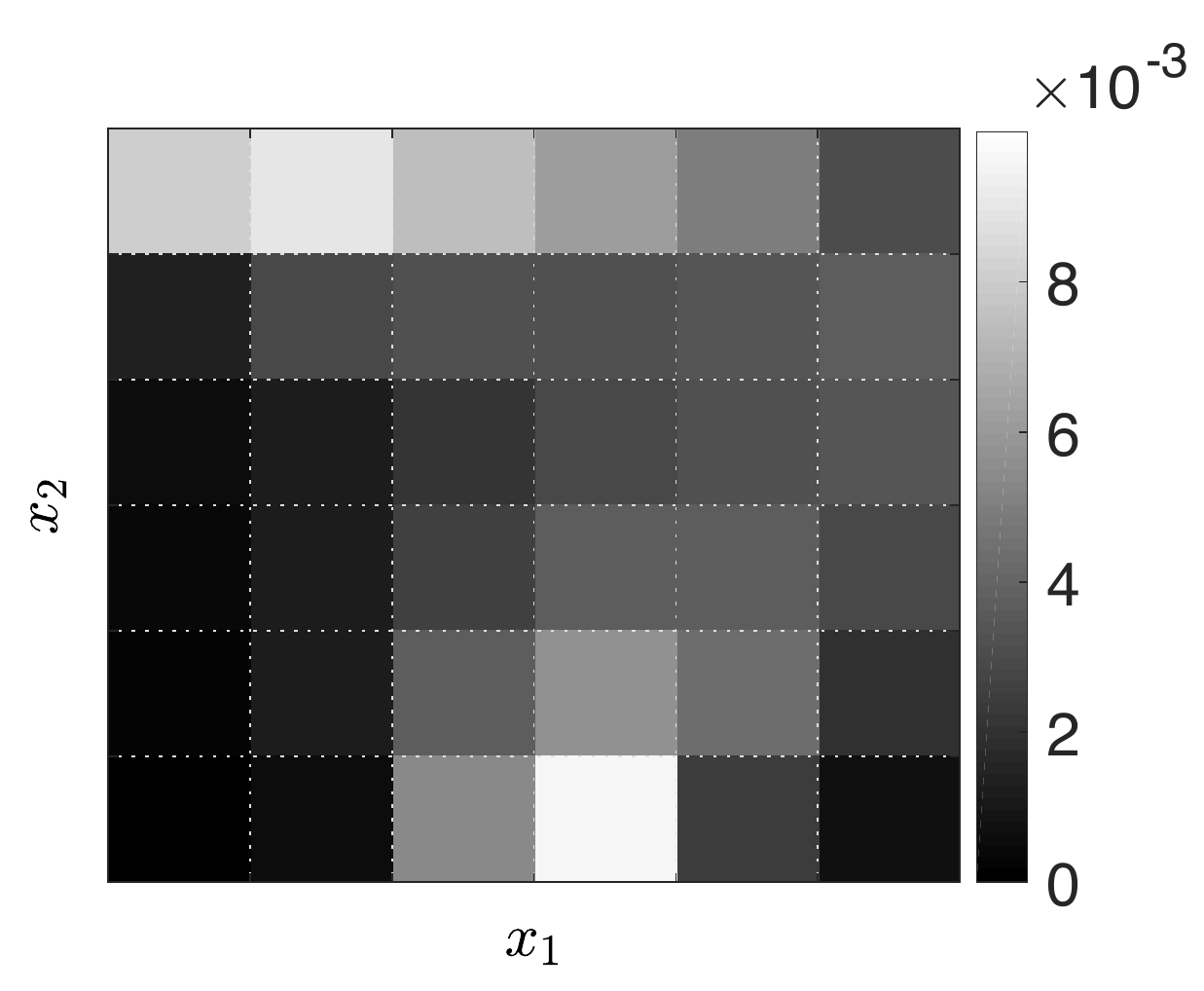}
	\hspace{0.0cm} 
    \includegraphics[width=4.2cm]{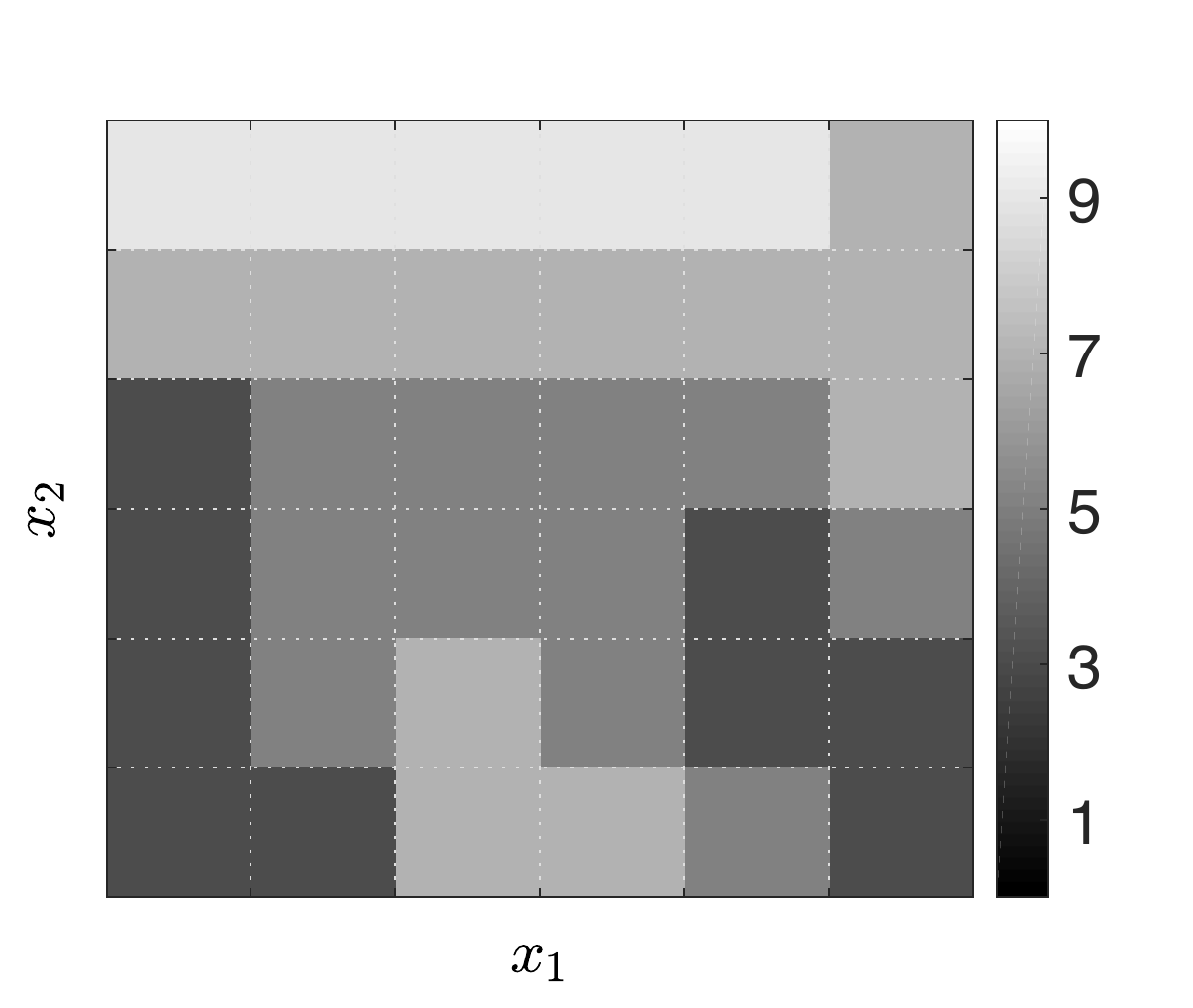}
    \includegraphics[width=4.2cm]{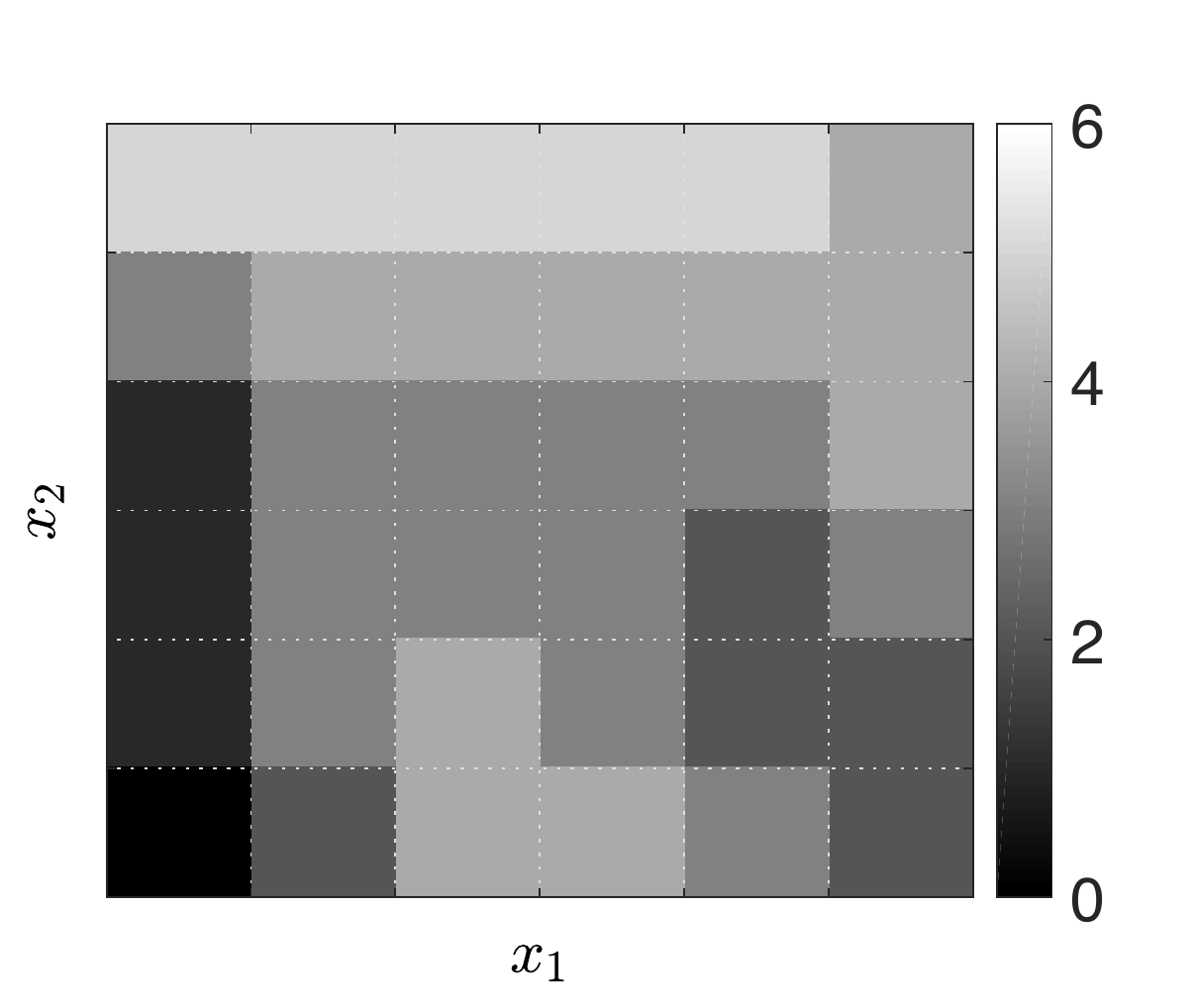}} 
   \caption{
  	The mean of the first order ANOVA expansion term  $\| \E[u_j] \|$, the maximum polynomial order $p_j$, and the number of points selected as the reduced basis $N_{r_j}$ {on each subdomain $D_j$ with $\xi_j$ parameter for $6\times 6$ and $\nu=1/2$. The importance of the subdomains near the top boundary layer and the bottom inflow discontinuity is automatically identified by the algorithm, as is evident in all three figures.} } \label{fig:cgridAdapt_nxn}
\end{figure}

We now investigate the performance of our algorithm 
in higher dimensional parametric space, considering  
$6\times 6$, $8\times 8$, and $10\times 10$ subdomain 
partitions.  
The reduced basis tolerance 
$\epsilon_{RB}$ is varied from $10^{-2.5}$ to $10^{-4}$, 
and $\epsilon_{A}$ and $\epsilon_{p}$ are taken to be proportional 
to $\epsilon_{RB}$,  
$\epsilon_{A} = \epsilon_{p} = \epsilon_{RB}/2$. 
%
In Figure \ref{fig:cgridAdapt_nxn} shows the solution 
on a $6\times 6$ partition with 
diffusion parameter $\nu = 1/2$ using tolerance $\epsilon_{RB} = 10^{-4}$. 
On each subdomain $D_j$ with parameter $\xi_j$, 
the figures plot 
the magnitude of the first order ANOVA terms $\|\E[u_j]\|$, 
the adapted maximum polynomial order $p_j$ in direction $j$,  
and the number of reduced basis $N_{r_j}$ selected in direction $j$. 
They show that 
the algorithm automatically 
increases the polynomial order and 
selects more reduced basis solutions on the 
subdomains near the top boundary layer and 
the bottom inflow discontinuity, 
where the ANOVA term is relatively large. 

\begin{figure}
    \centerline{  {\footnotesize \hspace{0.7cm}  $6 \times 6$ \hspace{3.2cm} $8 \times 8$ \hspace{3.2cm}   $ 10\times 10$ } } 
    \centerline{ \rotatebox{90}{\hspace{1.2cm} {\footnotesize $\nu = 1/2$ } }
    \includegraphics[width=4.0cm]{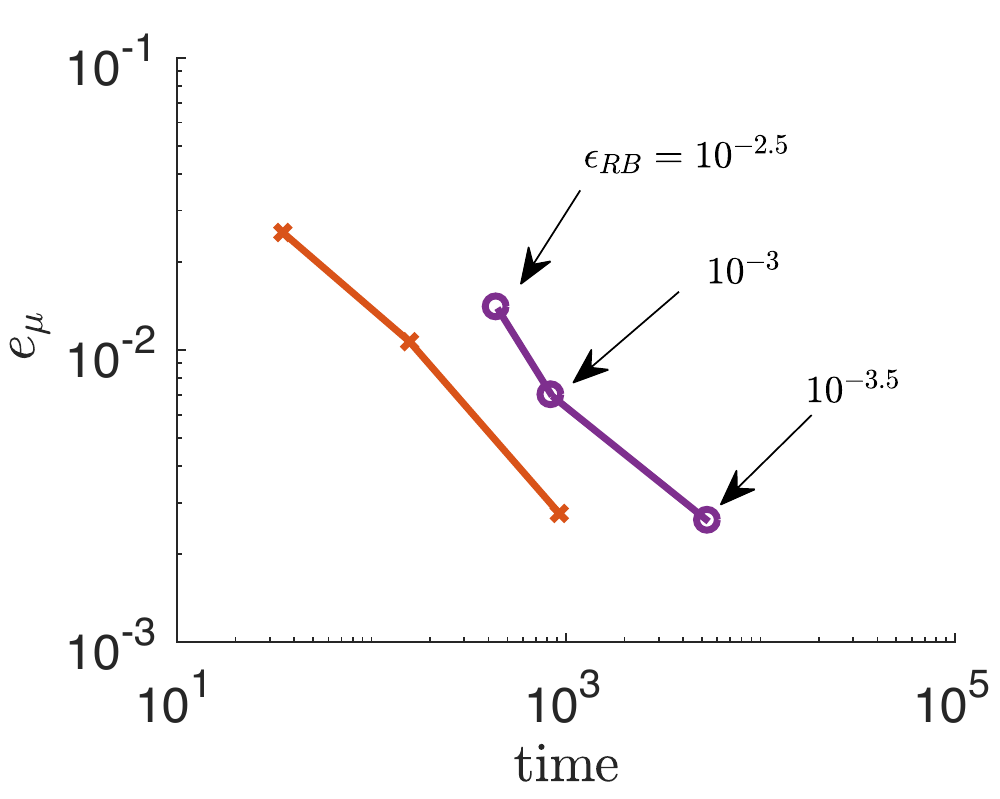}
    \includegraphics[width=4.0cm]{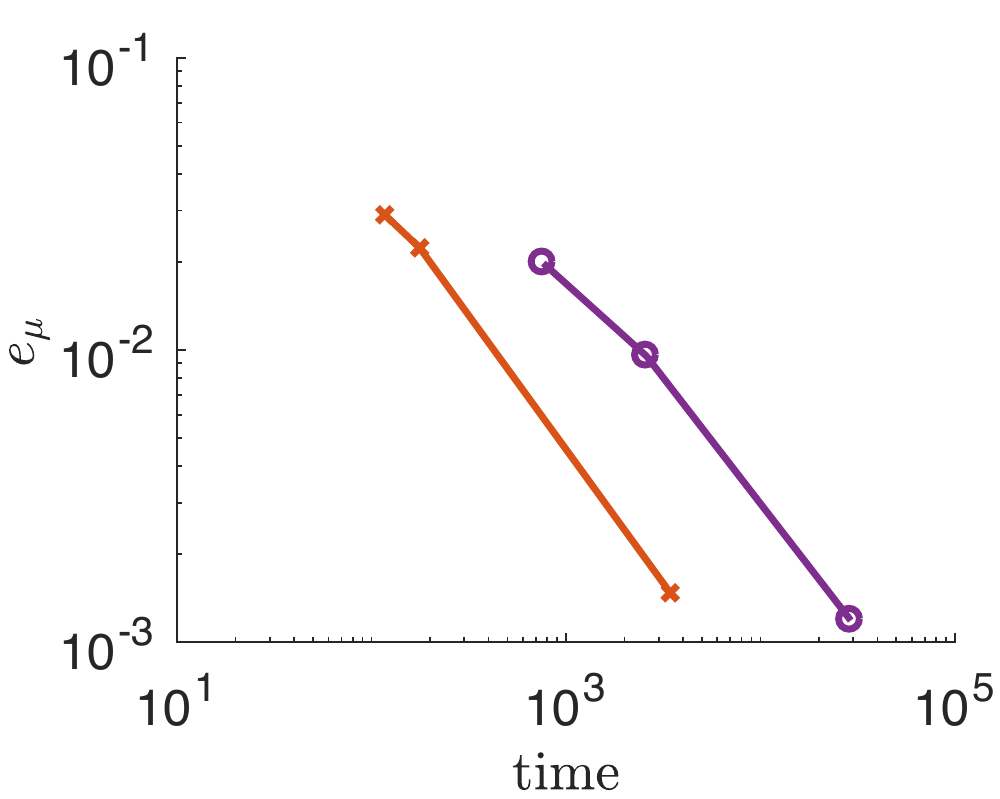}
    \includegraphics[width=4.0cm]{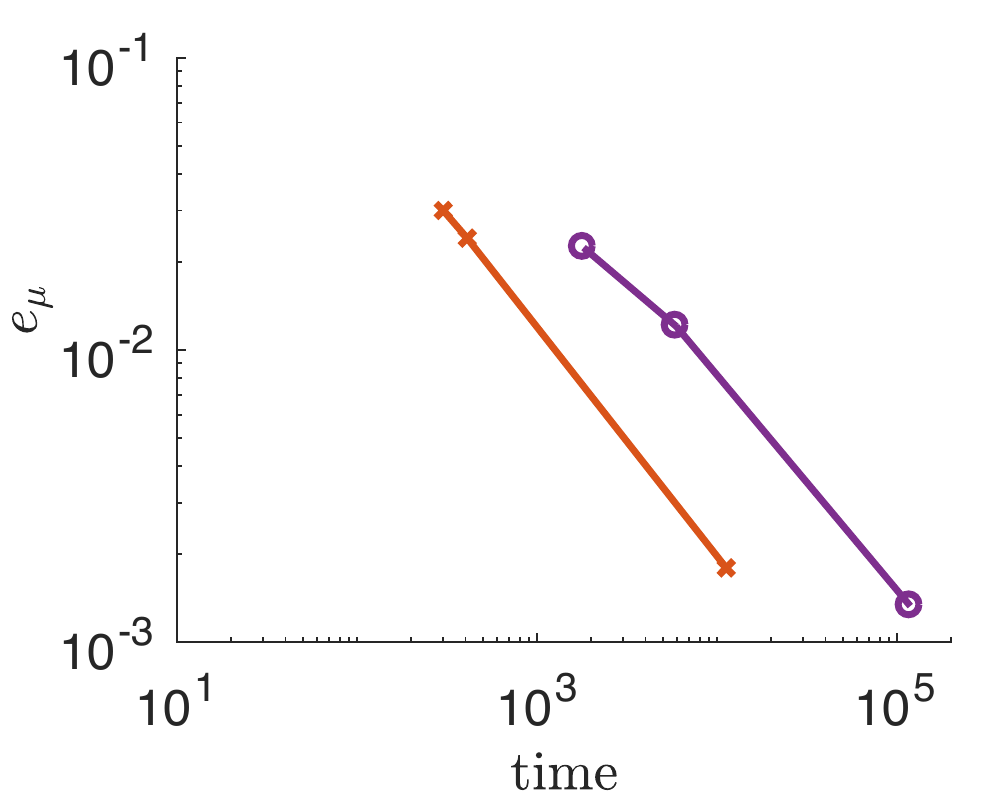}
	}
    \centerline{ \rotatebox{90}{\hspace{1.2cm} {\footnotesize $\nu = 1/20$ } }
    \includegraphics[width=4.0cm]{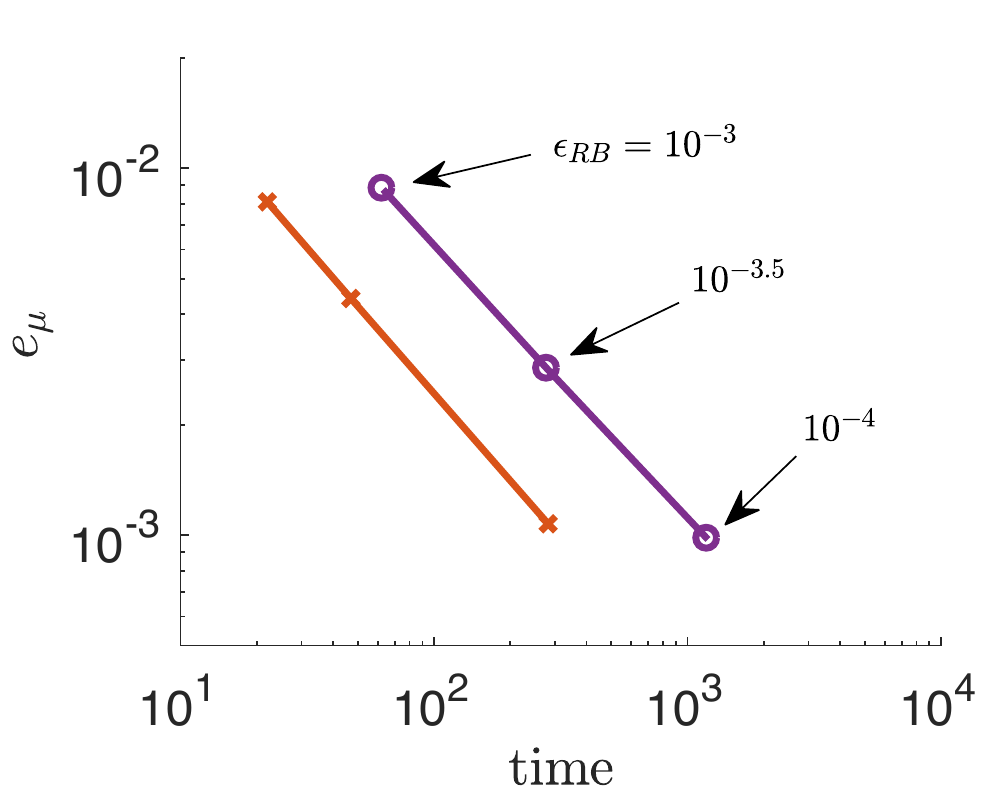}
    \includegraphics[width=4.0cm]{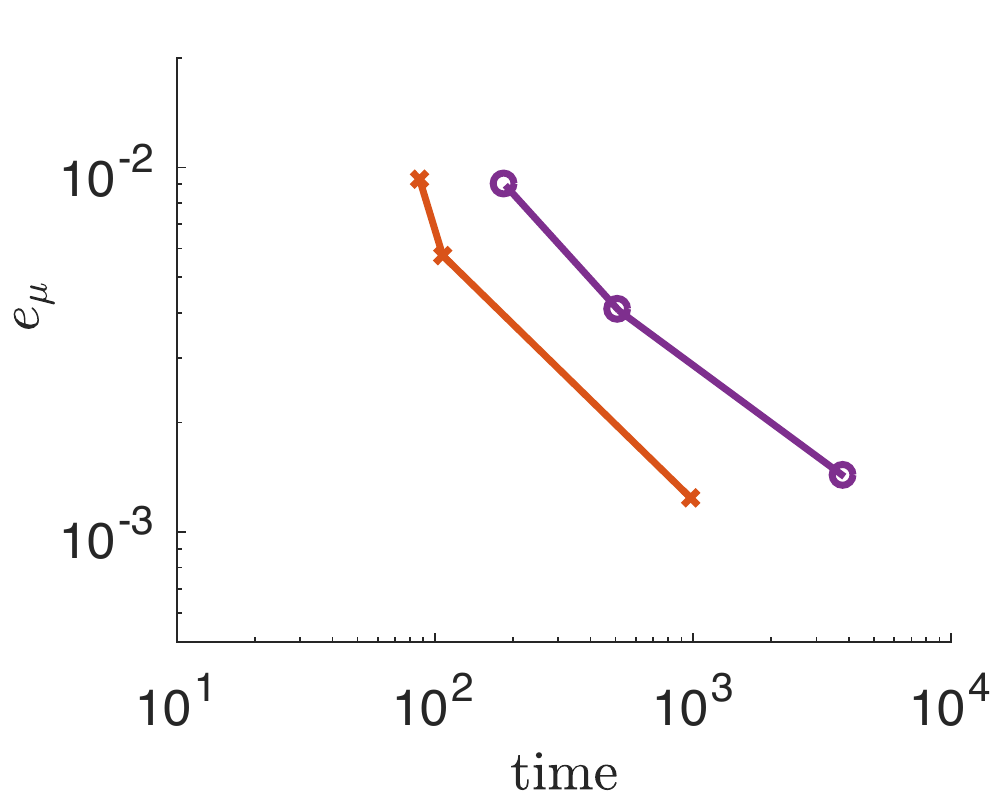}
    \includegraphics[width=4.0cm]{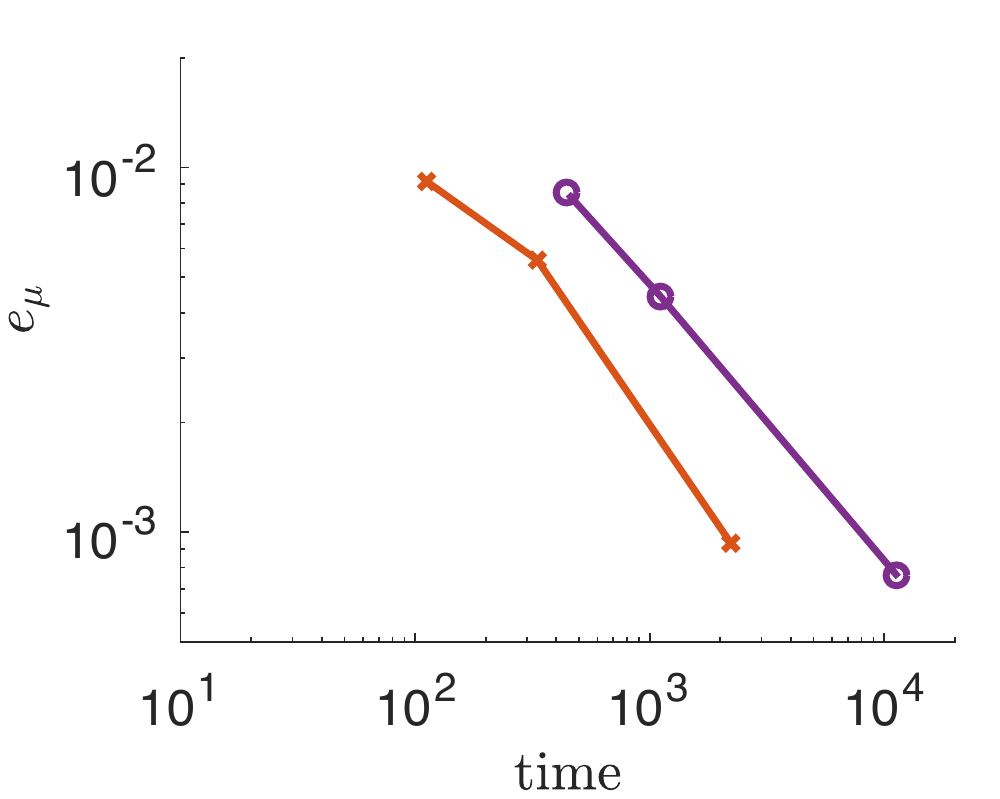}
	}
    \centerline{ 
    \includegraphics[width=9.5cm]{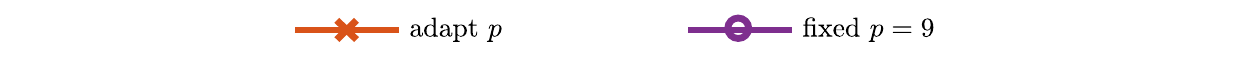}}
   \caption{ Error in the mean and standard deviation versus computational time in seconds. The solution is computed with fixed $p=9$ and adaptive $p$ with $\epsilon_{RB} = 10^{-3}$ to $10^{-4}$.  } 
\label{fig:errvsTime_nxn}
\end{figure}

Let us compare the efficiency of 
the reduced basis algorithm 
using a fixed polynomial order $p=9$ and 
adaptive determination of $p$.  
Table \ref{Tbl:nRBerr_nxn} shows 
the reduced basis size and 
the error in the moments 
for diffusion parameter $\nu=1/20$ (top) 
and $\nu=1/2$ (bottom). 
With a slightly smaller value of 
$\epsilon_{RB}$, 
the adaptive $p$ approach yields 
a reduced basis set 
that achieves the same level of accuracy 
with the fixed $p$ approach. 
%
{ 
Figure \ref{fig:errvsTime_nxn} shows accuracy versus computational 
time for both fixed $p$ and adaptively determined $p$ 
(These results were obtained in Matlab 2016 on a 3.33GHz 
Intel core i5 processor). 
It can be seen that changing $p$ adaptively significantly reduces 
the computational cost.}
%
The computational cost saving 
of the adaptive $p$ approach is primarily due to 
its use of fewer search points 
$C(M,\ell)$ in \eqref{eq:CostGreedy}, 
as shown in Table \ref{Tbl:nRB_nxn_}. 
We compare the total number of points 
on which the reduced basis is searched 
between using the PCM-ANOVA collocation points 
in the entire $M$-dimensions 
with a fixed $p=9$, 
adapting directions with a fixed $p=9$, and 
adapting both directions and $p$. 
Compared to considering a tensor product, 
which yields $p^M$ points, 
the number of PCM-ANOVA points $C(M,\ell)$ is 
reduced by the truncation 
in the interaction level with $\ell=2$. 
However, 
the number is large due to 
the high dimensionality of size $M$, 
and the adaptive approach in 
both dimension and $p$ reduces 
the computational cost. 
For the case of $\nu=1/20$, 
due to its strong anisotropy, 
the cost saving is significant. 
When $\nu=1/2$, 
although the cost saving from 
the truncation of directions is less than for $\nu=1/20$, 
the procedure of adaptive $p$ 
reduces the cost further by increasing $p$ 
only among the effective dimensions.

\begin{table}
\centering
\caption{The size of reduced basis $N_r$ and the error in the moments $e_{\mu}$ and $e_{\sigma}$ 
for  $N_D \times N_D$ partitions and diffusion parameter $\nu=1/2$ (top) and $\nu=1/20$ (bottom), comparing the reduced basis algorithm with adaptive $p$ and fixed $p=9$. Corresponding computational time is plotted in Figure \ref{fig:errvsTime_nxn}.}
\label{Tbl:nRBerr_nxn} 
\begin{tabular}{c|c|ccc|ccc}
\hline 
\multicolumn{2}{c|}{$\nu=1/2$ }  & \multicolumn{3}{c|}{adapt $p$} & \multicolumn{3}{c}{fixed $p=9$}  \\ \hline  
\multicolumn{2}{c|}{$\epsilon_{RB}$}&  $10^{-2.5}$ & $10^{-3}$ & $10^{-3.5}$ &  $10^{-2.5}$ & $10^{-2.9}$ & $10^{-3.4}$ \\ \hline 
\multirow{3}{*}{ {$6\times 6\,\,\,$}} & 
$N_r$ &3&34&157  &  3&38&191  \\ 
 & $e_{\mu}$ &  3.156e-02 & 1.067e-02  & 2.756e-03 
 &   3.507e-02 &  1.391e-02 &  2.700e-03  \\ 
 & $e_{\sigma}$ & 7.634e-01 &  2.632e-01  & 5.282e-02 
 &   7.244e-01 &  2.785e-01 &  5.440e-02  \\ \hline 
\multirow{3}{*}{ {$8\times 8\,\,$}} & 
$N_r$ &3&26&197 &  3&37&291  \\ 
 & $e_{\mu}$ &  3.381e-02 & 2.234e-02 &  1.474e-03  
 &   3.626e-02 & 1.980e-02  & 1.377e-03   \\ 
 & $e_{\sigma}$ & 7.767e-01 & 5.135e-01 &  5.480e-02  
 &  7.411e-01   & 3.996e-01 &  5.487e-02  \\ \hline 
\multirow{3}{*}{ {$10\times 10\,$}} & 
$N_r$ &3&35&248 &  3&43&385 \\ 
 & $e_{\mu}$ &   3.578e-02 &  2.411e-02  & 1.793e-03  
 &   3.751e-02 &   2.239e-02 &  1.680e-03  \\ 
 & $e_{\sigma}$ & 7.858e-01  & 4.973e-01 &  6.340e-02 
 &    7.533e-01 &  4.140e-01 &  6.272e-02  \\ \hline 
\multicolumn{8}{c}{ }  \\ 
\end{tabular} 
\begin{tabular}{c|c|ccc|ccc}
\hline 
\multicolumn{2}{c|}{$\nu=1/20$}  & \multicolumn{3}{c|}{adapt $p$} & \multicolumn{3}{c}{fixed $p=9$}  \\ \hline  
\multicolumn{2}{c|}{$\epsilon_{RB}$}&  $10^{-3}$ & $10^{-3.5}$ & $10^{-4}$ &  $10^{-3}$ & $10^{-3.4}$ & $10^{-3.9}$ \\ \hline 
 \multirow{3}{*}{ {$6\times 6\,\,\,$}} & 
$N_r$ &  13&24&85  	&  13&29&89  \\ 
  & $e_{\mu}$ &  8.095e-03 & 4.415e-03 &  1.073e-03 
  &   8.749e-03 &  4.897e-03 &  1.086e-03 \\ 
  & $e_{\sigma}$ & 6.279e-01 &  2.424e-01 &  2.305e-02 
  &   6.335e-01 &  2.448e-01 &  2.515e-02 \\ \hline
\multirow{3}{*}{ {$8\times 8\,\,$}} & 
$N_r$ &15&30&118 &  17&31&119  \\ 
  & $e_{\mu}$ &  9.287e-03 & 5.736e-03 & 1.242e-03 
  &   8.942e-03 & 5.899e-03 & 1.348e-03 \\ 
  & $e_{\sigma}$ & 6.777e-01 & 2.670e-01 & 3.429e-02 
  &   6.195e-01  & 2.676e-01 & 3.582e-02 \\ \hline
 \multirow{3}{*}{ {$10\times 10\,$}} & 
$N_r$ &16&34&137 &  21&35&153 \\ 
  & $e_{\mu}$ &   9.159e-03 &  5.574e-03 &  9.333e-04  
  &   8.447e-03 &  6.093e-03 &  9.052e-04 \\ 
  & $e_{\sigma}$ & 7.037e-01  & 3.208e-01 &  2.363e-02 
  &   6.196e-01 &  3.237e-01 &  2.448e-02 \\ \hline 
\end{tabular} 
\end{table} 

\begin{table}
\centering
\caption{Comparison of the number of search points $C(M,\ell)$ \eqref{eq:ANOVAcost} using the entire directions in $M$ dimensions with a fixed polynomial order $p=9$, adapting direction with a fixed $p=9$, and adapting both direction and $p$. }
\begin{tabular}{c|c|c|c|c|c|c|c}
\multicolumn{2}{c}{ $ $ } & \multicolumn{3}{|c|}{$8\times 8$} & \multicolumn{3}{c}{$10\times 10$} \\ \hline  
\multirow{2}{*}{$\nu$} & \multirow{2}{*}{$\epsilon_{RB}$}
 & $M=64$ &  \multicolumn{2}{c|}{adapt $M$}  & $M=100$ & \multicolumn{2}{c}{adapt $M$}  \\ \cline{3-8}
 & & $p=9$ & $p=9$  &  adapt $p$  & $p=9$ & $p=9$ &  adapt $p$ \\ \hline  
1/2 & $10^{-3}$ & \multirow{4}{*}{
$\begin{matrix}  C(64,2)\\=129536	   \end{matrix}$
} & 24704 & 542  &  \multirow{4}{*}{
$\begin{matrix}  C(100,2)\\=317600 \end{matrix}$
} & 38880 &  512 \\ 
1/2 & $10^{-3.5}$ &  &  99072 & 13924  & & 223904 & 25410  \\ 
1/20 & $10^{-3}$&  & 2304 & 122  & & 3680 &  156  \\  
1/20 &  $10^{-3.5}$&  & 6336 &  246  &  & 10592 & 288   \\  
\hline  
\end{tabular} 
\label{Tbl:nRB_nxn_} 
\end{table}


\section{Conclusion} 

In this paper, we developed an adaptive reduced basis method 
for anisotropic parametrized stochastic PDEs 
based on the PCM and ANOVA approximation. 
The proposed algorithm 
enhances the reduced basis collocation methods 
by first obtaining the effective dimensions and 
then increasing the resolution in those directions. 
The dimensional adaptive procedure 
based on the ANOVA decomposition 
reduces the computational cost significantly 
for an anisotropic problem. 
In addition, 
the ANOVA indicator rearranges both of the search directions 
and the reduced basis solution, 
so that 
the reduced basis is sorted 
in a descending order of error 
with lower computational cost 
than the cost that 
typically occurs in greedy type approaches. 
%
The second adaptive procedure increases 
the resolution 
in each direction among the effective dimensions. 
This eliminates the stagnation of the error 
that may occur by fixing the polynomial order in PCM. 
Moreover,  
{this adaptive procedure 
can be terminated when} 
the ANOVA term converges. 
The proposed adaptive reduced basis algorithm is implemented 
on the stationary stochastic convection-diffusion problem. 
The algorithm effectively detects 
the directions in the stochastic space 
that are associated with the irregularity of the solution, 
for instance, discontinuous boundary conditions 
and steep boundary layers. 
By emphasizing those directions 
and by increasing the number of 
collocation points until necessary, 
the algorithm reduces the computational cost 
while improving accuracy in the moments.

Our future work includes analysis of 
the proposed adaptive reduced basis algorithm 
using reduced basis error indicators 
based on posterior error estimates of parametrized PDEs 
\cite{RBM_Patera1} and combining them 
with analytical studies of 
adaptive PCM \cite{XWGK_JCP05,Foo} and 
anchored ANOVA method \cite{Handy_SISC_ANOVA_2012}. 
In addition, 
we aim to further investigate 
the effect of the quadrature rules, 
e.g., Clenshaw-Curtis and Gauss abscissae, 
as candidate points 
to the accuracy of the reduced basis solution. 
Since different choice for the interpolation 
and anchor points in PCM-ANOVA 
yields significantly distinct set of collocation points, 
strategies to utilize the error estimates involving different 
interpolation points \cite{CCGL_Lloyd,CompareCCGL_1D} and 
to choose the anchor point that minimizes the error 
estimates of the ANOVA approximation 
\cite{Handy_ANOVA_2010,Handy_SISC_ANOVA_2012} 
may be useful. 
Finally, 
extending our algorithm to 
nonlinear and time-dependent problems, for example, 
stochastic Navier-Stokes equation, is another future direction. 
Various techniques including 
linearized iterative solvers such as Newton iterations 
\cite{RBM_time,FEM_Elman} 
and empirical interpolation methods \cite{RBM_YMaday1,VirginiaRBM}
can be employed.


\acknowledgements 
This work was supported by the U. S. Department of Energy Office of Advanced Scientific Computing Research, Applied Mathematics program under Award Number DE-SC0009301, and by the U. S. National Science Foundation under grant DMS1418754. 

\bibliographystyle{acm}
\bibliography{RBM_mePCMANOVA.bib}

\end{document}